\documentclass[12pt]{article}
\usepackage[pagebackref]{hyperref}

\usepackage[utf8]{inputenc}
\usepackage[T1]{fontenc}

\usepackage[]{amsfonts}
\usepackage{amssymb}
\usepackage{amsmath,amsthm}
\def\couleur(#1 #2 #3)
	{
\special{ps::[end]}}

\def\bx#1{\setbox1=\hbox{\kern3pt{#1}\kern3pt}			
 \dimen1=\ht1 \advance\dimen1 by 3pt \dimen2=\dp1 \advance\dimen2 by 3pt
 \setbox1=\hbox{\vrule height\dimen1 depth\dimen2\box1\vrule}%
 \setbox1=\vbox{\hrule\box1\hrule}%
 \advance\dimen1 by .4pt \ht1=\dimen1
 \advance\dimen2 by .4pt \dp1=\dimen2 \box1\relax}

\def\wbb#1{\kern#1em}
\def\vci{\vrule  width.02em height1.47ex depth-.0ex}		
\def\11{{\rm\wbb{.2}\vci\wbb{-.37}1}}

\def\underset#1#2{\mathrel{\mathop{\kern0pt #2}\limits_{#1}}}

\def\overset#1#2{\mathrel{\mathop{\kern0pt #2}\limits^{#1}}}

\def\Supp{\mathop{\rm Supp}\nolimits}

\newtheorem{thm}{Theorem}[section]
\newtheorem{lem}[thm]{Lemma}
\newtheorem{prop}[thm]{Proposition}
\newtheorem{cor}[thm]{Corollary}
\newtheorem{defin}[thm]{Definition}
\newtheorem{rem}[thm]{Remark}

\begin{document}

\title{On the $\displaystyle L^{r}$ Hodge theory in complete non compact riemannian manifolds.}

\author{Eric Amar}

\date{ }
\maketitle
 \ \par 
\ \par 
\renewcommand{\abstractname}{Abstract}

\begin{abstract}
We study solutions for the Hodge laplace equation $\Delta u=\omega
 $ on $p$ forms with $\displaystyle L^{r}$ estimates for $\displaystyle
 r>1.$ Our main hypothesis is that $\Delta $ has a spectral gap
 in $\displaystyle L^{2}.$ We use this to get \emph{non classical}
 $\displaystyle L^{r}$ Hodge decomposition theorems. An interesting
 feature is that to prove these decompositions we never use the
 boundedness of the Riesz transforms in $\displaystyle L^{s}.$\ \par 
These results are based on a generalisation of the Raising Steps
 Method to complete non compact riemannian manifolds.\ \par 
\end{abstract}

\tableofcontents
\ \par 

\section{Introduction.}
\quad In the sequel, a riemannian manifold $\displaystyle (M,g)$ means
 a ${\mathcal{C}}^{\infty }$ smooth connected riemannian manifold
 of dimension $\displaystyle \geq 3.$\ \par 
\quad In this work we study the problem of $\displaystyle L^{r}$ estimates
 of the Laplace equation $\displaystyle \Delta u=\omega $ for
 the Hodge laplacian on $p$-forms and the Hodge decomposition
 theorems on complete non compact riemannian manifolds.\ \par 
\quad This problem was studied by a several authors, in particular
 N. Lohou\'e in~\cite{Lohoue16} (see also the references therein).
 Also the problem of Hodge decompositions has a long history
 and for the recent developments one can see the papers by X.D.
 Li~\cite{XDLi09}, \cite{XDLi10}, \cite{XDLi010} and also related
 to several complex variables~\cite{XDLi2010} (see also the references
 therein).\ \par 
\quad In all those works the boundedness of the Riesz transforms are
 explicitely used and in this work, where the Hodge decompositions
 are \emph{not the classical ones}, we shall see that it is not the case.\ \par 
\quad Let me describe the method we shall use.\ \par 
\quad Suppose you are interested by solving an equation $\displaystyle
 Du=\omega ,$ in a manifold $M$ with estimates of type Lebesgue
 $\displaystyle L^{r}$ or Sobolev $\displaystyle W^{d,r}$ ; you
 know how to solve it globally with "threshold" estimates $\displaystyle
 L^{s}\rightarrow L^{s}$ and locally with estimates $\displaystyle
 L^{r}\rightarrow L^{t}$ with a strict increase of the regularity,
 for instance $\displaystyle \frac{1}{t}=\frac{1}{r}-\delta ,\
 \delta >0$ for any $\displaystyle r\leq s,$ then the Raising
 Steps Method (RSM for short) gives a \emph{global} solution
 $v$ of $\displaystyle Dv=\omega $ which is essentially in $\displaystyle
 L^{t}(M)$ for $\omega \in L^{r}(M).$\ \par 
\quad I introduced this method in~\cite{AmarSt13} to get solutions
 for the $\bar \partial $ equation with good estimates in relatively
 compact domains in Stein manifold. I extend it to linear partial
 differential operator $D$ of any finite order $m$ in~\cite{HodgeCompact15}
 and I apply it to study the Poisson equation for the Hodge laplacian
 on forms in spaces $\displaystyle L^{r}(M)$ where $\displaystyle
 (M,g)$ is a compact riemannian manifold. This gave $\displaystyle
 L^{r}$ Hodge decomposition theorems as was done by C. Scott~\cite{Scott95},
 but by an entirely different approach.\ \par 
\ \par 
\quad The aim of this work is to extend it to the case of complete
 non compact riemannian manifold, and, as we shall see, at no
 point we shall use the boundedness of the Riesz transforms.\ \par 

\subsection{Solutions of the Poisson equation for the Hodge laplacian.}
\quad Let $\displaystyle (M,g)$ be a ${\mathcal{C}}^{\infty }$ smooth
 connected riemannian manifold with metric tensor $g$ and $\displaystyle
 n=\mathrm{d}\mathrm{i}\mathrm{m}M\geq 3\ ;$ let $d$ be the exterior
 derivative, $\displaystyle d^{*}$ its formal adjoint with respect
 to the Riemannian volume measure $\displaystyle dv_{g}={\sqrt{\mathrm{d}\mathrm{e}\mathrm{t}g}}dx,$
 where $\displaystyle dx$ is the Lebesgue measure in the chart
 $\displaystyle x,$ and $\Delta =\Delta _{p}:=dd^{*}+d^{*}d$
 the Hodge laplacian acting on $p$ forms. Let $\displaystyle
 L^{r}_{p}(M)$ be the space of $p$ forms on $M$ in the Lebesgue
 space $\displaystyle L^{r}(M).$\ \par 
\quad We get the local solution of the Hodge Laplacian $\Delta u=\omega
 $ in a ball $\displaystyle B(x,R)$ in $\displaystyle (M,g)$
 with a radius $\displaystyle R(x)$ small enough to make this
 ball "not too different"  to a ball in the euclidean space ${\mathbb{R}}^{n}\
 ;$ this "admissible" radius is a special case of the "harmonic
 radius" of Hebey and Herzlich~\cite{HebeyHerzlich97}. If $\omega
 $ is a $p$ form in $L^{r}(B(x,R))$ then we get a $p$ form solution
 $u$ in the Sobolev space $\displaystyle W^{2,r}(B(x,r))$ of
 the ball, hence in $\displaystyle L^{t}(B(x,R))$ with $\displaystyle
 \frac{1}{t}=\frac{1}{r}-\frac{2}{n}$ by the Sobolev embeddings.
       This is done classically by use of the Newtonian potential.
 So the first assumption for the RSM is true : we have locally
 a strict increase of the regularity.\ \par 
\quad In order to get global solutions we need to cover the manifold
 $M$ with our "admissible balls" and for this we use a classical
 "Vitali type covering" with a uniformly finite overlap. We shall
 denote it by ${\mathcal{C}}.$\ \par 
\ \par 
\quad When comparing non compact $M$ to the compact case treated in~\cite{HodgeCompact15},
 we have two  important issues :\ \par 
\quad (i) the "admissible" radius may go to $0$ at infinity, which
 is the case, for instance, if the canonical volume measure $\displaystyle
 dv_{g}$ of $\displaystyle (M,g)$ is finite and $M$ is not compact ;\ \par 
\quad (ii) if $\displaystyle dv_{g}$ is not finite, which is the case,
 for instance, if the "admissible"  radius is bounded below,
 then $p$ forms in $\displaystyle L^{t}_{p}(M)$ are generally
 not in $\displaystyle L^{r}_{p}(M)$ for $\displaystyle r<t.$\ \par 
\quad We address these problems by use of adapted weights on $\displaystyle
 (M,g).$ These weights are relative to the covering ${\mathcal{C}}$
 : they are positive functions which vary slowly on the balls
 of the covering ${\mathcal{C}}.$\ \par 
\quad To deal with the problem (i) we shall use a weight\ \par 
\quad \quad \quad \begin{equation}  w_{0}(x)=R(x)^{-2k}\label{HC30}\end{equation}\ \par 
for an adapted integer $k,$ where $\displaystyle R(x)$ is the
 admissible  radius at the point $\displaystyle x\in M.$\ \par 
\quad To deal with the problem (ii) we shall use a weight $\alpha (x)$
 which is in $\displaystyle L^{\mu }(M)$ with $\displaystyle
 \mu :=\frac{2t}{2-t},$ for a $\displaystyle t<2,$ i.e.\ \par 
\quad \quad \quad \begin{equation}  \gamma (w,t):=\int_{M}{w^{\frac{2t}{2-t}}}dv_{g}<\infty
 .\label{HC31}\end{equation}\ \par 
This is done to get $L^{2}_{p}(M)\subset L^{t}_{p}(M,\alpha ).$\ \par 
\quad Our Hodge decompositions are \emph{not the classical ones} because
 we do not use the laplacian adapted to those weights, but we
  always use the standard laplacian.\ \par 
\quad We define the Sobolev spaces $\displaystyle W^{d,r}_{p}(M)$ of
 $\displaystyle (M,g)$ following E. Hebey~\cite{Hebey96}, and we set\ \par 

\begin{defin}
~\label{HC33}We shall define the Sobolev exponents $\displaystyle
 S_{k}(r)$ by $\displaystyle \frac{1}{S_{k}(r)}:=\frac{1}{r}-\frac{k}{n}.$
\end{defin}
Then our first result is a "twisted" Calderon Zygmund inequalities
 (CZI) with weight, different from results in~\cite{GuneysuPigola}
 because we have weights and our forms are \emph{not} asked to
 have compact support.\ \par 

\begin{thm}
Let $\displaystyle (M,g)$ be a complete riemannian manifold.
 Let $w$ be a weight relative to the ${\mathcal{C}}_{\epsilon
 }$ associated covering $\displaystyle \lbrace B(x_{j},5r(x_{j}))\rbrace
 _{j\in {\mathbb{N}}}$ and set $\displaystyle w_{0}:=R(x)^{-2}.$
 Let $\displaystyle u\in L^{r}_{p}(M,ww_{0}^{r})$ such that $\Delta
 u\in L^{r}_{p}(M,w)$ ; then there are constants $\displaystyle
 C_{1},C_{2}$ depending only on $n=\mathrm{d}\mathrm{i}\mathrm{m}_{{\mathbb{R}}}M,\
 r$ and $\epsilon $ such that:\par 
\quad \quad \quad $\displaystyle \ {\left\Vert{u}\right\Vert}_{W^{2,r}(M,w)}\leq
 C_{1}{\left\Vert{u}\right\Vert}_{L^{r}(M,ww_{0}^{r})}+C_{2}{\left\Vert{\Delta
 u}\right\Vert}_{L^{r}(M,w)}.$\par 
Moreover we have for $\displaystyle t=S_{2}(r)$ that $\displaystyle
 u\in L^{t}_{p}(M,w^{t})$ with $\displaystyle \ {\left\Vert{u}\right\Vert}_{L^{t}(M,w^{t})}\leq
 c{\left\Vert{u}\right\Vert}_{W^{2,r}(M,w^{r}w_{0}^{t})}.$
\end{thm}
\quad We set, for a weight $\alpha ,\ $${\mathcal{H}}^{r}_{p}(M,\alpha
 ):=L^{r}_{p}(M,\alpha )\cap \mathrm{k}\mathrm{e}\mathrm{r}\Delta
 _{p},$ the space of harmonic $p$ forms in $\displaystyle L^{r}(M,\alpha
 ).$\ \par 
\quad This is our main hypothesis :\ \par 
(HL2,p) $\Delta =\Delta _{p}$ has a spectral gap in $\displaystyle
 L^{2}_{p}(M),$ i.e. there is no spectrum of $\Delta _{p}$ in
 an open interval $\displaystyle (0,\eta )$ with $\displaystyle \eta >0.$\ \par 
This assumption allows us to use $\displaystyle L^{2}_{p}(M)$
 as a threshold for the Raising Steps Method.\ \par 
\quad The (HL2,p) assumption is known to be true in the case of the
 hyperbolic manifold ${\mathbb{H}}^{2n}$ of dimension $\displaystyle
 2n$ for any value of $p\in \lbrace 0,\ 2n\rbrace .$ For $\displaystyle
 p\neq n$ the space ${\mathcal{H}}_{p}^{2}$ is reduced to $\displaystyle
 0.$ For ${\mathbb{H}}^{2n+1}$ the (HL2,p) is valid for $\displaystyle
 p\neq n$ and $\displaystyle p\neq n+1$ and, out of these two
 cases, the space ${\mathcal{H}}_{p}^{2}$ is reduced to $\displaystyle
 0$ as was proved by Donnelly~\cite{Donnelly81}.\ \par 
\quad When $\displaystyle \mathrm{R}\mathrm{i}\mathrm{c}(M)\geq -c^{2}$
 and $M$ is open at infinity then $\displaystyle 0\notin \mathrm{S}\mathrm{p}\Delta
 _{\mathrm{0}}$ by a result of Buser, see Lott~\cite{Lott96},
 proposition 6, p. 353, hence (HL2,0) is true.  If $M$ is a normal
 covering of a compact manifold $X$ with covering group $\Gamma
 ,$ then $\displaystyle 0\notin \mathrm{S}\mathrm{p}\Delta _{\mathrm{0}}$
 iff $\Gamma $ is not amenable by a result of Brooks, see Lott~\cite{Lott96},
 corollary 3, p. 354, for precise references. Hence (HL2,0) is
 true if $\Gamma $ is not amenable.\ \par 
\quad For $r=2,$ there is the orthogonal projection $H$ from $\displaystyle
 L^{2}_{p}(M)$ on ${\mathcal{H}}^{2}_{p}(M)\ ;$ we shall prove
 that this projection extends to $\displaystyle L^{r}(M,w_{0}^{r}),$
 with $\displaystyle w_{0}:=R(x)^{-2k}$ and $\displaystyle R(x)$
 the admissible radius at $\displaystyle x\in M,$ as in~(\ref{HC30}),
 i.e.\ \par 
\quad \quad \quad \begin{equation}  \forall r\leq 2,\ H\ :\ L^{r}(M,w_{0}^{r})\rightarrow
 {\mathcal{H}}^{2}_{p}(M)\label{HC32}\end{equation}\ \par 
boundedly and we get the following results on solutions of the
 Poisson equation.\ \par 

\begin{thm}
~\label{CL28} Suppose that $\displaystyle (M,g)$ is a complete
 riemannian manifold ; let $\displaystyle r<2$ and choose a weight
 $\alpha \in L^{\infty }(M)$ verifying $\gamma (\alpha ,r)<\infty
 .$ Set $\displaystyle t:=\min (2,S_{2}(r)).$ If $\displaystyle
 t<2,$ take the weight $\alpha \in L^{\infty }(M)$ verifying
 also $\displaystyle \gamma (\alpha ,t)<\infty .$ Suppose we
 have conditions (HL2,p).\par 
\quad Take $k$ big enough so that the threshold $\displaystyle S_{k}(r)\geq
 2,$ and set $\displaystyle w_{0}(x):=R(x)^{-2k},$ then for any
 $\omega \in L^{r}_{p}(M,w_{0}^{r})$ verifying $\displaystyle
 H\omega =0,$ for the orthogonal projection $H$ defined in corollary~\ref{CF9},
 there is a $\displaystyle u\in W^{2,r}_{p}(M,\alpha )\cap L^{t}_{p}(M,\alpha
 ),$ such that $\Delta u=\omega .$\par 
Moreover the solution $u$ is given linearly with respect to $\omega .$
\end{thm}
\quad Here $k$ was chosen such that $\displaystyle S_{k}(r)\geq 2$
 in order to use $\displaystyle L^{2}_{p}(M)$ as a threshold
 for the Raising Steps Method.\ \par 
\quad Setting $\displaystyle r'$ for the conjugate exponent for $\displaystyle
 r,\ \frac{1}{r'}+\frac{1}{r}=1,$ by duality from theorem~\ref{CL28},
 we get\ \par 

\begin{thm}
Suppose that $\displaystyle (M,g)$ is a complete riemannian manifold
 ; suppose we have $\displaystyle r<2$ and (HL2,p), then with
 $\displaystyle k::S_{k}(r)\geq 2,$ and $\displaystyle w_{0}(x):=R(x)^{-k},$
 for any $\varphi \in L^{2}_{p}(M)\cap L^{r'}_{p}(M),\ H\varphi
 =0,$ there is a $\displaystyle u\in L^{r'}(M,w_{0}^{r})$ such
 that $\displaystyle \Delta u=\varphi .$ This solution is linear
 with respect to $\varphi .$\par 
If we add the hypothesis that the $\epsilon _{0}$ admissible
 radius is bounded below, we get\par 
\quad \quad \quad $\displaystyle u:=(T-C)^{*}\varphi ,\ u\in W^{2,r'}_{p}(M)$ and
 $u$ verifies $\displaystyle \Delta u=\varphi .$
\end{thm}
\quad By theorem 1.3 in Hebey~\cite{Hebey96}, we have that the harmonic
 radius $\displaystyle r_{H}(1+\epsilon ,\ 2,0)$ is bounded below
 if the Ricci curvature $\displaystyle Rc$ verifies $\displaystyle
 {\left\Vert{\nabla Rc}\right\Vert}_{\infty }<\infty $  and the
 injectivity radius is bounded below. This implies that the $\epsilon
 $ admissible radius is also bounded below.\ \par 

\subsection{Hodge decomposition in $\displaystyle L^{r}$ spaces.
 Known results.}
\quad In 1949, Kodaira~\cite{Kodaira49} proved that the $\displaystyle
 L^{2}$-space of $p$-forms on $\displaystyle (M,g)$ has the orthogonal
 decomposition :\ \par 
\quad \quad \quad $L^{2}_{p}(M)={\mathcal{H}}^{2}_{p}\oplus {\overline{d{\mathcal{D}}_{p-1}(M)}}\oplus
 {\overline{d^{*}{\mathcal{D}}_{p+1}(M)}},$\ \par 
and in 1991 Gromov~\cite{Gromov91} proved a strong $\displaystyle
 L^{2}$ Hodge decomposition, under the hypothesis (HL2,p) :\ \par 
\quad \quad \quad $\displaystyle L^{2}_{p}(M)={\mathcal{H}}^{2}_{p}\oplus dW^{1,2}_{p-1}(M)\oplus
 d^{*}W^{1,2}_{p+1}(M).$\ \par 
\quad In 1995  Scott~\cite{Scott95} proved a strong $\displaystyle
 L^{r}$ Hodge decomposition but on \emph{compact} riemannian manifold\ \par 
\quad \quad \quad $\displaystyle \forall r>1,\ L^{r}_{p}(M)={\mathcal{H}}^{r}_{p}\oplus
 dW^{1,r}_{p-1}(M)\oplus d^{*}W^{1,r}_{p+1}(M).$\ \par 
\quad Let $\displaystyle d_{\varphi }^{*}$ be the formal adjoint of
 $d$ relatively to the measure $\displaystyle d\mu (x)=e^{-\varphi
 (x)}dv_{g}(x),$ where $\displaystyle \varphi \in {\mathcal{C}}^{2}(M),$
 and let $\displaystyle \Delta _{\varphi ,p}:=dd_{\varphi }^{*}+d_{\varphi
 }^{*}d$ acting on $p$ forms. Setting $\Delta =\mathrm{T}\mathrm{r}\nabla
 ^{2}$ the covariant Laplace Beltrami operator acting on $p$
 forms and $\displaystyle L=\Delta -\nabla \varphi \cdot \nabla
 ,$ then, in 2009 X-D. Li~\cite{XDLi09} proved, among other nice
 results, a strong $\displaystyle L^{r}$ Hodge decomposition
 on complete non compact riemannian manifold :\ \par 

\begin{thm}
(X-D. Li) Let $\displaystyle r>1,\ r'=\frac{r}{r-1}.$ Let $\displaystyle
 (M,g)$ be a complete riemannian manifold, $\varphi \in {\mathcal{C}}^{2}(M),$
 and $\displaystyle d\mu (x)=e^{-\varphi (x)}dv_{g}(x).$ Suppose
 that the Riesz transforms $\displaystyle d\Delta _{\varphi ,p}^{-1/2}$
 and $\displaystyle d^{*}\Delta _{\varphi ,p}^{-1/2}$ are bounded
 in $\displaystyle L^{r}$ and $\displaystyle L^{r'},$ and the
 Riesz potential is bounded in $\displaystyle L^{r}.$ Suppose
 also that $\displaystyle (M,g)$ is $L$ stochastically complete,
 then the strong $\displaystyle L^{r}$ Hodge direct sum decomposition
 holds on $p$ forms :\par 
\quad \quad \quad $L^{r}_{p}(M,\mu )={\mathcal{H}}^{r}_{p}(M,\mu )\oplus dW^{1,r}_{p-1}(M,\mu
 )\oplus d_{\varphi }^{*}W^{1,r}_{p+1}(M,\mu ).$
\end{thm}
\quad These results are valid for the family of weights $\displaystyle
 \varphi \in {\mathcal{C}}^{2}(M)$ and for the Hodge laplacian
 associated to them, in the Witten sense~\cite{Witten82}. Nevertheless
 it is worthwhile to notice that, even in the classical case
 $\varphi \equiv 0,$ this result was new at the time it was proved,
 2007, by X-D. Li.\ \par 

\subsection{Non classical Hodge decomposition in $\displaystyle
 L^{r}$ spaces. Main results.}
\quad The results of X-D. Li are based on the boundedness of the Riesz
 transforms in $\displaystyle L^{r}$ and $\displaystyle L^{r'}$
 and the results we get use mainly the spectral gap hypothesis
 (HL2,p). X-D. Li was already concerned by the fact that the
 bottom of the spectrum of $\Delta $ should be strictly positive
 ; the difference here is that we allow an eigenvalue $0$ but
 a gap without spectrum after it, which gives the possible existence
 of non trivial harmonic functions in $\displaystyle L^{2}.$
 This is the meaning of (HL2,p).\ \par 
\quad In this way our results may appear to be the natural generalisation
 of Gromov results from $\displaystyle L^{2}$ to $\displaystyle
 L^{r}.$ On the other hand our results are proved only in the
 case $\varphi =0.$\ \par 
\quad Our decompositions are \emph{non classical} because we use weights
 to get estimates, but we use the usual laplacian, \emph{not
 the Witten laplacian adapted to these weights.}\ \par 
\quad We shall need the following definition.\ \par 

\begin{defin}
Let $\alpha $ be a weight on $M,$ we define the space $\displaystyle
 \tilde W^{2,r}_{p}(M,\alpha )$ to be\par 
\quad \quad \quad $\displaystyle \tilde W^{2,r}_{p}(M,\alpha ):=\lbrace u\in L^{r}_{p}(M,\alpha
 )::\Delta u\in L^{r}_{p}(M,\alpha )\rbrace $\par 
with the norm\par 
\quad \quad \quad $\displaystyle \ {\left\Vert{u}\right\Vert}_{\tilde W^{2,r}_{p}(M,\alpha
 )}:={\left\Vert{u}\right\Vert}_{L^{r}_{p}(M,\alpha )}+{\left\Vert{\Delta
 u}\right\Vert}_{L^{r}_{p}(M,\alpha )}.$
\end{defin}
\quad To get these decomposition theorems we shall apply our results
 on solutions of the Poisson equation.\ \par 

\begin{thm}
Let $\displaystyle (M,g)$ be a complete riemannian manifold.
 Let $\displaystyle r<2$ and take a weight $\alpha \in L^{\infty
 }(M)$ be such that $\gamma (\alpha ,r)<\infty $ ; with $\displaystyle
 k::S_{k}(r)\geq 2,$ set $\displaystyle w_{0}=R(x)^{-2k},$ and
 suppose we have hypothesis (HL2,p). We have the direct decomposition
 given by linear operators :\par 
\quad \quad \quad $\ L^{r}_{p}(M,w_{0}^{r})={\mathcal{H}}_{p}^{2}\oplus \Delta
 (W^{2,r}_{p}(M,\alpha )).$\par 
With $\displaystyle r'>2,$ the conjugate exponent to $\displaystyle
 r,$ we have the weaker decomposition, still given by linear operators :\par 
\quad \quad \quad $L^{r'}_{p}(M)\cap L^{2}_{p}(M)={\mathcal{H}}_{p}^{2}\cap {\mathcal{H}}^{r'}_{p}+\Delta
 (\tilde W^{2,r'}_{p}(M)).$
\end{thm}
\quad Because $\displaystyle H:L^{r}(M,w_{0}^{r})\rightarrow {\mathcal{H}}^{2}_{p}(M)$
 boundedly by~(\ref{HC32}), where $H$ is the orthogonal projection
 from $\displaystyle L^{2}_{p}(M)$ on $\displaystyle {\mathcal{H}}^{2}_{p}(M),$
 this explain the appearance of $\displaystyle L^{2}_{p}(M)$
 and $\displaystyle {\mathcal{H}}^{2}_{p}(M)$ in the second part
 of the previous theorem.\ \par 
\quad To replace $\displaystyle \tilde W^{2,r'}_{p}(M))$ by $\displaystyle
 W^{2,r'}_{p}(M,\alpha ))$ the price is the hypothesis that the
 $\epsilon _{0}$ admissible radius is bounded below. So we get\ \par 

\begin{cor}
Suppose the admissible radius is bounded below and suppose also
 hypothesis (HL2,p). Take $\displaystyle r'>2,$ then we have
 the direct decomposition given by linear operators\par 
\quad \quad \quad $\displaystyle L^{r'}_{p}(M)\cap L^{2}_{p}(M)={\mathcal{H}}_{p}^{2}\cap
 {\mathcal{H}}^{r'}_{p}\oplus \Delta (W^{2,r'}_{p}(M)).$
\end{cor}
\quad As a corollary we get\ \par 

\begin{cor}
~\label{HC29}Let $\displaystyle r<2$ and choose a weight $\alpha
 \in L^{\infty }(M)$ such that $\gamma (\alpha ,r)<\infty $ ;
 with $\displaystyle k::S_{k}(r)\geq 2,$ set $\displaystyle w_{0}=R(x)^{-2k},$
 and suppose we have hypothesis (HL2,p). We have the direct decompositions
 given by linear operators\par 
\quad \quad \quad $\displaystyle L^{r}_{p}(M,w_{0}^{r})={\mathcal{H}}_{p}^{2}\oplus
 d(W^{1,r}_{p}(M,\alpha ))\oplus d^{*}(W^{1,r}_{p}(M,\alpha )).$\par 
With $\displaystyle r'>2$ the conjugate exponent of $\displaystyle
 r,$ and adding the hypothesis that the $\epsilon _{0}$ admissible
 radius is bounded below, we get\par 
\quad \quad \quad $\displaystyle L^{r'}_{p}(M)\cap L^{2}_{p}(M)={\mathcal{H}}_{p}^{2}\cap
 {\mathcal{H}}^{r'}_{p}\oplus d(W^{1,r'}_{p}(M))\oplus d^{*}(W^{1,r'}_{p}(M)).$
\end{cor}
\quad We also have weak $\displaystyle L^{r}$ Hodge decompositions,
 where $\displaystyle d^{*}$ is the adjoint of $d$ with respect
 to the \emph{usual volume measure,} not the weighted one, despite
 the weight appearing here.\ \par 
\quad We shall need another hypothesis :\ \par 
(HWr) if the space ${\mathcal{D}}_{p}(M)$ is dense in $\displaystyle
 W^{2,r}_{p}(M).$\ \par 
We already know that (HWr) is true if :\ \par 
\quad $\bullet $ either : the injectivity radius is strictly positive
 and the Ricci curvature is bounded~[\cite{Hebey96} theorem 2.8, p. 12].\ \par 
\quad $\bullet $ or : $M$ is geodesically complete with a bounded curvature
 tensor~[\cite{GuneysuPigola} theorem 1.1 p.3].\ \par 

\begin{thm}
Suppose that $\displaystyle (M,g)$ is a complete riemannian manifold,
 fix $\displaystyle r<2$ and choose a bounded weight $\alpha
 $ with $\displaystyle \gamma (\alpha ,r)<\infty .$\par 
Take $k$ with $\displaystyle S_{k}(r)\geq 2,$ and set the weight
 $\displaystyle w_{0}:=R(x)^{-2k}.$ Suppose we have (HL2,p) and
 (HW2). Then    $L^{r}_{p}(M,\alpha )={\mathcal{H}}_{p}^{r}(M,\alpha
 )\oplus {\overline{\Delta ({\mathcal{D}}_{p}(M))}},$ the closure
 being taken in $L^{r}(M,\alpha ).$
\end{thm}
\quad We also have a weak $\displaystyle L^{r}$ Hodge decomposition
 without hypothesis (HWr):\ \par 

\begin{thm}
Suppose that $\displaystyle (M,g)$ is a complete riemannian manifold
 and suppose we have (HL2,p). Fix $\displaystyle r<2$ and take
 a weight $\alpha $ verifying $\gamma (\alpha ,r)<\infty .$ Then we have:\par 
\quad \quad \quad $L^{r}_{p}(M,\alpha )={\mathcal{H}}_{p}^{r}(M,\alpha )\oplus
 {\overline{d({\mathcal{D}}_{p-1}(M))}}\oplus {\overline{d^{*}({\mathcal{D}}_{p+1}(M))}},$\par
 
the closures being taken in $L^{r}(M,\alpha ).$
\end{thm}
\quad For the case $\displaystyle r>2$ we need a stronger hypothesis,
 namely that the $\epsilon _{0}$ admissible radius is bounded
 below. Then we get a \emph{classical weak Hodge decompositions.}\ \par 

\begin{thm}
Suppose that $\displaystyle (M,g)$ is a complete riemannian manifold
 and suppose the $\epsilon _{0}$ admissible radius is bounded
 below and (HWr) and suppose also hypothesis (HL2,p). Fix $\displaystyle
 r>2,$ then we have\par 
\quad \quad \quad $\displaystyle L^{r}_{p}(M)={\mathcal{H}}_{p}^{r}(M)\oplus {\overline{\Delta
 ({\mathcal{D}}_{p}(M))}}.$\par 
Without (HWr) we still get\par 
\quad \quad \quad $L^{r}_{p}(M)={\mathcal{H}}_{p}^{r}(M)\oplus {\overline{d({\mathcal{D}}_{p-1}(M))}}\oplus
 {\overline{d^{*}({\mathcal{D}}_{p+1}(M))}}.$\par 
All the closures being taken in $L^{r}(M).$
\end{thm}

\begin{rem}
By theorem 1.3 in Hebey~\cite{Hebey96}, we have that the harmonic
 radius $\displaystyle r_{H}(1+\epsilon ,\ 2,0)$ is bounded below
 if the Ricci curvature $\displaystyle Rc$ verifies $\displaystyle
 \ {\left\Vert{\nabla Rc}\right\Vert}_{\infty }<\infty $  and
 the injectivity radius is bounded below. This implies that the
 $\epsilon $ admissible radius is also bounded below.\par 
Moreover if we add the hypothesis that the Ricci curvature $\displaystyle
 Rc$ is bounded below then by Proposition 2.10 in Hebey~\cite{Hebey96},
 we have hypothesis (HWr).
\end{rem}
\quad These results are based on the raising steps method :\ \par 

\begin{thm}
(Raising Steps Method) Let $\displaystyle (M,g)$ be a riemannian
 manifold and take $w$ a weight relative to the Vitali covering
 $\displaystyle \lbrace B(x_{j},5r(x_{j}))\rbrace _{j\in {\mathbb{N}}}.$\par 
For any $\displaystyle r\leq 2,$ any threshold $\displaystyle
 s\geq r,$ take $\displaystyle k\in {\mathbb{N}}$ such that $\displaystyle
 t_{k}:=S_{k}(r)\geq s$ then, with $\displaystyle w_{0}(x):=w(x)R(x)^{-2k},$\par
 
\quad \quad \quad $\displaystyle \forall \omega \in L^{r}_{p}(M,w_{0}^{r}),\ \exists
 v\in L^{r}_{p}(M,w^{r})\cap L^{s_{1}}_{p}(M,w^{s_{1}})\cap W^{2,r}(M,w^{r}),\
 \exists \tilde \omega \in L^{s}_{p}(M,w^{s})::\Delta v=\omega
 +\tilde \omega $\par 
with $\displaystyle s_{1}=S_{2}(r)$ and we have the control of the norms :\par 
\quad \quad \quad $\displaystyle \forall q\in \lbrack r,s_{1}\rbrack ,\ {\left\Vert{v}\right\Vert}_{L^{q}_{p}(M,w^{q})}\leq
 C_{q}{\left\Vert{\omega }\right\Vert}_{L^{r}_{p}(M,w_{0}^{r})}\
 ;\ {\left\Vert{v}\right\Vert}_{W^{2,r}_{p}(M,w^{r})}\leq C_{r}{\left\Vert{\omega
 }\right\Vert}_{L^{r}_{p}(M,w_{0}^{r})}\ ;$\par 
\quad \quad \quad \quad \quad \quad \quad \quad \quad \quad \quad \quad $\displaystyle \ {\left\Vert{\tilde \omega }\right\Vert}_{L^{s}_{p}(M,w^{s})}\leq
 C_{s}{\left\Vert{\omega }\right\Vert}_{L^{r}_{p}(M,w_{0}^{r})}.$\par 
Moreover $v$ and $\tilde \omega $ are linear in $\omega .$\par 
If $M$ is complete and $\omega $ is of compact support, so are
 $v$ and $\tilde \omega .$
\end{thm}
\quad I thank the referee for his pertinent questions and remarks making
 precise the meaning of these non classical Hodge decompositions.\ \par 
\ \par 
\quad This work will be presented in the following way.\ \par 
\quad In section 2 we define the admissible balls, the admissible 
 radius and the basic facts relative to them.\ \par 
\quad In section 3 we use a Vitali type covering lemma with our admissible
  balls and we prove that its overlap is finite.\ \par 
\quad In section 4 we define the Sobolev spaces, following E. Hebey~\cite{Hebey96}.\
 \par 
\quad In section 5 we prove the local estimates for the Hodge Laplacian.
 This is essentially standard by use of classical results from
 Gilbarg and Trudinger~\cite{GuilbargTrudinger98}.\ \par 
\quad In section 6 we develop the Raising Steps Method in the non compact
 case. The useful weights are defined here.\ \par 
This is the basis of our results.\ \par 
\quad In section 7 we prove Calderon Zygmund inequalities with weights.\ \par 
\quad In section 8 we deduce the applications to the Poisson equation
 associated to the Hodge Laplacian.\ \par 
\quad In section 9 we use these solutions to get non classical strong
 $\displaystyle L^{r}$ Hodge decomposition theorems. We also
 get non classical weak $\displaystyle L^{r}$ Hodge decomposition
 theorems.\ \par 

\section{Basic facts.}

\begin{defin}
Let $(M,g)$ be a riemannian manifold and $\displaystyle x\in
 M.$ We shall say that the geodesic ball $\displaystyle B(x,R)$
 is $\epsilon $ {\bf admissible} if there is a chart $\displaystyle
 \varphi \ :\ (x_{1},...,x_{n})$ defined on it with\par 
\quad 1) $\displaystyle (1-\epsilon )\delta _{ij}\leq g_{ij}\leq (1+\epsilon
 )\delta _{ij}$ in $\displaystyle B(x,R)$ as bilinear forms,\par 
\quad 2) $\displaystyle \ \sum_{\left\vert{\beta }\right\vert =1}{\sup
 \ _{i,j=1,...,n,\ y\in B_{x}(R)}\left\vert{\partial ^{\beta
 }g_{ij}(y)}\right\vert }\leq \epsilon .$
\end{defin}
\ \par 

\begin{defin}
~\label{CL26}Let $\displaystyle x\in M,$ we set $\displaystyle
 R'(x)=\sup \ \lbrace R>0::B(x,R)\ is\ \epsilon \ admissible\rbrace
 .$ We shall say that $\displaystyle R_{\epsilon }(x):=\min \
 (1,R'(x))$ is the $\epsilon $ {\bf admissible radius} at $\displaystyle x.$
\end{defin}
\quad Our admissible radius is smaller than the harmonic radius $\displaystyle
 r_{H}(1+\epsilon ,\ 1,\ 0)$ defined in Hebey~[\cite{Hebey96}, p. 4].\ \par 
\quad By theorem 1.3 in Hebey~\cite{Hebey96}, we have that the harmonic
 radius $\displaystyle r_{H}(1+\epsilon ,\ 2,0)$ is bounded below
 if the Ricci curvature $\displaystyle Rc$ verifies $\displaystyle
 {\left\Vert{\nabla Rc}\right\Vert}_{\infty }<\infty $  and the
 injectivity radius is bounded below. This implies easily that
 the $\epsilon $ admissible radius is also bounded below.\ \par 

\begin{rem}
By its very definition, we always have $\displaystyle R_{\epsilon }(x)\leq 1.$
\end{rem}
\quad Of course, without any extra hypotheses on the riemannian manifold
 $M,$ we have $\displaystyle \forall \epsilon >0,\ \forall x\in
 M,$ taking $\displaystyle g_{ij}(x)=\delta _{ij}$ in a chart
 on $\displaystyle B(x,R)$ and the radius $R$ small enough, the
 ball $\displaystyle B(x,R)$ is $\epsilon $ admissible.\ \par 
We shall use the following lemma.\ \par 

\begin{lem}
~\label{CF10}Let $(M,g)$ be a riemannian manifold then with $\displaystyle
 R(x)=R_{\epsilon }(x)=$ the $\epsilon $ admissible radius at
 $\displaystyle x\in M$ and $\displaystyle d(x,y)$ the riemannian
 distance on $\displaystyle (M,g)$ we get :\par 
\quad \quad \quad $\displaystyle d(x,y)\leq \frac{1}{4}(R(x)+R(y))\Rightarrow R(x)\leq 4R(y).$
\end{lem}
\quad Proof.\ \par 
Let $\displaystyle x,y\in M::d(x,y)\leq \frac{1}{4}(R(x)+R(y))$
 and suppose for instance that $\displaystyle R(x)\geq R(y).$
 Then $\displaystyle y\in B(x,R(x)/2)$ hence we have $\displaystyle
 B(y,R(x)/4)\subset B(x,\frac{3}{4}R(x)).$ But by the definition
 of $\displaystyle R(x),$ the ball $\displaystyle B(x,\frac{3}{4}R(x))$
 is admissible  and this implies that the ball $\displaystyle
 B(y,R(x)/4)$ is also admissible  for exactly the same constants
 and the same chart ; this implies that $\displaystyle R(y)\geq
 R(x)/4.$ $\blacksquare $\ \par 

\section{Vitali covering.~\label{CL24}}

\begin{lem}
~\label{CF0}Let ${\mathcal{F}}$ be a collection of balls $\displaystyle
 \lbrace B(x,r(x))\rbrace $ in a metric space, with $\forall
 B(x,r(x))\in {\mathcal{F}},\ 0<r(x)\leq R.$ There exists a disjoint
 subcollection ${\mathcal{G}}$ of ${\mathcal{F}}$ with the following
 property :\par 
\quad \quad every ball $B$ in ${\mathcal{F}}$ intersects a ball $C$ in ${\mathcal{G}}$
 and $\displaystyle B\subset 5C.$
\end{lem}
This is a well known lemma, see for instance~\cite{EvGar92},
 section 1.5.1].\ \par 
\ \par 
\quad So fix $\epsilon >0$ and let $\displaystyle \forall x\in M,\
 r(x):=R_{\epsilon }(x)/120,\ $where $\displaystyle R_{\epsilon
 }(x)$ is the admissible  radius at $\displaystyle x,$ we built
 a Vitali covering with the collection ${\mathcal{F}}:=\lbrace
 B(x,r(x))\rbrace _{x\in M}.$ So lemma~\ref{CF0} gives a disjoint
 subcollection ${\mathcal{G}}$ such that every ball $B$ in ${\mathcal{F}}$
 intersects a ball $C$ in ${\mathcal{G}}$ and we have $\displaystyle
 B\subset 5C.$ We set ${\mathcal{G}}':=\lbrace x_{j}\in M::B(x_{j},r(x_{j}))\in
 {\mathcal{G}}\rbrace $ and ${\mathcal{C}}_{\epsilon }:=\lbrace
 B(x,5r(x)),\ x\in {\mathcal{G}}'\rbrace $ : we shall call ${\mathcal{C}}_{\epsilon
 }$ the $\epsilon $ {\bf admissible  covering} of $\displaystyle (M,g).$\ \par 
\quad Then we have :\ \par 

\begin{prop}
~\label{CF2}Let $(M,g)$ be a riemannian manifold, then the overlap
 of the $\epsilon $ admissible covering ${\mathcal{C}}_{\epsilon
 }$ is less than $\displaystyle T=\frac{(1+\epsilon )^{n/2}}{(1-\epsilon
 )^{n/2}}(120)^{n},$ i.e.\par 
\quad \quad \quad $\forall x\in M,\ x\in B(y,5r(y))$ where $B(y,r(y))\in {\mathcal{G}}$
 for at most $T$ such balls.\par 
So we have\par 
\quad \quad \quad $\forall f\in L^{1}(M),\ \sum_{j\in {\mathbb{N}}}{\int_{B_{j}}{\left\vert{f(x)}\right\vert
 dv_{g}(x)}}\leq T{\left\Vert{f}\right\Vert}_{L^{1}(M)}.$
\end{prop}
\quad Proof.\ \par 
Let $B_{j}:=B(x_{j},r(x_{j}))\in {\mathcal{G}}$ and suppose that
 $\displaystyle x\in \bigcap_{j=1}^{k}{B(x_{j},5r(x_{j}))}.$ Then we have\ \par 
\quad \quad \quad $\displaystyle \forall j=1,...,k,\ d(x,x_{j})\leq 5r(x_{j})$\ \par 
hence\ \par 
\quad \quad \quad $\displaystyle d(x_{j},x_{l})\leq d(x_{j},x)+d(x,x_{l})\leq 5(r(x_{j})+r(x_{l}))\leq
 \frac{1}{4}(R(x_{j})+R(x_{l}))\Rightarrow R(x_{j})\leq 4R(x_{l})$\ \par 
and by exchanging $\displaystyle x_{j}$ and $\displaystyle x_{l},\
 R(x_{l})\leq 4R(x_{j}).$\ \par 
So we get\ \par 
\quad \quad \quad $\displaystyle \forall j,l=1,...,k,\ r(x_{j})\leq 4r(x_{l}),\
 r(x_{l})\leq 4r(x_{j}).$\ \par 
Now the ball $\displaystyle B(x_{j},5r(x_{j})+5r(x_{l}))$ contains
 $\displaystyle x_{l}$ hence the ball $\displaystyle B(x_{j},5r(x_{j})+6r(x_{l}))$
 contains the ball $\displaystyle B(x_{l},r(x_{l})).$ But, because
 $\displaystyle r(x_{l})\leq 4r(x_{j}),$ we get\ \par 
\quad \quad \quad $\displaystyle B(x_{j},5r(x_{j})+6{\times}4r(x_{j}))=B(x_{j},r(x_{j})(5+24))\supset
 B(x_{l},r(x_{l})).$\ \par 
The balls in ${\mathcal{G}}$ being disjoint, we get, setting
 $\displaystyle B_{l}:=B(x_{l},\ r(x_{l})),$\ \par 
\quad \quad \quad $\displaystyle \sum_{j=1}^{k}{\mathrm{V}\mathrm{o}\mathrm{l}(B_{l})}\leq
 \mathrm{V}\mathrm{o}\mathrm{l}(B(x_{j},29r(x_{j}))).$\ \par 
\quad The Lebesgue measure read in the chart $\varphi $  and the canonical
 measure $dv_{g}$ on $\displaystyle B(x,R_{\epsilon }(x))$ are
 equivalent ; precisely because of condition 1) in the admissible
 ball definition, we get that :\ \par 
\quad \quad \quad $\displaystyle (1-\epsilon )^{n}\leq \left\vert{\mathrm{d}\mathrm{e}\mathrm{t}g}\right\vert
 \leq (1+\epsilon )^{n},$\ \par 
and the measure $dv_{g}$ read in the chart $\varphi $ is $dv_{g}={\sqrt{\left\vert{\mathrm{d}\mathrm{e}\mathrm{t}g_{ij}}\right\vert
 }}d\xi ,$ where $\displaystyle d\xi $ is the Lebesgue measure
 in ${\mathbb{R}}^{n}.$ In particular :\ \par 
\quad \quad \quad $\displaystyle \forall x\in M,\ \mathrm{V}\mathrm{o}\mathrm{l}(B(x,\
 R_{\epsilon }(x)))\leq (1+\epsilon )^{n/2}\nu _{n}R^{n},$\ \par 
where $\nu _{n}$ is the euclidean volume of the unit ball in
 ${\mathbb{R}}^{n}.$\ \par 
\quad Now because $\displaystyle R(x_{j})$ is the admissible  radius
 and $\displaystyle 4{\times}29r(x_{j})<R(x_{j}),$ we have\ \par 
\quad \quad \quad $\displaystyle \mathrm{V}\mathrm{o}\mathrm{l}(B(x_{j},29r(x_{j})))\leq
 29^{n}(1+\epsilon )^{n/2}v_{n}r(x_{j})^{n}.$\ \par 
On the other hand we have also\ \par 
\quad \quad \quad $\displaystyle \mathrm{V}\mathrm{o}\mathrm{l}(B_{l})\geq v_{n}(1-\epsilon
 )^{n/2}r(x_{l})^{n}\geq v_{n}(1-\epsilon )^{n/2}4^{-n}r(x_{j})^{n},$\ \par 
hence\ \par 
\quad \quad \quad $\displaystyle \sum_{j=1}^{k}{(1-\epsilon )^{n/2}4^{-n}r(x_{j})^{n}}\leq
 29^{n}(1+\epsilon )^{n/2}r(x_{j})^{n},$\ \par 
so finally\ \par 
\quad \quad \quad $\displaystyle k\leq (29{\times}4)^{n}\frac{(1+\epsilon )^{n/2}}{(1-\epsilon
 )^{n/2}},$\ \par 
which means that $\displaystyle T\leq \frac{(1+\epsilon )^{n/2}}{(1-\epsilon
 )^{n/2}}(120)^{n}.$\ \par 
\quad Saying that any $\displaystyle x\in M$ belongs to at most $T$
 balls of the covering $\displaystyle \lbrace B_{j}\rbrace $
 means that $\sum_{j\in {\mathbb{N}}}{{\11}_{B_{j}}(x)}\leq T,$
 and this implies easily that :\ \par 
\quad \quad \quad $\displaystyle \forall f\in L^{1}(M),\ \sum_{j\in {\mathbb{N}}}{\int_{B_{j}}{\left\vert{f(x)}\right\vert
 dv_{g}(x)}}\leq T{\left\Vert{f}\right\Vert}_{L^{1}(M)}.$ $\blacksquare $\ \par 

\begin{lem}
~\label{HCN48} Let $(M,g)$ be a non compact connected complete
 riemannian manifold and ${\mathcal{C}}:=\lbrace B_{j}\rbrace
 _{j\in {\mathbb{N}}}$ a Vitali covering of $M$ with balls of
 radius less than $\displaystyle \delta >0.$ For any compact
 set $K$ in $M$ covered by ${\mathcal{O}}:=\bigcup_{k\in F_{K}}{B_{k}},$
 with $\displaystyle F_{K}$ finite, we can find a compact set
 $\displaystyle K'\supset K$ such that $\displaystyle \partial
 K'$ can be  covered by elements of ${\mathcal{C}}$ not intersecting
 $\bar {\mathcal{O}}.$
\end{lem}
\quad Proof.\ \par 
If this was not the case then there is a compact $K$ covered
 by $\displaystyle {\mathcal{O}}:=\bigcup_{k\in F_{K}}{B_{k}}$
 and such that for any compact $\displaystyle K'\supset K$ and
 any covering of $\displaystyle \partial K'$ by elements $\displaystyle
 B_{k}$ of ${\mathcal{C}},$ then $B_{k}\cap \bar {\mathcal{O}}\neq
 \emptyset .$ Because the balls have radius less than $\delta
 ,$ this means that $\displaystyle \partial K'$ is at most at
 a distance $\displaystyle 2\delta $ of ${\mathcal{O}}$ hence
 $M$ is bounded, hence the completeness of $M$ implies that $M$
 is compact. $\blacksquare $\ \par 
\quad Clearly the assumption that the radii are uniformly bounded is
 necessary as the example of ${\mathbb{R}}^{n}$ shows.\ \par 

\section{Sobolev spaces.~\label{CL27}}
\quad We have to define the Sobolev spaces in our setting, following
 E. Hebey~[\cite{Hebey96}, p. 10].\ \par 
First define the covariant derivatives by $\displaystyle (\nabla
 u)_{j}:=\partial _{j}u$ in local coordinates, while the components
 of $\nabla ^{2}u$ are given by\ \par 
\quad \quad \quad \begin{equation}  (\nabla ^{2}u)_{ij}=\partial _{ij}u-\Gamma
 ^{k}_{ij}\partial _{k}u,\label{HCF40}\end{equation}\ \par 
with the convention that we sum over repeated index. The Christoffel
 $\displaystyle \Gamma ^{k}_{ij}$ verify~\cite{BerGauMaz71} :\ \par 
\quad \quad \quad \begin{equation}  \Gamma ^{k}_{ij}=\frac{1}{2}g^{il}(\frac{\partial
 g_{kl}}{\partial x^{j}}+\frac{\partial g_{lj}}{\partial x^{k}}-\frac{\partial
 g_{jk}}{\partial x^{l}}).\label{HCF39}\end{equation}\ \par 
If $\displaystyle k\in {\mathbb{N}}$ and $r\geq 1$ are given,
 we denote by ${\mathcal{C}}^{r}_{k}(M)$ the space of smooth
 functions $u\in {\mathcal{C}}^{\infty }(M)$ such that $\displaystyle
 \ \left\vert{\nabla ^{j}u}\right\vert \in L^{r}(M)$ for $\displaystyle
 j=0,...,k.$ Hence\ \par 
\quad \quad \quad $\displaystyle {\mathcal{C}}^{r}_{k}(M):=\lbrace u\in {\mathcal{C}}^{\infty
 }(M),\ \forall j=0,...,k,\ \int_{M}{\left\vert{\nabla ^{j}u}\right\vert
 ^{r}dv_{g}}<\infty \rbrace $\ \par 
Now we have~\cite{Hebey96}\ \par 

\begin{defin}
The Sobolev space $\displaystyle W^{k,r}(M)$ is the completion
 of ${\mathcal{C}}^{r}_{k}(M)$ with respect to the norm :\par 
\quad \quad \quad $\displaystyle {\left\Vert{u}\right\Vert}_{W^{k,r}(M)}=\sum_{j=0}^{k}{{\left({\int_{M}{\left\vert{\nabla
 ^{j}u}\right\vert ^{r}dv_{g}}}\right)}^{1/r}}.$
\end{defin}
\quad We shall be interested only by $\displaystyle k\leq 2$ and we
 extend in a natural way this definition to the case of $p$ forms.\ \par 
Let the Sobolev exponents $\displaystyle S_{k}(r)$ as in the
 definition~\ref{HC33}, then the $k$ th Sobolev embedding is
 true if we have\ \par 
\quad \quad \quad $\displaystyle \forall u\in W^{k,r}(M),\ u\in L^{S_{k}(r)}(M).$\ \par 
This is the case in ${\mathbb{R}}^{n},$ or if $M$ is compact,
 or if $M$ has a Ricci curvature bounded from below and $\displaystyle
 \inf \ _{x\in M}v_{g}(B_{x}(1))\geq \delta >0,$ due to Varopoulos~\cite{Varopoulos89},
 see~\cite{Hebey96} theorem 3.14, p. 31].\ \par 

\begin{lem}
~\label{HCS43}We have the Sobolev comparison estimates where
 $\displaystyle B(x,R)$ is a $\epsilon $ admissible ball in $M$
 and $\varphi \ :\ B(x,R)\rightarrow {\mathbb{R}}^{n}$  is the
 admissible chart relative to $\displaystyle B(x,R),$\par 
\quad \quad \quad $\displaystyle \forall u\in W^{2,r}(B(x,R)),\ {\left\Vert{u}\right\Vert}_{W^{2,r}(B(x,R))}\leq
 (1+\epsilon C){\left\Vert{u\circ \varphi ^{-1}}\right\Vert}_{W^{2,r}(\varphi
 (B(x,R)))},$\par 
and, with $\displaystyle B_{e}(0,t)$ the euclidean ball in ${\mathbb{R}}^{n}$
 centered at $0$ and of radius $\displaystyle t,$\par 
\quad \quad \quad $\displaystyle {\left\Vert{v}\right\Vert}_{W^{2,r}(B_{e}(0,(1-\epsilon
 )R))}\leq (1+2C\epsilon ){\left\Vert{u}\right\Vert}_{W^{2,r}(B(x,R))}.$
\end{lem}
\quad Proof.\ \par 
We have to compare the norms of $\displaystyle u,\ \nabla u,\
 \nabla ^{2}u$ with the corresponding ones for $\displaystyle
 v:=u\circ \varphi ^{-1}$ in ${\mathbb{R}}^{n}.$\ \par 
First we have because $\displaystyle (1-\epsilon )\delta _{ij}\leq
 g_{ij}\leq (1+\epsilon )\delta _{ij}$ in $\displaystyle B(x,R)$ :\ \par 
\quad \quad \quad $\displaystyle B_{e}(0,(1-\epsilon )R)\subset \varphi (B(x,R))\subset
 B_{e}(0,(1+\epsilon )R).$\ \par 
\quad Because $\displaystyle \sum_{\left\vert{\beta }\right\vert =1}{\sup
 \ _{i,j=1,...,n,\ y\in B_{x}(R)}\left\vert{\partial ^{\beta
 }g_{ij}(y)}\right\vert }\leq \epsilon $ in $\displaystyle B(x,R),$
 we have the estimates, with $\displaystyle \forall y\in B(x,R),\
 z:=\varphi (y),$\ \par 
\quad \quad \quad $\displaystyle \forall y\in B(x,R),\ \left\vert{u(y)}\right\vert
 =\left\vert{v(z)}\right\vert ,\ \ \left\vert{\nabla u(y)}\right\vert
 \leq (1+C\epsilon )\left\vert{\partial v(z)}\right\vert .$\ \par 
\quad Because of~(\ref{HCF39}) and~(\ref{HCF40}) we get\ \par 
\quad \quad \quad $\displaystyle \forall y\in B(x,R),\ \left\vert{\nabla ^{2}u(y)}\right\vert
 \leq \left\vert{\partial ^{2}v(z)}\right\vert +\epsilon C\left\vert{\partial
 v(z)}\right\vert .$\ \par 
\quad Integrating this we get\ \par 
\quad \quad \quad $\displaystyle {\left\Vert{\nabla ^{2}u}\right\Vert}_{L^{r}(B(x,R))}\leq
 {\left\Vert{\left\vert{\partial ^{2}v}\right\vert +\epsilon
 C\left\vert{\partial v}\right\vert }\right\Vert}_{L^{r}(B_{e}(0,(1+\epsilon
 )R))}\leq {\left\Vert{\partial ^{2}v}\right\Vert}_{L^{r}(B_{e}(0,(1+\epsilon
 )R))}+C\epsilon {\left\Vert{\partial v}\right\Vert}_{L^{r}(B_{e}(0,(1+\epsilon
 )R))},$\ \par 
and\ \par 
\quad \quad \quad $\displaystyle {\left\Vert{\nabla u}\right\Vert}_{L^{r}(B(x,R))}\leq
 (1+C\epsilon ){\left\Vert{\partial v}\right\Vert}_{L^{r}(B_{e}(0,(1+\epsilon
 )R))}.$\ \par 
We also have the reverse estimates\ \par 
\quad \quad \quad $\displaystyle {\left\Vert{\partial ^{2}v}\right\Vert}_{L^{r}(B_{e}(0,(1-\epsilon
 )R))}\leq {\left\Vert{\nabla ^{2}u}\right\Vert}_{L^{r}(B(x,R))}+C\epsilon
 {\left\Vert{\nabla u}\right\Vert}_{L^{r}(B(x,R))},$\ \par 
and\ \par 
\quad \quad \quad $\displaystyle {\left\Vert{\partial v}\right\Vert}_{L^{r}(B_{e}(0,(1-\epsilon
 )R))}\leq (1+C\epsilon ){\left\Vert{\nabla u}\right\Vert}_{L^{r}(B(x,R))}.$\
 \par 
So, using that\ \par 
\quad \quad \quad $\displaystyle {\left\Vert{u}\right\Vert}_{W^{2,r}(B(x,R))}={\left\Vert{\nabla
 ^{2}u}\right\Vert}_{L^{r}(B(x,R))}+{\left\Vert{\nabla u}\right\Vert}_{L^{r}(B(x,R))}+{\left\Vert{u}\right\Vert}_{L^{r}(B(x,R))},$\
 \par 
we get\ \par 
\quad \quad \quad $\displaystyle {\left\Vert{u}\right\Vert}_{W^{2,r}(B(x,R))}^{r}\leq
 {\left\Vert{\partial ^{2}v}\right\Vert}_{L^{r}(B_{e}(0,(1+\epsilon
 )R))}+C\epsilon {\left\Vert{\partial v}\right\Vert}_{L^{r}(B_{e}(0,(1+\epsilon
 )R))}+(1+C\epsilon ){\left\Vert{\partial v}\right\Vert}_{L^{r}(B_{e}(0,(1+\epsilon
 )R))}+$\ \par 
\quad \quad \quad \quad \quad \quad \quad \quad \quad \quad \quad \quad \quad \quad \quad $\displaystyle +{\left\Vert{v}\right\Vert}_{L^{r}(B_{e}(0,(1+\epsilon
 )R))}\leq $\ \par 
\quad \quad \quad \quad \quad \quad \quad \quad \quad \quad \quad $\displaystyle \leq (1+2\epsilon C){\left\Vert{v}\right\Vert}_{W^{2,r}(B_{e}(0,(1+\epsilon
 )R))}.$\ \par 
Of course all these estimates can be reversed so we also have\ \par 
\quad \quad \quad $\displaystyle {\left\Vert{v}\right\Vert}_{W^{2,r}(B_{e}(0,(1-\epsilon
 )R))}\leq (1+2C\epsilon ){\left\Vert{u}\right\Vert}_{W^{2,r}(B(x,R))}.$\ \par 
\quad This ends the proof of the lemma. $\blacksquare $\ \par 

\begin{lem}
~\label{3S4}Let $\displaystyle B:=B(x,R)$ be a $\epsilon $ admissible
 ball in $M,$ we have the punctual estimates in $B$\par 
(i)            $\displaystyle \exists C>0,\ \forall \chi \in
 {\mathcal{C}}^{2}(B),\ \forall u\in {\mathcal{C}}^{2}_{p}(B),\
 \left\vert{\nabla (\chi u)}\right\vert \leq (1+C\epsilon )(\left\vert{\chi
 }\right\vert \left\vert{\nabla u}\right\vert +\left\vert{\nabla
 \chi }\right\vert \left\vert{u}\right\vert ).$\par 
(ii)           $\exists C>0,\ \forall \chi \in {\mathcal{C}}^{2}(B),\
 \forall u\in {\mathcal{C}}^{2}_{p}(B),\ \left\vert{\nabla ^{2}(\chi
 u)}\right\vert \leq (1+C\epsilon )(\left\vert{\chi }\right\vert
 \left\vert{\nabla ^{2}u}\right\vert +\left\vert{\nabla ^{2}\chi
 }\right\vert \left\vert{u}\right\vert +\left\vert{\nabla \chi
 }\right\vert \left\vert{\nabla u}\right\vert ).$
\end{lem}
\quad Proof.\ \par 
We have to compare the modulus of $\displaystyle u,\ \nabla u,\
 \nabla ^{2}u$ with the corresponding ones for $\displaystyle
 v:=u\circ \varphi ^{-1}$ in ${\mathbb{R}}^{n}.$\ \par 
First we have because $\displaystyle (1-\epsilon )\delta _{ij}\leq
 g_{ij}\leq (1+\epsilon )\delta _{ij}$ in $\displaystyle B(x,R)$ :\ \par 
\quad \quad \quad $\displaystyle B_{e}(0,(1-\epsilon )R)\subset \varphi (B(x,R))\subset
 B_{e}(0,(1+\epsilon )R).$\ \par 
Because $\displaystyle \sum_{\left\vert{\beta }\right\vert =1}{\sup
 \ _{i,j=1,...,n,\ y\in B_{x}(R)}\left\vert{\partial ^{\beta
 }g_{ij}(y)}\right\vert }\leq \epsilon $ in $\displaystyle B(x,R),$
 we have the estimates, with $\displaystyle \forall y\in B(x,R),\
 z:=\varphi (y),$\ \par 
\quad \quad \quad $\displaystyle \forall y\in B(x,R),\ \left\vert{u(y)}\right\vert
 =\left\vert{v(z)}\right\vert ,\ \ \left\vert{\nabla u(y)}\right\vert
 \leq (1+C\epsilon )\left\vert{\partial v(z)}\right\vert .$\ \par 
Now replacing $u$ by $\displaystyle \chi u,$ clearly we have
 for $\displaystyle \chi \circ \varphi ^{-1}v\circ \varphi ^{-1}$
 what we want, just using Leibnitz rule. Then the computations
 above gives the existence of a new constant $C$ such that\ \par 
\quad \quad \quad $\displaystyle \left\vert{\nabla (\chi u)}\right\vert \leq (1+C\epsilon
 )(\left\vert{\chi }\right\vert \left\vert{\nabla u}\right\vert
 +\left\vert{\nabla \chi }\right\vert \left\vert{u}\right\vert )$\ \par 
at all point of $\displaystyle B$ which gives \emph{(i)}.\ \par 
\ \par 
\quad Because of~(\ref{HCF39}) and~(\ref{HCF40}) we get\ \par 
\quad \quad \quad $\displaystyle \forall y\in B(x,R),\ \left\vert{\nabla ^{2}u(y)}\right\vert
 \leq \left\vert{\partial ^{2}v(z)}\right\vert +\epsilon C\left\vert{\partial
 v(z)}\right\vert .$\ \par 
Now replacing $u$ by $\displaystyle \chi u,$ clearly we have
 for $\displaystyle \chi \circ \varphi ^{-1}v\circ \varphi ^{-1}$
 what we want, again just using Leibnitz rule. Then the computations
 above gives the existence of a new constant $C$ such that\ \par 
\quad \quad \quad $\displaystyle \left\vert{\nabla ^{2}(\chi u)}\right\vert \leq
 (1+C\epsilon )(\left\vert{\chi \nabla ^{2}u}\right\vert +\left\vert{\nabla
 ^{2}\chi u}\right\vert +\left\vert{\nabla \chi }\right\vert
 \left\vert{\nabla u}\right\vert )$\ \par 
at all point of $\displaystyle B$ which gives \emph{(ii)} and
 ends the proof of this lemma. $\blacksquare $\ \par 
\ \par 
\quad We have to study the behavior of the Sobolev embeddings w.r.t.
 the radius. Set $\displaystyle B_{R}:=B_{e}(0,R).$\ \par 

\begin{lem}
~\label{3S3}We have, with $\displaystyle s=S_{1}(r),\ t=S_{2}(r),$\par 
(i)            $\displaystyle \forall R,\ 0<R\leq 1,\ \forall
 u\in W^{2,r}(B_{R}),\ {\left\Vert{u}\right\Vert}_{L^{t}(B_{R})}\leq
 CR^{-2}\ {\left\Vert{u}\right\Vert}_{W^{2,r}(B_{R})}$\par 
and\par 
(ii)           $\displaystyle \forall R,\ 0<R\leq 1,\ \forall
 u\in W^{2,r}(B_{R}),\ {\left\Vert{\partial u}\right\Vert}_{L^{s}(B_{R})}\leq
 CR^{-1}\ {\left\Vert{u}\right\Vert}_{W^{2,r}(B_{R})}$\par 
the constant $C$ depending only on $\displaystyle n,\ r.$
\end{lem}
\quad Proof.\ \par 
We start with $\displaystyle R=1,$ then we have by Sobolev embeddings
 with $\displaystyle t=S_{2}(r),$\ \par 
\quad \quad \quad \begin{equation}  \forall v\in W^{2,r}(B_{1}),\ {\left\Vert{v}\right\Vert}_{L^{t}(B_{1})}\leq
 C{\left\Vert{v}\right\Vert}_{W^{2,r}(B_{1})}\label{aS0}\end{equation}\ \par 
where $\displaystyle C$ depends only on $n.$ For $\displaystyle
 u\in W^{2,r}(B_{R})$ we set\ \par 
\quad \quad \quad $\displaystyle \forall x\in B_{1},\ y:=Rx\in B_{R},\ v(x):=u(y).$\ \par 
Then we have\ \par 
\quad \quad \quad $\displaystyle \partial v(x)=\partial u(y){\times}\frac{\partial
 y}{\partial x}=R\partial u(y);\ \partial ^{2}v(x)=\partial ^{2}u(y){\times}(\frac{\partial
 y}{\partial x})^{2}=R^{2}\partial ^{2}u(y).$\ \par 
So we get, because the jacobian for this change of variables
 is $\displaystyle R^{-n},$\ \par 
\quad \quad \quad $\displaystyle {\left\Vert{\partial v}\right\Vert}_{L^{r}(B_{1})}^{r}=\int_{B_{1}}{\left\vert{\partial
 v(x)}\right\vert ^{r}dm(x)}=\int_{B_{R}}{\left\vert{\partial
 u(y)}\right\vert ^{r}\frac{R^{r}}{R^{n}}dm(x)}=R^{r-n}{\left\Vert{\partial
 u}\right\Vert}_{L^{r}(B_{R})}^{r}.$\ \par 
So\ \par 
\quad \quad \quad \begin{equation} {\left\Vert{ \partial u}\right\Vert}_{L^{r}(B_{R})}=R^{-1+n/r}{\left\Vert{\partial
 v}\right\Vert}_{L^{r}(B_{1})}.\label{HS50}\end{equation}\ \par 
The same way we get\ \par 
\quad \quad \quad \begin{equation} {\left\Vert{ \partial ^{2}u}\right\Vert}_{L^{r}(B_{R})}=R^{-2+n/r}{\left\Vert{\partial
 ^{2}v}\right\Vert}_{L^{r}(B_{1})}\label{HS51}\end{equation}\ \par 
and of course\ \par 
\quad \quad \quad $\displaystyle {\left\Vert{u}\right\Vert}_{L^{r}(B_{R})}=R^{n/r}{\left\Vert{v}\right\Vert}_{L^{r}(B_{1})}.$\
 \par 
So with~\ref{aS0} we get\ \par 
\quad \quad \quad \begin{equation} {\left\Vert{ u}\right\Vert}_{L^{t}(B_{R})}=R^{n/t}{\left\Vert{v}\right\Vert}_{L^{t}(B_{1})}\leq
 CR^{n/t}{\left\Vert{v}\right\Vert}_{W^{k,r}(B_{1})}.\label{3S1}\end{equation}\
 \par 
But\ \par 
\quad \quad \quad $\displaystyle {\left\Vert{u}\right\Vert}_{W^{2,r}(B_{R})}:={\left\Vert{u}\right\Vert}_{L^{r}(B_{R})}+{\left\Vert{\partial
 u}\right\Vert}_{L^{r}(B_{R})}+{\left\Vert{\partial ^{2}u}\right\Vert}_{L^{r}(B_{R})},$\
 \par 
and\ \par 
\quad \quad \quad $\displaystyle {\left\Vert{v}\right\Vert}_{W^{2,r}(B_{1})}:={\left\Vert{v}\right\Vert}_{L^{r}(B_{1})}+{\left\Vert{\partial
 v}\right\Vert}_{L^{r}(B_{1})}+{\left\Vert{\partial ^{2}v}\right\Vert}_{L^{r}(B_{1})},$\
 \par 
so\ \par 
\quad \quad \quad $\displaystyle {\left\Vert{v}\right\Vert}_{W^{2,r}(B_{1})}:=R^{-n/r}{\left\Vert{u}\right\Vert}_{L^{r}(B_{R})}+R^{1-n/r}{\left\Vert{\partial
 u}\right\Vert}_{L^{r}(B_{R})}+R^{2-n/r}{\left\Vert{\partial
 ^{2}u}\right\Vert}_{L^{r}(B_{R})}.$\ \par 
Because we have $\displaystyle R\leq 1,$ we get\ \par 
\quad \quad \quad $\displaystyle {\left\Vert{v}\right\Vert}_{W^{2,r}(B_{1})}\leq
 R^{-n/r}({\left\Vert{u}\right\Vert}_{L^{r}(B_{R})}+{\left\Vert{\partial
 u}\right\Vert}_{L^{r}(B_{R})}+{\left\Vert{\partial ^{2}u}\right\Vert}_{L^{r}(B_{R})})=R^{-n/r}{\left\Vert{u}\right\Vert}_{W^{2,r}(B_{R})}.$\
 \par 
Putting it in~(\ref{3S1}) we get\ \par 
\quad \quad \quad $\displaystyle {\left\Vert{u}\right\Vert}_{L^{t}(B_{R})}\leq
 CR^{n/t}{\left\Vert{v}\right\Vert}_{W^{k,r}(B_{1})}\leq CR^{-n(\frac{1}{r}-\frac{1}{t})}{\left\Vert{u}\right\Vert}_{W^{2,r}(B_{R})}.$\
 \par 
But, because $\displaystyle t=S_{2}(r),$ we get $\displaystyle
 (\frac{1}{r}-\frac{1}{t})=\frac{2}{n}$ and\ \par 
\quad \quad \quad $\displaystyle {\left\Vert{u}\right\Vert}_{L^{t}(B_{R})}\leq
 CR^{-2}{\left\Vert{u}\right\Vert}_{W^{2,r}(B_{R})}.$\ \par 
\ \par 
\quad To have the \emph{(ii)} we proceed exactly the same way.\ \par 
We start with $\displaystyle R=1,$ then we have by Sobolev embeddings
 with $\displaystyle s=S_{1}(r),$\ \par 
\quad \quad \quad $\displaystyle \forall v\in W^{2,r}(B_{1}),\ {\left\Vert{\partial
 v}\right\Vert}_{L^{s}(B_{1})}\leq C{\left\Vert{v}\right\Vert}_{W^{2,r}(B_{1})}$\
 \par 
and this leads as above to\ \par 
\quad \quad \quad $\displaystyle {\left\Vert{\partial u}\right\Vert}_{L^{t}(B_{R})}\leq
 CR^{-1}{\left\Vert{u}\right\Vert}_{W^{2,r}(B_{R})}.$ $\blacksquare $\ \par 
\quad The constant $C$ depends only on $\displaystyle n,r.$\ \par 

\begin{lem}
~\label{3S2}Let $x\in M$ and $\displaystyle B(x,R)$ be a $\epsilon
 $ admissible ball ; we have, with $\displaystyle s=S_{1}(r),\ t=S_{2}(r),$\par 
(i)            $\displaystyle \forall u\in W^{2,r}(B(x,R)),\
 {\left\Vert{u}\right\Vert}_{L^{t}(B(x,R))}\leq CR^{-2}\ {\left\Vert{u}\right\Vert}_{W^{2,r}(B(x,R))},$\par
 
and\par 
(ii)           $\displaystyle \forall u\in W^{2,r}(B(x,R)),\
 {\left\Vert{\nabla u}\right\Vert}_{L^{s}(B(x,R))}\leq CR^{-1}\
 {\left\Vert{u}\right\Vert}_{W^{2,r}(B(x,R))},$\par 
the constant $\displaystyle C$ depending only on $\displaystyle
 n,\ r$ and $\displaystyle \epsilon .$
\end{lem}
\quad Proof.\ \par 
This is true in ${\mathbb{R}}^{n}$ by lemma~\ref{3S3} so we can
 apply the comparison lemma~\ref{HCS43}. $\blacksquare $\ \par 
\ \par 

\begin{lem}
~\label{HC42} Let $\displaystyle B:=B(0,R)$ be the ball in ${\mathbb{R}}^{n}$
 of center $0$ and radius $\displaystyle R\leq 1$ and $\displaystyle
 B'=B(0,R/2).$ Let $\displaystyle u\in L^{r}(B)$ such that $\displaystyle
 \Delta u\in L^{r}(B)$ then we have\par 
\quad \quad \quad $\displaystyle u\in W^{2,r}(B'),\ {\left\Vert{u}\right\Vert}_{W^{2,r}(B')}\leq
 c_{1}R^{-2}{\left\Vert{u}\right\Vert}_{L^{r}(B)}+c_{2}{\left\Vert{\Delta
 u}\right\Vert}_{L^{r}(B)},$\par 
where the constants $\displaystyle c_{1},\ c_{2}$ depend only
 on $\displaystyle n,r.$
\end{lem}
\quad Proof.\ \par 
We start with $\displaystyle R=1,$ then we have by the classical
 CZI for the usual laplacian $\Delta _{{\mathbb{R}}}$ in ${\mathbb{R}}^{n},$
 [ \cite{GuilbargTrudinger98}, Th. 9.11, p. 235]:\ \par 
\quad \quad \quad \begin{equation}  \forall v\in L^{r}(B),\ \Delta v\in L^{r}(B),\
 {\left\Vert{v}\right\Vert}_{W^{2,r}(B')}\leq c_{1}{\left\Vert{v}\right\Vert}_{L^{r}(B)}+c_{2}{\left\Vert{\Delta
 v}\right\Vert}_{L^{r}(B)},\label{HC41}\end{equation}\ \par 
the constants $\displaystyle c_{1},\ c_{2}$ depending only on
 $\displaystyle n,\ r.$\ \par 
\quad To go to any $R$ we take $u$ with the hypotheses of the lemma
 and we make the change of variables\ \par 
\quad \quad \quad $\displaystyle y=Rx,\ dm(y)=R^{n}dm(x),\ v(x):=u(Rx).$\ \par 
We set $\displaystyle B_{R}:=B_{e}(0,R)$ then we get\ \par 
\quad \quad \quad $\displaystyle {\left\Vert{v}\right\Vert}_{L^{r}(B_{1})}^{r}:=\int_{B_{1}}{\left\vert{v(x)}\right\vert
 ^{r}dm(x)}=\int_{B_{R}}{\left\vert{v(\frac{y}{R})}\right\vert
 ^{r}R^{-n}dm(y)}=\int_{B_{R}}{\left\vert{u(y)}\right\vert ^{r}R^{-n}dm(y)}=R^{-n}{\left\Vert{u}\right\Vert}_{L^{r}(B_{R})}^{r}.$\
 \par 
And\ \par 
\quad \quad \quad $\displaystyle \partial _{i}v(x)=\partial _{i}u(Rx)R,\ \partial
 _{ij}v(x)=\partial _{ij}u(Rx)R^{2},$\ \par 
hence\ \par 
\quad \quad \quad $\displaystyle {\left\Vert{\partial _{i}v}\right\Vert}_{L^{r}(B_{1})}^{r}:=\int_{B_{1}}{\left\vert{\partial
 _{i}v(x)}\right\vert ^{r}dm(x)}=$\ \par 
\quad \quad \quad \quad \quad $\displaystyle =\int_{B_{R}}{\left\vert{\partial _{i}v}\right\vert
 (\frac{y}{R})^{r}R^{-n}dm(y)}=\int_{B_{R}}{R^{r}\left\vert{\partial
 _{i}u(y)}\right\vert ^{r}R^{-n}dm(y)}=R^{r-n}{\left\Vert{\partial
 _{i}u}\right\Vert}_{L^{r}(B_{R})}^{r}.$\ \par 
\quad \quad \quad $\displaystyle {\left\Vert{\partial _{ij}v}\right\Vert}_{L^{r}(B_{1})}^{r}:=\int_{B_{1}}{\left\vert{\partial
 _{ij}v(x)}\right\vert ^{r}dm(x)}=$\ \par 
\quad \quad \quad \quad \quad $\displaystyle =\int_{B_{R}}{\left\vert{\partial _{ij}v}\right\vert
 (\frac{y}{R})^{r}R^{-n}dm(y)}=\int_{B_{R}}{R^{2r}\left\vert{\partial
 _{ij}u(y)}\right\vert ^{r}R^{-n}dm(y)}=R^{2r-n}{\left\Vert{\partial
 _{ij}u}\right\Vert}_{L^{r}(B_{R})}^{r}.$\ \par 
So we get by~(\ref{HC41})\ \par 
\quad \quad \quad $\displaystyle (R/2)^{1-n/r}{\left\Vert{\partial _{i}u}\right\Vert}_{L^{r}(B_{R/2})}\leq
 c_{1}(r)R^{-n/r}{\left\Vert{u}\right\Vert}_{L^{r}(B_{R})}+c_{2}(R)R^{2-n/r}{\left\Vert{\Delta
 u}\right\Vert}_{L^{r}(B_{R})}$\ \par 
hence\ \par 
\quad \quad \quad $\displaystyle {\left\Vert{\partial _{i}u}\right\Vert}_{L^{r}(B_{R/2})}\leq
 2^{-1+n/r}(c_{1}(r)R^{-1}{\left\Vert{u}\right\Vert}_{L^{r}(B_{R})}+c_{2}(r)R{\left\Vert{\Delta
 u}\right\Vert}_{L^{r}(B_{R})}).$\ \par 
And the same way\ \par 
\quad \quad \quad $\displaystyle {\left\Vert{\partial _{ij}u}\right\Vert}_{L^{r}(B_{R/2})}\leq
 2^{-2+n/r}(c_{1}(r)R^{-2}{\left\Vert{u}\right\Vert}_{L^{r}(B_{R})}+c_{2}(r){\left\Vert{\Delta
 u}\right\Vert}_{L^{r}(B_{R})}).$\ \par 
So we get finally\ \par 
\quad \quad \quad $\displaystyle {\left\Vert{u}\right\Vert}_{W^{2,r}(B_{R/2})}\leq
 c_{1}(n,r)R^{-2}{\left\Vert{u}\right\Vert}_{L^{r}(B_{R})}+c_{2}(n,r){\left\Vert{\Delta
 u}\right\Vert}_{L^{r}(B_{R})}),$\ \par 
where the constants $\displaystyle c_{1},\ c_{2}$ depend only
 on $\displaystyle n,r.$ $\blacksquare $\ \par 

\section{Local estimates for the laplacian.}
\quad All these local estimates are quite well known. I reprove them
 here to precise the notations and the dependences of the constants.\ \par 

\begin{lem}
~\label{CF1}Let $(M,g)$ be a riemannian manifold. For $x\in M,\
 \epsilon >0,$ we take a $\epsilon $ admissible  ball $B_{x}(R).$
 Then there is a $0<\epsilon _{0}\leq \epsilon ,$ hence a $R=R_{\epsilon
 _{0}}(x)>0,$ and a constant $C$ depending only on $n=\mathrm{d}\mathrm{i}\mathrm{m}_{{\mathbb{R}}}M,\
 r$ and $\epsilon _{0}$ such that :\par 
\quad \quad \quad $\displaystyle \forall \omega \in L_{p}^{r}(B_{x}(R)),\ \exists
 u\in W_{p}^{2,r}(B_{x}(R))::\Delta u=\omega ,\ {\left\Vert{u}\right\Vert}_{W^{2,r}(B_{x}(R))}\leq
 C{\left\Vert{\omega }\right\Vert}_{L^{r}(B_{x}(R))}.$\par 
Moreover $u$ is linear in $\omega .$
\end{lem}
\quad Proof.\ \par 
For $x\in M$ and $\epsilon >0,$ we take a $\epsilon $ admissible
  ball $B_{x}(R)$ and a chart $\varphi \ :\ (y_{1},...,y_{n}),$
 which means:\ \par 
\quad 1) $(1-\epsilon )\delta _{ij}\leq g_{ij}\leq (1+\epsilon )\delta
 _{ij}$ in $B_{x}(R)$ as bilinear forms,\ \par 
\quad 2) $\sum_{\left\vert{\beta }\right\vert =1}{\sup \ _{i,j=1,...,n,\
 y\in B_{x}(R)}\left\vert{\partial ^{\beta }g_{ij}(y)}\right\vert
 }\leq \epsilon .$\ \par 
Of course the operator $d$ on $p$ forms is local and so is $d^{*}$
 as a first order differential operator.\ \par 
\quad So the Hodge laplacian $\Delta _{\varphi }$ read by $\varphi
 $ in $U:=\varphi (B_{x}(R))$ is still a second order partial
 differential system of operators and with $\Delta _{{\mathbb{R}}}$
 the usual laplacian in ${\mathbb{R}}^{n}$ acting on forms in
 $U,$ we set: $A\omega _{\varphi }:=\Delta _{\varphi }\omega
 _{\varphi }-\Delta _{{\mathbb{R}}}\omega _{\varphi },$ where
 $\omega _{\varphi }$ is the $p$ form $\omega $ read in the chart
 $(B_{x}(R),\varphi )$ and $A$ is a matrix valued second order
 operator with ${\mathcal{C}}^{\infty }$ smooth coefficients
 such that $A\ :=\Delta _{\varphi }-\Delta _{{\mathbb{R}}}:\
 W^{2,r}(U)\rightarrow L^{r}(U).$\ \par 
\quad This difference $A$ is controlled by the derivatives of the metric
 tensor up to order $1;$ for instance for a function $f$ we have
 in the chart $\varphi $\ \par 
\quad \quad \quad $\displaystyle \Delta _{\varphi }f=\frac{1}{{\sqrt{\mathrm{d}\mathrm{e}\mathrm{t}(g_{ij})}}}\partial
 _{i}(g^{ij}{\sqrt{\mathrm{d}\mathrm{e}\mathrm{t}(g_{ij})}}\partial
 _{j}f)=g^{ij}\partial ^{2}_{ij}f+Y_{0}f,$\ \par 
where $Y_{0}$  is a first order differential operator depending
 on $g$ and its first derivatives;\ \par 
more generally for a $p$ form $u,$ still in the chart $\varphi
 ,$ the formula 21.23, p. 169 in~\cite{TaylorGD} gives  $\Delta
 _{\varphi }u=g^{ij}\partial ^{2}_{ij}u+Y_{p}u,$ where $Y_{p}$
 is a first order differential operator.\ \par 
So $\Delta _{\varphi }$ depends on the first order derivatives
 of $g,$ hence the difference $A:=\Delta _{\varphi }-\Delta _{{\mathbb{R}}},$
 where  $\Delta _{{\mathbb{R}}}u=\delta ^{ij}\partial ^{2}_{ij}u,$
 is controlled by the first order derivatives of $\displaystyle
 g,$ which, by condition 2), can be made as small as we wish. So we have\ \par 
\quad \quad \quad \begin{equation}  \ \forall y\in U,\ \left\vert{A(u)(y)}\right\vert
 \leq \left\vert{(g^{ij}(y)-\delta ^{ij})\partial ^{2}_{ij}u(y)}\right\vert
 +\left\vert{E(u)(y)}\right\vert ,\label{HC7}\end{equation}\ \par 
where $\displaystyle E$ is a first order partial differential
 operator whose coefficients depend on the first order derivatives
 of $g,$ and are $0$ for $y=x.$ So $\ \left\vert{E(u)(y)}\right\vert
 \leq \eta \left\vert{\nabla u(y)}\right\vert $ for $y\in U,$
 where $\eta $ is a continuous function of the metric $g$ and
 $\nabla g$ only; since $\ \left\vert{\nabla g}\right\vert \leq
 \epsilon ,\ \eta $ may be chosen to depend on $\epsilon >0$
 only and $\eta (0)=0.$\ \par 
Hence, integrating~(\ref{HC7}), we get:\ \par 
\quad \quad \quad $\displaystyle {\left\Vert{Au}\right\Vert}_{L^{r}(U)}\leq {\left\Vert{\nabla
 g}\right\Vert}_{L^{\infty }(U)}{\left\Vert{u}\right\Vert}_{W^{2,r}(U)}+\eta
 (\epsilon ){\left\Vert{\nabla u}\right\Vert}_{L^{r}(U)}.$\ \par 
\quad So, because ${\left\Vert{\nabla u}\right\Vert}_{L^{r}(U)}$ is
 controlled by ${\left\Vert{u}\right\Vert}_{W^{2,r}(U)},$ there
 is a $0\leq c(\epsilon ),\ c(0)=0$ and $c$ continuous at $0,$
 such that ${\left\Vert{Au}\right\Vert}_{L^{r}(U)}\leq c(\epsilon
 ){\left\Vert{u}\right\Vert}_{W^{2,r}(U)}.$\ \par 
\quad Let $\gamma $ be a $p$ form in $U\subset {\mathbb{R}}^{n}.$ We
 know that $\Delta _{{\mathbb{R}}}$ operates component-wise on
 the $p$ form $\gamma ,$ so we have:\ \par 
\quad \quad \quad $\displaystyle \forall \gamma \in L^{r}_{p}(U),\ \exists v_{0}\in
 W^{2,r}_{p}(U)::\Delta _{{\mathbb{R}}}v_{0}=\gamma ,\ {\left\Vert{v_{0}}\right\Vert}_{W^{2,r}(U)}\leq
 C{\left\Vert{\gamma }\right\Vert}_{L^{r}(U)},$\ \par 
simply setting the component of $\displaystyle v_{0}$ to be the
 Newtonian potential of the corresponding component of $\gamma
 $ in $U.$ These non trivial estimates are coming from Gilbarg
 and Trudinger~[\cite{GuilbargTrudinger98}, Th 9.9, p. 230] and
 the constant $C=C(n,r)$ depends only on $n$ and $r.$\ \par 
\quad So we get $\Delta _{{\mathbb{R}}}v_{0}+Av_{0}=\gamma +\gamma _{1},$ with\ \par 
\quad \quad \quad $\displaystyle \gamma _{1}=Av_{0}\Rightarrow {\left\Vert{\gamma
 _{1}}\right\Vert}_{L^{r}(U)}\leq c{\left\Vert{v_{0}}\right\Vert}_{W^{2,r}(U)}\leq
 cC{\left\Vert{\gamma }\right\Vert}_{L^{r}(U)}.$\ \par 
\quad We solve again\ \par 
\quad \quad \quad $\displaystyle \exists v_{1}\in W^{2,r}_{p}(U)::\Delta _{{\mathbb{R}}}v_{1}=\gamma
 _{1},\ {\left\Vert{v_{1}}\right\Vert}_{W^{2,r}(U)}\leq C{\left\Vert{\gamma
 _{1}}\right\Vert}_{L^{r}(U)}=C^{2}c{\left\Vert{\gamma }\right\Vert}_{L^{r}(U)},$\
 \par 
and we set\ \par 
\quad \quad \quad $\displaystyle \gamma _{2}:=Av_{1}\Rightarrow {\left\Vert{\gamma
 _{2}}\right\Vert}_{L^{r}(U)}\leq c{\left\Vert{v_{1}}\right\Vert}_{W^{2,r}(U)}\leq
 C{\left\Vert{\gamma _{1}}\right\Vert}_{L^{r}(U)}\leq C^{2}c^{2}{\left\Vert{\gamma
 }\right\Vert}_{L^{r}(U)}.$\ \par 
\quad And by induction:\ \par 
\quad \quad \quad $\displaystyle \forall k\in {\mathbb{N}},\ \gamma _{k}:=Av_{k-1}\Rightarrow
 {\left\Vert{\gamma _{k}}\right\Vert}_{L^{r}(U)}\leq c{\left\Vert{v_{k-1}}\right\Vert}_{W^{2,r}(U)}\leq
 C{\left\Vert{\gamma _{k-1}}\right\Vert}_{L^{r}(U)}\leq C^{k}c^{k}{\left\Vert{\gamma
 }\right\Vert}_{L^{r}(U)}$\ \par 
and\ \par 
\quad \quad \quad $\displaystyle \exists v_{k}\in W^{2,r}_{p}(U)::\Delta _{{\mathbb{R}}}v_{k}=\gamma
 _{k},\ {\left\Vert{v_{k}}\right\Vert}_{W^{2,r}(U)}\leq C{\left\Vert{\gamma
 _{k}}\right\Vert}_{L^{r}(U)}\leq C^{k+1}c^{k}{\left\Vert{\gamma
 }\right\Vert}_{L^{r}(U)}.$\ \par 
Now we set $v:=\sum_{j\in {\mathbb{N}}}{(-1)^{j}v_{j}}.$ This
 series converges in norm $W^{2,r}(U),$ provided that we choose
 $\epsilon _{0}\leq \epsilon $ small enough to have $c(\epsilon
 _{0})C^{2}<1,$ and we get:\ \par 
\quad \quad \quad $\displaystyle \Delta _{\varphi }v=\Delta _{{\mathbb{R}}}v+Av=\sum_{j\in
 {\mathbb{N}}}{(-1)^{j}(\Delta _{{\mathbb{R}}}v_{j}+Av_{j})}=\gamma ,$\ \par 
the last series converging in $L_{p}^{r}(U).$\ \par 
\quad Going back to the manifold $M$ with $\gamma :=\omega _{\varphi
 }$ and setting $u_{\varphi }:=v,$ we get the right estimates:\ \par 
\quad \quad \quad $\displaystyle \exists u\in W^{2,r}(B_{x}(R))::\Delta u=\omega
 \ in\ B_{x}(R),\ {\left\Vert{u}\right\Vert}_{W^{2,r}(B_{x}(R))}\leq
 C{\left\Vert{\omega }\right\Vert}_{L^{r}(B_{x}(R))},$\ \par 
because, by use of the comparison lemma already seen, the Sobolev
 spaces for $U$ go to the analogous Sobolev spaces for $B_{x}(R)$
 in $M.$ Moreover $C$ depends only on $n,\ r,\ \epsilon _{0}$
 and $u$ is linear in $\omega .$ $\blacksquare $\ \par 

\begin{lem}
~\label{lE1} Let $B_{R}:=B(0,R)$ be the ball in ${\mathbb{R}}^{n}$
 of center $0$ and radius $R\leq 1$ and $B'_{R}=B(0,3R/4).$ Suppose
 we have, with a constant $C$ depending only on $n,r$ and the
 ${\mathcal{C}}^{1}$ bound of the coefficients of $\Delta _{\varphi }$\par 
\quad \quad \quad $\displaystyle \forall v\in W^{2,r}(B_{1}),\ {\left\Vert{v}\right\Vert}_{W^{2,r}(B_{1}')}\leq
 C({\left\Vert{v}\right\Vert}_{L^{r}(B_{1})}+{\left\Vert{\Delta
 _{\varphi }v}\right\Vert}_{L^{r}(B_{1})}).$\par 
Let $u\in L^{r}(B_{R})$ such that $\Delta _{\varphi }u\in L^{r}(B_{R})$
 and  then we have:\par 
\quad \quad \quad $\displaystyle u\in W^{2,r}(B_{R}'),\ {\left\Vert{u}\right\Vert}_{W^{2,r}(B_{R}')}\leq
 c_{1}R^{-2}{\left\Vert{u}\right\Vert}_{L^{r}(B_{R})}+c_{2}{\left\Vert{\Delta
 _{\varphi }u}\right\Vert}_{L^{r}(B_{R})},$\par 
where the constants $c_{1},c_{2}$ depend only on $n,r$ and the
 ${\mathcal{C}}^{1}$ bound of the coefficients of $\Delta _{\varphi }.$
\end{lem}
\quad Proof.\ \par 
We start with $R=1,\ B:=B(0,1),$ then we have by assumption:\ \par 
\quad \quad \quad $\displaystyle \forall v\in W^{2,r}(B_{1}),\ {\left\Vert{v}\right\Vert}_{W^{2,r}(B_{1}')}\leq
 C({\left\Vert{v}\right\Vert}_{L^{r}(B_{1})}+{\left\Vert{\Delta
 _{\varphi }v}\right\Vert}_{L^{r}(B_{1})}),$\ \par 
the constants $\displaystyle C$ depending only on $n,r$ and the
 ${\mathcal{C}}^{1}$ bound of the coefficients of $\Delta _{\varphi
 }.$ It remains to make the simple change of variables $y=Rx,\
 dm(y)=R^{n}dm(x),\ v(x):=u(Rx)$ and to notice that $\partial
 _{j}v(x)=R\partial _{j}(u)(Rx),\ \partial ^{2}_{ij}v(x)=R^{2}\partial
 ^{2}_{ij}(u)(Rx)$ in the integrals defining the $L^{r}$ norm
 to get the result. $\blacksquare $\ \par 

\begin{lem}
~\label{lE0}Let $\Delta _{\varphi }$ be a second order elliptic
 matrix operator with ${\mathcal{C}}^{\infty }$ coefficients
 operating on  $p$ forms $v$ defined in $U\subset {\mathbb{R}}^{n}.$
 Let $B:=B(0,R)$ a ball in ${\mathbb{R}}^{n},\ B':=B(0,3R/4)$
 and suppose that $B\Subset U.$ Then we have an interior estimate:
 there are constants $c_{1},c_{2}$ depending only on $n=\mathrm{d}\mathrm{i}\mathrm{m}_{{\mathbb{R}}}M,\
 r$ and the ${\mathcal{C}}^{1}$norm of the coefficients of $\Delta
 _{\varphi }$ in $\bar B$ such that\par 
\quad \quad \quad \begin{equation}  \forall v\in W_{p}^{2,r}(B),\ {\left\Vert{v}\right\Vert}_{W^{2,r}(B')}\leq
 c_{1}R^{-2}{\left\Vert{v}\right\Vert}_{L^{r}(B)}+c_{2}{\left\Vert{\Delta
 _{\varphi }v}\right\Vert}_{L^{r}(B)}.\label{CM11}\end{equation}
\end{lem}
\quad Proof.\ \par 
For a $0$ form, this lemma is exactly theorem 9.11, in~\cite{GuilbargTrudinger98}
 plus lemma~\ref{lE1} to get the dependence in $R.$\ \par 
For $p$ forms we cannot avoid the use of deep results on elliptic
 systems of equations.\ \par 
\quad Let $\displaystyle v$ be a $p$ form in $B\subset {\mathbb{R}}^{n}.$
 We use the interior estimates in~[\cite{Morrey66}, \S  6.2,
 Thm 6.2.6]. In our context, second order elliptic system, and
 with our notations, with $r>1,$ we get:\ \par 
\quad \quad \quad $\displaystyle \exists C>0,\ \forall v\in W_{p}^{2,r}(B),\ {\left\Vert{v}\right\Vert}_{W^{2,r}(B')}\leq
 c_{1}R^{-2}{\left\Vert{v}\right\Vert}_{L^{r}(B)}+c_{2}{\left\Vert{\Delta
 _{\varphi }v}\right\Vert}_{L^{r}(B)},$\ \par 
already including the dependence in $R.$\ \par 
\quad The constants $c_{1},c_{2}$ depend only on $r,\ n:=\mathrm{d}\mathrm{i}\mathrm{m}M$
 and the bounds and moduli of continuity of all the coefficients
 of the matrix $\Delta _{\varphi }.$ (In Morrey's book, p. 213:
 the constant depends \emph{only on} $E$ \emph{and on} $E'.$)\ \par 
\quad In particular, if $\Delta _{\varphi }$ has its coefficients near
 those of $\Delta _{{\mathbb{R}}}$ in the ${\mathcal{C}}^{1}$
 norm, then the constants $c_{1},c_{2}$ are near the ones obtained
 for $\Delta _{{\mathbb{R}}}.$ $\blacksquare $\ \par 

\begin{lem}
~\label{HF0}Let $(M,g)$ be a riemannian manifold. For $x\in M,\
 \epsilon >0,$ we take a $\epsilon $ admissible  ball $B_{x}(R).$
 We have a local Calderon Zygmund inequality on the manifold
 $M.$ There are constants $\displaystyle c_{1},c_{2}$ depending
 only on $n=\mathrm{d}\mathrm{i}\mathrm{m}_{{\mathbb{R}}}M,\
 r$ and $\epsilon $ such that:\par 
\quad \quad \quad $\displaystyle \forall u\in W^{2,r}(B_{x}(R)),\ {\left\Vert{u}\right\Vert}_{W^{2,r}(B_{x}(R/2))}\leq
 c_{1}R^{-2}{\left\Vert{u}\right\Vert}_{L^{r}(B_{x}(R))}+c_{2}{\left\Vert{\Delta
 u}\right\Vert}_{L^{r}(B_{x}(R))}.$
\end{lem}
\quad Proof.\ \par 
We transcribe the problem in ${\mathbb{R}}^{n}$ by use of a coordinates
 path $(V,\varphi ).$ The Hodge laplacian is the second order
 elliptic matrix operator $\Delta _{\varphi }$ with ${\mathcal{C}}^{\infty
 }$ coefficients operating in $\varphi (V)\subset {\mathbb{R}}^{n}.$
 By the choice of a $\epsilon $ admissible ball $B_{x}(R),$ and
 with $R':=3R/4,$ we have:\ \par 
\quad \quad \quad $\displaystyle U':=\varphi (B_{x}(R'))\subset B_{e}(0,(1+\epsilon
 )R'),\ U:=\varphi (B_{x}(R))\subset B_{e}(0,(1+\epsilon )R)\subset
 \varphi (V).$\ \par 
We apply lemma~\ref{lE0}, to the euclidean balls $B':=B_{e}(0,(1+\epsilon
 )R'),\ B:=B_{e}(0,(1+\epsilon )R)$ and we get, with $u_{\varphi
 }$ the $p$ form $u$ read in the chart $(V,\varphi )$\ \par 
\quad \quad \quad $\displaystyle \ {\left\Vert{u_{\varphi }}\right\Vert}_{W^{2,r}(B')}\leq
 c_{1}R^{-2}{\left\Vert{u_{\varphi }}\right\Vert}_{L^{r}(B)}+c_{2}{\left\Vert{\Delta
 _{\varphi }u_{\varphi }}\right\Vert}_{L^{r}(B)}.$\ \par 
The fact that the coefficients of $\Delta _{\varphi }$ are $\epsilon
 $ near, in the ${\mathcal{C}}^{1}$ norm, of those of $\Delta
 _{{\mathbb{R}}},$ by condition 2) in the definition of the $\epsilon
 $ admissible ball, implies that the constants $c_{1},c_{2}$
 depend only on $n,r$ and $\epsilon .$\ \par 
\quad The Lebesgue measure on $U$ and the canonical measure $dv_{g}$
 on $B_{x}(R)$ are equivalent; precisely because of condition
 1) we get that $(1-\epsilon )^{n}\leq \left\vert{\mathrm{d}\mathrm{e}\mathrm{t}g}\right\vert
 \leq (1+\epsilon )^{n},$ and the measure $dv_{g}$ read in the
 chart $\varphi $ is $dv_{g}={\sqrt{\left\vert{\mathrm{d}\mathrm{e}\mathrm{t}g_{ij}}\right\vert
 }}d\xi ,$ where $d\xi $ is the Lebesgue measure in ${\mathbb{R}}^{n}.$
 So the Lebesgue estimates and the Sobolev estimates up to order
 2 on $U$ are valid in $B_{x}(R)$ up to a constant depending
 only on $n,r$ and $\displaystyle \epsilon $ by lemma~\ref{HCS43}.
 In particular:\ \par 
\quad \quad \quad \begin{equation}  \forall x\in M,\ \mathrm{V}\mathrm{o}\mathrm{l}(B_{x}(R))\leq
 (1+\epsilon )^{n/2}\nu _{n}R^{n},\label{HC23}\end{equation}\ \par 
where $\nu _{n}$ is the euclidean volume of the unit ball in
 ${\mathbb{R}}^{n}.$\ \par 
\quad So passing back to $M,$ we get, with $A:=B_{x}((1+2\epsilon )R)\supset
 \varphi ^{-1}(B),\ A':=\varphi ^{-1}(B')$\ \par 
\quad \quad \quad $\displaystyle {\left\Vert{u}\right\Vert}_{W^{2,r}(A')}\leq c_{1}R^{-2}{\left\Vert{u}\right\Vert}_{L^{r}(A)}+c_{2}{\left\Vert{\Delta
 u}\right\Vert}_{L^{r}(A)}).$\ \par 
Now we notice that $A'\supset B_{x}(R')$ so a fortiori:\ \par 
\quad \quad \quad $\displaystyle {\left\Vert{u}\right\Vert}_{W^{2,r}(B_{x}(R'))}\leq
 c_{1}R^{-2}{\left\Vert{u}\right\Vert}_{L^{r}(B_{x}((1+2\epsilon
 )R)}+c_{2}{\left\Vert{\Delta u}\right\Vert}_{L^{r}(B_{x}((1+2\epsilon
 )R)}.$\ \par 
Finally, choosing $\epsilon \leq 1/4$ we get:\ \par 
\quad \quad \quad $\displaystyle {\left\Vert{u}\right\Vert}_{W^{2,r}(B_{x}(R'))}\leq
 c_{1}R^{-2}{\left\Vert{u}\right\Vert}_{L^{r}(B_{x}(3R/2)}+c_{2}{\left\Vert{\Delta
 u}\right\Vert}_{L^{r}(B_{x}(3R/2)},$\ \par 
i.e. the $CZI$ local interior inequalities on $B_{x}(3R/4)\subset
 B_{x}(3R/2)\subset M.$ So changing $R$ for $3R/2$ we proved
 the lemma. $\blacksquare $\ \par 

\section{The raising steps method.}
\quad Let $(M,g)$ be a riemannian manifold. From now on we take $\epsilon
 =\epsilon _{0}$ with $\epsilon _{0}$ given by lemma~\ref{CF1}
 and we take the $\epsilon _{0}$ admissible  radius and the Vitali
 covering $\displaystyle \lbrace B(x_{j},5r(x_{j}))\rbrace _{j\in
 {\mathbb{N}}}$ associated to it.\ \par 

\begin{defin}
~\label{CL25}Let $(M,g)$ be a riemannian manifold. A weight relative
 to the covering $\displaystyle \lbrace B(x_{j},5r(x_{j}))\rbrace
 _{j\in {\mathbb{N}}}$ is a function $\displaystyle w(x)>0$ on
 $M$ such that :\par 
\quad there are two constants $\displaystyle 0<c_{iw}\leq 1\leq c_{sw}$
 such that, setting\par 
\quad \quad \quad $\displaystyle \forall j\in {\mathbb{N}},\ B_{j}:=B(x_{j},5r(x_{j})),\
 w_{j}:=\frac{1}{v_{g}(B_{j})}\int_{B_{j}}{w(x)dv_{g}(x)},$\par 
we have $\displaystyle \forall j\in {\mathbb{N}},\ \forall x\in
 B_{j},\ c_{iw}w_{j}\leq w(x)\leq c_{sw}w_{j}.$ By smoothing
 $w$ if necessary, we shall also suppose that $w\in {\mathcal{C}}^{\infty
 }(M).$\par 
\quad As an example we have the constant weight, $\displaystyle \forall
 x\in M,\ w(x)=1.$
\end{defin}
This means that $w$ varies slowly on $\displaystyle B_{j}.$\ \par 
\quad So let $\displaystyle w(x)>0$ be any weight we say that $\displaystyle
 \omega \in L^{r}_{p}(M,w),$ if :\ \par 
\quad \quad \quad $\displaystyle {\left\Vert{\omega }\right\Vert}_{L^{r}(M,w)}^{r}:=\int_{M}{\left\vert{\omega
 (x)}\right\vert ^{r}w(x)dv_{g}(x)}<\infty .$\ \par 

\subsection{The raising steps method.}
\quad We shall use the following lemma.\ \par 

\begin{lem}
~\label{HC34} For $\chi \in {\mathcal{D}}(M)$ and $\displaystyle
 u\in W^{2,r}_{p}(M),$ set $\displaystyle B(\chi ,u):=\Delta
 (\chi u)-\chi \Delta (u).$ We have :\par 
\quad \quad \quad $\displaystyle \ \left\vert{B(\chi ,u)}\right\vert \leq \left\vert{\Delta
 \chi }\right\vert \left\vert{u}\right\vert +2\left\vert{\nabla
 \chi }\right\vert \left\vert{\nabla u}\right\vert .$
\end{lem}
\quad Proof.\ \par 
Exactly as for Proposition G.III.6 in~\cite{BerGauMaz71} we have
 in an exponential chart at a point $\displaystyle x\in M,$\ \par 
\quad \quad \quad $\displaystyle u=\sum_{J,\left\vert{J}\right\vert =p}{u_{J}dx^{J},}\
 g^{ij}(x)=\delta _{ij}$ and the basis $\displaystyle \lbrace
 \frac{\partial }{\partial x_{j}}\rbrace _{j=1,...,n}$ is orthogonal.\ \par 
In this chart and \emph{at the point} $x$ we have that the laplacian
 is diagonal so\ \par 
\quad \quad \quad $\displaystyle \Delta u(x)=\sum_{J,\left\vert{J}\right\vert =p}{\frac{\partial
 ^{2}u_{J}}{\partial x_{j}^{2}}(x)dx^{J}}$\ \par 
hence, for any $\displaystyle x\in M,$\ \par 
\quad \quad \quad $\displaystyle B(\chi ,u)(x)=\Delta \chi (x)u(x)-2\sum_{J,\ \left\vert{J}\right\vert
 =p}{(\sum_{j=1}^{n}{\frac{\partial u_{J}}{\partial x_{j}}\frac{\partial
 \chi }{\partial x_{j}}})dx^{J}}.$\ \par 
So we get\ \par 
\quad \quad \quad $\displaystyle \left\vert{B(\chi ,u)}\right\vert \leq \left\vert{\Delta
 \chi }\right\vert \left\vert{u}\right\vert +2\left\vert{\nabla
 \chi }\right\vert \left\vert{\nabla u}\right\vert .$ $\blacksquare $\ \par 

\begin{lem}
~\label{5S4}Let $w$ be a weight relative to the covering ${\mathcal{C}}_{\epsilon
 }$ and set $\displaystyle w_{j}$ as in definition~\ref{CL25}.
 If  $v:=\sum_{j\in {\mathbb{N}}}{\chi _{j}u_{j}}$ then we have\par 
(i)               $\displaystyle {\left\Vert{v}\right\Vert}_{L^{s}(M,w^{s})}^{s}\leq
 T^{s}c_{sw}^{s}\sum_{j\in {\mathbb{N}}}{w_{j}^{s}{\left\Vert{u_{j}}\right\Vert}_{L^{s}(B_{j})}^{s}}.$\par
 
(ii)              $\displaystyle {\left\Vert{\nabla v}\right\Vert}_{L^{s}(M,w^{s})}^{s}\leq
 2^{s/s'}(1+C\epsilon )T^{s}c_{sw}^{s}\sum_{j\in {\mathbb{N}}}{w_{j}^{s}(R_{j}^{-s}{\left\Vert{u_{j}}\right\Vert}_{L^{s}(B_{j})}^{s}+{\left\Vert{\nabla
 u_{j}}\right\Vert}_{L^{s}(B_{j})}^{s})}.$\par 
(iii)     $\displaystyle {\left\Vert{\nabla ^{2}v}\right\Vert}_{L^{s}(M,w^{s})}^{s}\leq
 3^{s/s'}(1+C\epsilon )T^{s}c_{sw}^{s}\sum_{j\in {\mathbb{N}}}{w_{j}^{s}(R_{j}^{-2s}{\left\Vert{u_{j}}\right\Vert}_{L^{s}(B_{j})}^{s}+{\left\Vert{\nabla
 ^{2}u_{j}}\right\Vert}_{L^{s}(B_{j})}^{s}+R_{j}^{-s}{\left\Vert{\nabla
 u_{j}}\right\Vert}_{L^{s}(B_{j})}^{s})}.$
\end{lem}
\quad Proof.\ \par 
We have for \emph{(i)}\ \par 
\quad \quad \quad ${\left\Vert{v}\right\Vert}_{L^{s}(M,w^{s})}^{s}=\int_{M}{\left\vert{v}\right\vert
 ^{s}w^{s}dv_{g}}\leq \sum_{k\in {\mathbb{N}}}{\int_{B_{k}}{\left\vert{\sum_{j\in
 {\mathbb{N}}}{\chi _{j}u_{j}}}\right\vert ^{s}w^{s}dv_{g}}}.$\ \par 
But the support of $\displaystyle \chi _{j}$ is in $\displaystyle
 B_{j}$ and the overlap of the covering is less that $T$ so let\ \par 
\quad \quad \quad $\displaystyle I(k):=\lbrace B_{j}::B_{j}\cap B_{k}\neq \emptyset
 \rbrace $\ \par 
then $\displaystyle \mathrm{C}\mathrm{a}\mathrm{r}\mathrm{d}I(k)\leq
 T$ and we have\ \par 
\quad \quad \quad $\displaystyle {\left\Vert{v}\right\Vert}_{L^{s}(M,w^{s})}^{s}\leq
 \sum_{k\in {\mathbb{N}}}{\int_{B_{k}}{\left\vert{\sum_{j\in
 I(k)}{\chi _{j}u_{j}}}\right\vert ^{s}w^{s}dv_{g}}}.$\ \par 
We have, comparing the $l^{1}$and $l^{s}$ norms by H\"older inequalities,\ \par 
\quad \quad \quad $\displaystyle \left\vert{\sum_{j\in I(k)}{\chi _{j}u_{j}}}\right\vert
 ^{s}\leq T^{s-1}\sum_{j\in I(k)}{\left\vert{\chi _{j}u_{j}}\right\vert
 ^{s}}$\ \par 
so\ \par 
\quad \quad \quad \begin{equation} {\left\Vert{ v}\right\Vert}_{L^{s}(M,w^{s})}^{s}\leq
 T^{s-1}\sum_{k\in {\mathbb{N}}}{\sum_{j\in I(k)}{\int_{B_{k}}{\left\vert{\chi
 _{j}u_{j}}\right\vert ^{s}w^{s}dv_{g}}}}.\label{5S0}\end{equation}\ \par 
We still have, because $\displaystyle \chi _{j}$ is supported
 by $\displaystyle B_{j},$\ \par 
\quad \quad \quad $\displaystyle \sum_{j\in I(k)}{\int_{B_{k}}{\left\vert{\chi
 _{j}u_{j}}\right\vert ^{s}w^{s}dv_{g}}}=\sum_{j\in {\mathbb{N}}}{\int_{B_{k}}{\left\vert{\chi
 _{j}u_{j}}\right\vert ^{s}w^{s}dv_{g}}}$\ \par 
hence, exchanging the order of summation, all terms being positive,\ \par 
\quad \quad \quad $\displaystyle {\left\Vert{v}\right\Vert}_{L^{s}(M,w^{s})}^{s}\leq
 T^{s-1}\sum_{j\in {\mathbb{N}}}{\sum_{k\in {\mathbb{N}}}{\int_{B_{k}}{\left\vert{\chi
 _{j}u_{j}}\right\vert ^{s}w^{s}dv_{g}}}}.$\ \par 
The overlap being less than $T$ we get\ \par 
\quad \quad \quad $\displaystyle \sum_{k\in {\mathbb{N}}}{\int_{B_{k}}{\left\vert{\chi
 _{j}u_{j}}\right\vert ^{s}w^{s}dv_{g}}}\leq T\int_{M}{\left\vert{\chi
 _{j}u_{j}}\right\vert ^{s}w^{s}dv_{g}}$\ \par 
so\ \par 
\quad \quad \quad $\displaystyle {\left\Vert{v}\right\Vert}_{L^{s}(M,w^{s})}^{s}\leq
 T^{s-1}T\sum_{j\in {\mathbb{N}}}{\int_{M}{\left\vert{\chi _{j}u_{j}}\right\vert
 ^{s}w^{s}dv_{g}}}=T^{s}\sum_{j\in {\mathbb{N}}}{\int_{B_{j}}{\left\vert{\chi
 _{j}u_{j}}\right\vert ^{s}w^{s}dv_{g}}}.$\ \par 
With the constants $\displaystyle c_{sw}$ defined in definition~\ref{CL25},\
 \par 
\quad \quad \quad $\displaystyle {\left\Vert{v}\right\Vert}_{L^{s}(M,w^{s})}^{s}\leq
 T^{s}c_{sw}^{s}\sum_{j\in {\mathbb{N}}}{w_{j}^{s}{\left\Vert{\chi
 _{j}u_{j}}\right\Vert}_{L^{s}(B_{j})}^{s}}\leq T^{s}c_{sw}^{s}\sum_{j\in
 {\mathbb{N}}}{w_{j}^{s}{\left\Vert{u_{j}}\right\Vert}_{L^{s}(B_{j})}^{s}}$\
 \par 
hence we get the \emph{(i)} :\ \par 
\quad \quad \quad $\displaystyle {\left\Vert{v}\right\Vert}_{L^{s}(M,w^{s})}^{s}\leq
 T^{s}c_{sw}^{s}\sum_{j\in {\mathbb{N}}}{w_{j}^{s}{\left\Vert{u_{j}}\right\Vert}_{L^{s}(B_{j})}^{s}}.$\
 \par 
\ \par 
\quad For \emph{(ii)}.\ \par 
\quad \quad \quad $v:=\sum_{j\in {\mathbb{N}}}{\chi _{j}u_{j}}\Rightarrow \left\vert{\nabla
 v}\right\vert \leq (1+C\epsilon )\sum_{j\in {\mathbb{N}}}{(\left\vert{\chi
 _{j}}\right\vert \left\vert{\nabla u_{j}}\right\vert +\left\vert{\nabla
 \chi _{j}}\right\vert \left\vert{u_{j}}\right\vert )}$\ \par 
by lemma~\ref{3S4}.\ \par 
Because $\lbrace \chi _{j}\rbrace _{j\in {\mathbb{N}}}$ is a
 partition of unity relative to the covering $\lbrace B_{j}\rbrace
 _{j\in {\mathbb{N}}},$ we have\ \par 
\quad \quad \quad $\displaystyle \left\vert{\nabla \chi _{j}}\right\vert \leq \frac{1}{R_{j}}\
 ;\ \left\vert{\nabla ^{2}\chi _{j}}\right\vert \leq \frac{1}{R_{j}^{2}}.$\
 \par 
Hence for the first term, $\displaystyle A:=\sum_{j\in {\mathbb{N}}}{\left\vert{\chi
 _{j}}\right\vert \left\vert{\nabla u_{j}}\right\vert }$ we get,
 again exactly as above\ \par 
\quad \quad \quad $\displaystyle {\left\Vert{A}\right\Vert}_{L^{s}(M,w^{s})}^{s}\leq
 T^{s}c_{sw}^{s}\sum_{j\in {\mathbb{N}}}{w_{j}^{s}{\left\Vert{\nabla
 u_{j}}\right\Vert}_{L^{s}(B_{j})}^{s}}.$\ \par 
For the second one, $\displaystyle B:=\sum_{j\in {\mathbb{N}}}{\left\vert{\nabla
 \chi _{j}}\right\vert \left\vert{u_{j}}\right\vert }$ we get
 also as above, using the estimate $\displaystyle \ \left\vert{\nabla
 \chi _{j}}\right\vert \leq \frac{1}{R_{j}},$\ \par 
\quad \quad \quad $\displaystyle {\left\Vert{B}\right\Vert}_{L^{s}(M,w^{s})}^{s}\leq
 T^{s}c_{sw}^{s}\sum_{j\in {\mathbb{N}}}{R_{j}^{-s}w_{j}^{s}{\left\Vert{u_{j}}\right\Vert}_{L^{s}(B_{j})}^{s}}.$\
 \par 
Because $\displaystyle (a+b)^{s}\leq 2^{s/s'}(a^{s}+b^{s})$ we get\ \par 
\quad \quad \quad $\displaystyle {\left\Vert{\nabla v}\right\Vert}_{L^{s}(M,w^{s})}^{s}\leq
 2^{s/s'}(1+C\epsilon )T^{s}c_{sw}^{s}\sum_{j\in {\mathbb{N}}}{w_{j}^{s}(R_{j}^{-s}{\left\Vert{u_{j}}\right\Vert}_{L^{s}(B_{j})}^{s}+{\left\Vert{\nabla
 u_{j}}\right\Vert}_{L^{s}(B_{j})}^{s})}.$\ \par 
\ \par 
Finally for \emph{(iii)}. By lemma~\ref{3S4} \emph{(ii)}, we get\ \par 
\quad \quad \quad $\left\vert{\nabla ^{2}(v)}\right\vert \leq (1+C\epsilon )\sum_{j\in
 {\mathbb{N}}}{(\left\vert{\chi _{j}}\right\vert \left\vert{\nabla
 ^{2}u_{j}}\right\vert +\left\vert{\nabla ^{2}\chi _{j}}\right\vert
 \left\vert{u_{j}}\right\vert +\left\vert{\nabla \chi _{j}}\right\vert
 \left\vert{\nabla u_{j}}\right\vert )}$\ \par 
So we get, for the two first terms, as above\ \par 
\quad \quad $\displaystyle C:=(1+C\epsilon )\sum_{j\in {\mathbb{N}}}{(\left\vert{\chi
 _{j}}\right\vert \left\vert{\nabla ^{2}u_{j}}\right\vert )}\Rightarrow
 {\left\Vert{D}\right\Vert}_{L^{s}(M,w^{s})}^{s}\leq (1+C\epsilon
 )T^{s}c_{sw}^{s}\sum_{j\in {\mathbb{N}}}{w_{j}^{s}{\left\Vert{\nabla
 ^{2}u_{j}}\right\Vert}_{L^{s}(B_{j})}^{s}}.$\ \par 
And using the estimate $\displaystyle \ \left\vert{\nabla ^{2}\chi
 _{j}}\right\vert \leq \frac{1}{R_{j}^{2}},$\ \par 
\quad \quad \quad $\displaystyle D:=(1+C\epsilon )\sum_{j\in {\mathbb{N}}}{(\left\vert{\nabla
 ^{2}\chi _{j}}\right\vert \left\vert{u_{j}}\right\vert )}\Rightarrow
 {\left\Vert{C}\right\Vert}_{L^{s}(M,w^{s})}^{s}\leq (1+C\epsilon
 )T^{s}c_{sw}^{s}\sum_{j\in {\mathbb{N}}}{R_{j}^{-2s}w_{j}^{s}{\left\Vert{u_{j}}\right\Vert}_{L^{s}(B_{j})}^{s}}.$\
 \par 
For the third one, we get using the estimate $\displaystyle \
 \left\vert{\nabla \chi _{j}}\right\vert \leq \frac{1}{R_{j}},$\ \par 
\quad \quad \quad $\displaystyle E:=(1+C\epsilon )\sum_{j\in {\mathbb{N}}}{(\left\vert{\nabla
 \chi _{j}}\right\vert \left\vert{\nabla u_{j}}\right\vert )}\Rightarrow
 {\left\Vert{D}\right\Vert}_{L^{s}(M,w^{s})}^{s}\leq (1+C\epsilon
 )T^{s}c_{sw}^{s}\sum_{j\in {\mathbb{N}}}{w_{j}^{s}R_{j}^{-s}{\left\Vert{\nabla
 u_{j}}\right\Vert}_{L^{s}(B_{j})}^{s}}.$\ \par 
Adding this, we get\ \par 
$\displaystyle {\left\Vert{\nabla ^{2}v}\right\Vert}_{L^{s}(M,w^{s})}^{s}\leq
 3^{s/s'}(1+C\epsilon )T^{s}c_{sw}^{s}\sum_{j\in {\mathbb{N}}}{w_{j}^{s}(R_{j}^{-2s}{\left\Vert{u_{j}}\right\Vert}_{L^{s}(B_{j})}^{s}+{\left\Vert{\nabla
 ^{2}u_{j}}\right\Vert}_{L^{s}(B_{j})}^{s}+R_{j}^{-s}{\left\Vert{\nabla
 u_{j}}\right\Vert}_{L^{s}(B_{j})}^{s})}.$ $\blacksquare $\ \par 

\begin{lem}
~\label{5S6} Let $w$ be a weight relative to the covering ${\mathcal{C}}_{\epsilon
 }$ and set $\displaystyle w_{j}$ as in definition~\ref{CL25}. Suppose that\par 
\quad \quad \quad $\displaystyle \ I^{s}=\sum_{j\in {\mathbb{N}}}{w_{j}^{s}{\left\Vert{u_{j}}\right\Vert}_{L^{s}(B_{j})}^{s}},$\par
 
and, with $\displaystyle s\geq r,$\par 
\quad \quad \quad $\displaystyle w_{j}{\left\Vert{\chi _{j}u_{j}}\right\Vert}_{L^{s}(B_{j})}\leq
 w_{j}R_{j}^{-\gamma }c{\left\Vert{\omega }\right\Vert}_{L^{r}(B_{j})},$\par 
Then we have, with $\displaystyle \forall x\in M,\ \tilde w(x):=R(x)^{-\gamma
 }w(x),$\par 
\quad \quad \quad $\displaystyle \ I\leq c_{w}T^{s/r}{\left\Vert{\omega }\right\Vert}_{L^{r}(M,\tilde
 w^{r})}.$
\end{lem}
\quad Proof.\ \par 
By\ \par 
\quad \quad \quad $\displaystyle \sum_{j\in {\mathbb{N}}}{a_{j}^{s}}\leq (\sum_{j\in
 {\mathbb{N}}}{a_{j}^{r}})^{s/r}$ because $\displaystyle s\geq r,$\ \par 
we get\ \par 
\quad \quad \quad $\displaystyle \ I^{s}\leq {\left({\sum_{j\in {\mathbb{N}}}{w_{j}^{r}R_{j}^{-\gamma
 r}{\left\Vert{\omega }\right\Vert}^{r}_{L^{r}(B_{j})}}}\right)}^{s/r}.$\ \par 
By lemma~\ref{CF10} we have\ \par 
\quad \quad \quad $\displaystyle \forall x\in B_{j},\ d(x,x_{j})<R_{j}=5r(x_{j})\leq
 \frac{1}{4}R(x_{j})\leq \frac{1}{4}(R(x_{j})+R(x))\Rightarrow
 R(x)\leq 4R(x_{j}),$\ \par 
hence, because $\displaystyle r(x_{j})=\frac{R(x_{j})}{120}$
 and $\displaystyle R_{j}=5r(x_{j})=\frac{R(x_{j})}{24},$\ \par 
\quad \quad \quad \begin{equation}  \forall x\in B_{j},\ R(x)\leq 4R_{j}\Rightarrow
 R_{j}^{-\gamma r}\leq 96^{\gamma r}R(x)^{-\gamma r}.\label{5S2}\end{equation}\
 \par 
\quad But, by definition~\ref{CL25}, we have\ \par 
\quad \quad \quad $\displaystyle R_{j}^{-2r}w_{j}^{r}{\left\Vert{\omega }\right\Vert}^{r}_{L^{r}(B_{j})}\leq
 c_{iw}^{-r}R_{j}^{-\gamma r}\int_{B_{j}}{\left\vert{\omega }\right\vert
 ^{r}w^{r}dv_{g}}.$\ \par 
So\ \par 
\quad \quad \quad $\displaystyle \sum_{j\in {\mathbb{N}}}{R_{j}^{-\gamma r}w_{j}^{r}{\left\Vert{\omega
 }\right\Vert}^{r}_{L^{r}(B_{j})}}\leq c_{iw}^{-r}\sum_{j\in
 {\mathbb{N}}}{R_{j}^{-\gamma r}\int_{B_{j}}{\left\vert{\omega
 }\right\vert ^{r}w^{r}dv_{g}}}$\ \par 
and, by~(\ref{5S2}), we get\ \par 
\quad \quad \quad $\displaystyle \sum_{j\in {\mathbb{N}}}{R_{j}^{-\gamma r}w_{j}^{r}{\left\Vert{\omega
 }\right\Vert}^{r}_{L^{r}(B_{j})}}\leq 96^{\gamma r}c_{iw}^{-r}\sum_{j\in
 {\mathbb{N}}}{\int_{B_{j}}{\left\vert{\omega }\right\vert ^{r}R(x)^{-\gamma
 r}w^{r}dv_{g}}}.$\ \par 
\quad Set $\displaystyle \forall x\in M,\ \tilde w(x):=R(x)^{-\gamma }w(x).$\ \par 
Now, because the overlap is less that $T,$ by proposition~\ref{CF2},
 we get\ \par 
\quad \quad \quad $\displaystyle \sum_{j\in {\mathbb{N}}}{\int_{B_{j}}{\left\vert{\omega
 }\right\vert ^{r}\tilde w^{r}dv_{g}}}\leq 96^{\gamma r}T\int_{M}{\left\vert{\omega
 (x)}\right\vert ^{r}\tilde w(x)^{r}dv_{g}(x)}=T{\left\Vert{\omega
 }\right\Vert}_{L^{r}(M,\tilde w^{r})}^{r}.$\ \par 
\quad Putting this in $\displaystyle v,$ we get\ \par 
\quad \quad \quad $\displaystyle \ I^{s}\leq {\left({\sum_{j\in {\mathbb{N}}}{R_{j}^{-\gamma
 r}w_{j}^{r}{\left\Vert{\omega }\right\Vert}^{r}_{L^{r}(B_{j})}}}\right)}^{s/r}\leq
 96^{\gamma s}c_{iw}^{-s}C^{s}(T{\left\Vert{\omega }\right\Vert}_{L^{r}(M,\tilde
 w^{r})}^{r})^{s/r},$\ \par 
so, setting $\displaystyle c_{w}:=96^{\gamma }c_{iw}^{-1}C$ we get\ \par 
\quad \quad \quad $\displaystyle \ I\leq c_{w}T^{s/r}{\left\Vert{\omega }\right\Vert}_{L^{r}(M,\tilde
 w^{r})}.$ $\blacksquare $\ \par 
\ \par 
\quad With $\displaystyle R(x)$ the $\displaystyle \epsilon _{0}$ admissible
  radius at the point $\displaystyle x\in M,$ and ${\mathcal{C}}_{\epsilon
 _{0}}$ the $\epsilon _{0}$ admissible covering of $M,$ defined
 in section~\ref{CL24}, we shall prove now :\ \par 

\begin{thm}
~\label{CF3}(Raising Steps Method) Let $\displaystyle (M,g)$
 be a riemannian manifold and take $w$ a weight relative to the
 Vitali covering $\displaystyle \lbrace B(x_{j},5r(x_{j}))\rbrace
 _{j\in {\mathbb{N}}}.$\par 
For any $\displaystyle r\leq 2,$ any threshold $\displaystyle
 s\geq r,$ take $\displaystyle k\in {\mathbb{N}}$ such that $\displaystyle
 t_{k}:=S_{k}(r)\geq s$ then, with $\displaystyle w_{0}(x):=w(x)R(x)^{-2k},$\par
 
\quad $\displaystyle \forall \omega \in L^{r}_{p}(M,w_{0}^{r}),\ \exists
 v\in L^{r}_{p}(M,w^{r})\cap L^{s_{1}}_{p}(M,w^{s_{1}})\cap W^{2,r}(M,w^{r}),\
 \exists \tilde \omega \in L^{s}_{p}(M,w^{s})::\Delta v=\omega
 +\tilde \omega $\par 
with $\displaystyle s_{1}=S_{2}(r)$ and we have the control of the norms :\par 
\quad \quad \quad $\displaystyle \forall q\in \lbrack r,s_{1}\rbrack ,\ {\left\Vert{v}\right\Vert}_{L^{q}_{p}(M,w^{q})}\leq
 C_{q}{\left\Vert{\omega }\right\Vert}_{L^{r}_{p}(M,w_{0}^{r})}\
 ;\ {\left\Vert{v}\right\Vert}_{W^{2,r}_{p}(M,w^{r})}\leq C_{r}{\left\Vert{\omega
 }\right\Vert}_{L^{r}_{p}(M,w_{0}^{r})}\ ;$\par 
\quad \quad \quad \quad \quad \quad \quad \quad \quad \quad \quad \quad $\displaystyle {\left\Vert{\tilde \omega }\right\Vert}_{L^{s}_{p}(M,w^{s})}\leq
 C_{s}{\left\Vert{\omega }\right\Vert}_{L^{r}_{p}(M,w_{0}^{r})}.$\par 
Moreover $v$ and $\tilde \omega $ are linear in $\omega .$\par 
If $M$ is complete and $\omega $ is of compact support, so are
 $v$ and $\tilde \omega .$
\end{thm}
\quad Proof.\ \par 
To simplify notations we do not put the $p$ referring to the
 degree of the forms, i.e. we shall write $\displaystyle L^{r}$
 instead of $\displaystyle L^{r}_{p},\ W^{2,r}$ instead of $\displaystyle
 W^{2,r}_{p},$ etc...\ \par 
\quad Set $\displaystyle R_{j}:=5r(x_{j}),\ B_{j}:=B(x_{j},R_{j})$
 and apply lemma~\ref{CF1} to get, with $\displaystyle c=c(n,r,\epsilon
 _{0}),$\ \par 
\quad \quad \quad \begin{equation}  \exists u_{j}\in W^{2,r}(B_{j})::\Delta u_{j}=\omega
 ,\ {\left\Vert{u_{j}}\right\Vert}_{W^{2,r}(B_{j})}\leq c{\left\Vert{\omega
 }\right\Vert}_{L^{r}(B_{j})},\label{5S7}\end{equation}\ \par 
with $\displaystyle u_{j}$ linear in $\omega _{\mid B_{j}}.$\ \par 
So by lemma~\ref{3S2} we get, with $\displaystyle t=S_{1}(r),\
 s=S_{2}(r),$\ \par 
\quad \quad \quad $\displaystyle u_{j}\in L^{s}(B_{j}),\ {\left\Vert{u_{j}}\right\Vert}_{L^{s}(B_{j})}\leq
 CR_{j}^{-2}\ {\left\Vert{u}\right\Vert}_{W^{2,r}(B(x,R))}\leq
 cCR_{j}^{-2}{\left\Vert{\omega }\right\Vert}_{L^{r}(B_{j})}$\ \par 
and\ \par 
\quad \quad \quad $\displaystyle \nabla u_{j}\in L^{t}(B_{j}),\ {\left\Vert{\nabla
 u_{j}}\right\Vert}_{L^{t}(B_{j})}\leq CR_{j}^{-1}{\left\Vert{u}\right\Vert}_{W^{2,r}(B(x,R))}\leq
 CcR_{j}^{-1}{\left\Vert{\omega }\right\Vert}_{L^{r}(B_{j})}.$\ \par 
Hence, because $\displaystyle u_{j}\in L^{r}(B_{j}),$ we have
 by interpolation~\cite{BerghLofstrom76}, that $\displaystyle
 \forall s'\in \lbrack r,s\rbrack ,\ u_{j}\in L^{s'}(B_{j})$
 with $\displaystyle {\left\Vert{u_{j}}\right\Vert}_{L^{s'}(B_{j})}\leq
 cCR_{j}^{-2}{\left\Vert{\omega }\right\Vert}_{L^{r}(B_{j})}.$\ \par 
The same way, because $\displaystyle \nabla u_{j}\in L^{r}(B_{j}),$
 by interpolation we get $\displaystyle \forall t'\in \lbrack
 r,t\rbrack ,\ \nabla u_{j}\in L^{t'}(B_{j})$ with\ \par 
\quad \quad \quad \begin{equation} {\left\Vert{ \nabla u_{j}}\right\Vert}_{L^{t'}(B_{j})}\leq
 cCR_{j}^{-1}{\left\Vert{\omega }\right\Vert}_{L^{r}(B_{j})}.\label{5S1}\end{equation}\
 \par 
\quad Let $\displaystyle \lbrace \chi _{j}\rbrace _{j\in {\mathbb{N}}}$
 be a partition of unity associated to the covering $\displaystyle
 \lbrace B(x_{j},R_{j})\rbrace _{j\in {\mathbb{N}}}$ then we set\ \par 
\quad \quad \quad $\displaystyle v_{0}:=\sum_{j\in {\mathbb{N}}}{\chi _{j}u_{j}.}$\ \par 
Because the $\displaystyle u_{j}$ are linear in $\omega _{\mid
 B_{j}},\ v_{0}$ is linear in $\omega .$\ \par 
We have, because $\displaystyle {\left\Vert{\chi _{j}}\right\Vert}_{\infty
 }=1,$\ \par 
\quad \quad \quad $\displaystyle {\left\Vert{\chi _{j}u_{j}}\right\Vert}_{L^{s}(B_{j})}\leq
 cCR_{j}^{-2}{\left\Vert{\omega }\right\Vert}_{L^{r}(B_{j})},$\ \par 
and multiplying by the $\displaystyle w_{j},$ given in definition~\ref{CL25},\
 \par 
\quad \quad \quad \begin{equation}  w_{j}{\left\Vert{\chi _{j}u_{j}}\right\Vert}_{L^{s}(B_{j})}\leq
 w_{j}R_{j}^{-2}c{\left\Vert{\omega }\right\Vert}_{L^{r}(B_{j})}.\label{5S5}\end{equation}\
 \par 
\quad By lemma~\ref{5S4} \emph{(i)}, we have\ \par 
\quad \quad \quad $\displaystyle {\left\Vert{v_{0}}\right\Vert}_{L^{s}(M,w^{s})}^{s}\leq
 T^{s}c_{sw}^{s}\sum_{j\in {\mathbb{N}}}{w_{j}^{s}{\left\Vert{u_{j}}\right\Vert}_{L^{s}(B_{j})}^{s}}.$\
 \par 
Now, because of~(\ref{5S5}), we can apply lemma~\ref{5S6} with
 $\displaystyle I={\left\Vert{v_{0}}\right\Vert}_{L^{s}(M,w^{s})}$
 and $\displaystyle \gamma =2\ ;$ we get, with $\displaystyle
 c_{w}:=c_{w}:=96^{2}c_{iw}^{-1}C,$ and $\displaystyle \tilde
 w_{2}(x):=R(x)^{-2}w(x).$\ \par 
\quad \quad \quad $\displaystyle {\left\Vert{v_{0}}\right\Vert}_{L^{s}(M,w^{s})}\leq
 c_{w}T^{s/r}{\left\Vert{\omega }\right\Vert}_{L^{r}(M,\tilde
 w^{_{2}r})}.$\ \par 
\quad We also have $\displaystyle v_{0}\in L^{r}(M,w^{r})$ because
 $\displaystyle u_{j}\in W^{2,r}(B_{j})\Rightarrow u_{j}\in L^{r}(B_{j})$
 as well, this means that $\displaystyle v_{0}w\in L^{r}(M)\cap
 L^{s}(M)$ hence by interpolation we have that\ \par 
\quad \quad \quad $\displaystyle v_{0}w\in L^{t'}(M)\Rightarrow v_{0}\in L^{t'}(M,w^{t'})$
 for any $\displaystyle t'\in \lbrack r,s\rbrack $ with the same
 control of the norms.\ \par 
\ \par 
\quad Because $\displaystyle u_{j}\in W^{2,r}(B_{j})$ we shall apply
 the same procedure to $\displaystyle \nabla v_{0}$ by use of
 lemma~\ref{5S4} \emph{(ii)}, with $\displaystyle s=r,\ v=v_{0},\ $we get\ \par 
\quad \quad \quad \begin{equation} {\left\Vert{ \nabla v}\right\Vert}_{L^{r}(M,w^{r})}^{r}\leq
 2^{r/r'}(1+C\epsilon )T^{r}c_{sw}^{r}\sum_{j\in {\mathbb{N}}}{w_{j}^{r}(R_{j}^{-r}{\left\Vert{u_{j}}\right\Vert}_{L^{r}(B_{j})}^{r}+{\left\Vert{\nabla
 u_{j}}\right\Vert}_{L^{r}(B_{j})}^{r})}.\label{5S8}\end{equation}\ \par 
But, by~(\ref{5S7}),\ \par 
\quad \quad \quad $\displaystyle {\left\Vert{u_{j}}\right\Vert}_{W^{2,r}(B_{j})}\leq
 c{\left\Vert{\omega }\right\Vert}_{L^{r}(B_{j})}\Rightarrow
 {\left\Vert{\nabla u_{j}}\right\Vert}_{L^{r}(B_{j})}\leq c{\left\Vert{\omega
 }\right\Vert}_{L^{r}(B_{j})}.$\ \par 
To the first term of~(\ref{5S8}), $\displaystyle A:=\sum_{j\in
 {\mathbb{N}}}{w_{j}^{s}R_{j}^{-s}{\left\Vert{u_{j}}\right\Vert}_{L^{s}(B_{j})}^{s}},$
 we can apply lemma~\ref{5S6} with $\displaystyle s=r,\ I\rightarrow
 A$ and $\displaystyle C\rightarrow Tc_{sw},\ \gamma =1\ ,\ \tilde
 w_{1}(x):=R(x)^{-1}w(x),$ to get\ \par 
\quad \quad \quad $\displaystyle {\left\Vert{A}\right\Vert}_{L^{r}(M,w^{r})}\leq
 c_{w}T{\left\Vert{\omega }\right\Vert}_{L^{r}(M,\tilde w_{1}^{r})}.$\ \par 
To the second term of~(\ref{5S8}), $\displaystyle B:=\sum_{j\in
 {\mathbb{N}}}{w_{j}^{r}{\left\Vert{\nabla u_{j}}\right\Vert}_{L^{r}(B_{j})}^{r}}$
 we can apply lemma~\ref{5S6} with $\displaystyle s=r,\ I\rightarrow
 B,\ u_{j}\rightarrow \nabla u_{j}$ and $\displaystyle C\rightarrow
 Tc_{sw},\ \gamma =0\ ,\ \tilde w(x):=w(x),$ to get\ \par 
\quad \quad \quad $\displaystyle {\left\Vert{v}\right\Vert}_{L^{r}(M,w^{r})}\leq
 c_{w}T{\left\Vert{\omega }\right\Vert}_{L^{r}(M,w^{r})}.$\ \par 
Adding these terms, we get\ \par 
\quad \quad \quad $\displaystyle {\left\Vert{\nabla v}\right\Vert}_{L^{r}(M,w^{r})}\leq
 2^{1/r'}(1+C\epsilon )^{1/r}Tc_{sw}c_{w}T({\left\Vert{\omega
 }\right\Vert}_{L^{r}(M,\tilde w_{1}^{r})}+{\left\Vert{\omega
 }\right\Vert}_{L^{r}(M,w^{r})}).$\ \par 
\ \par 
\quad Again because $\displaystyle u_{j}\in W^{2,r}(B_{j})$ we shall
 apply the same procedure to $\displaystyle \nabla ^{2}v_{0}$
 by use of lemma~\ref{5S4} \emph{(iii)}, with $\displaystyle
 s=r,\ v=v_{0},\ \tilde w_{2}(x):=R^{-2}(x)w(x),$ we get\ \par 
\quad $\displaystyle {\left\Vert{\nabla ^{2}v_{0}}\right\Vert}_{L^{r}(M,w^{r})}^{r}\leq
 3^{r/r'}(1+C\epsilon )T^{r}c_{sw}^{r}\sum_{j\in {\mathbb{N}}}{w_{j}^{r}(R_{j}^{-2r}{\left\Vert{u_{j}}\right\Vert}_{L^{r}(B_{j})}^{r}+{\left\Vert{\nabla
 ^{2}u_{j}}\right\Vert}_{L^{r}(B_{j})}^{r}+R_{j}^{-r}{\left\Vert{\nabla
 u_{j}}\right\Vert}_{L^{r}(B_{j})}^{r})}.$\ \par 
But, by~(\ref{5S7}),\ \par 
\quad \quad \quad $\displaystyle {\left\Vert{u_{j}}\right\Vert}_{W^{2,r}(B_{j})}\leq
 c{\left\Vert{\omega }\right\Vert}_{L^{r}(B_{j})}\Rightarrow
 {\left\Vert{\nabla u_{j}}\right\Vert}_{L^{r}(B_{j})}\leq c{\left\Vert{\omega
 }\right\Vert}_{L^{r}(B_{j})}$\ \par 
and\ \par 
\quad \quad \quad $\displaystyle {\left\Vert{\nabla ^{2}u_{j}}\right\Vert}_{L^{r}(B_{j})}\leq
 c{\left\Vert{\omega }\right\Vert}_{L^{r}(B_{j})}.$\ \par 
So playing the same game for each term, we get\ \par 
\quad \quad \quad $\displaystyle {\left\Vert{\nabla ^{2}v_{0}}\right\Vert}_{L^{r}(M,w^{r})}^{r}\leq
 3^{1/r'}(1+C\epsilon )^{1/r}Tc_{sw}c_{w}T({\left\Vert{\omega
 }\right\Vert}_{L^{r}(M,\tilde w^{_{2}r})}+{\left\Vert{\omega
 }\right\Vert}_{L^{r}(M,\tilde w_{1}^{r})}+{\left\Vert{\omega
 }\right\Vert}_{L^{r}(M,w^{r})}).$\ \par 
Because we always have $\displaystyle R(x)\leq 1,$ we get that
 $\displaystyle {\left\Vert{\omega }\right\Vert}_{L^{r}(M,\tilde
 w_{2}^{r})}\geq {\left\Vert{\omega }\right\Vert}_{L^{r}(M,\tilde
 w_{1}^{r})}\geq {\left\Vert{\omega }\right\Vert}_{L^{r}(M,w^{r})}$
 so finally\ \par 
\quad \quad \quad $\displaystyle {\left\Vert{v_{0}}\right\Vert}_{L^{r}(M,w^{r})}\leq
 C_{0}{\left\Vert{\omega }\right\Vert}_{L^{r}(M,w^{r})}\ ;\ $\ \par 
\quad \quad \quad $\displaystyle {\left\Vert{\nabla v_{0}}\right\Vert}_{L^{r}(M,w^{r})}\leq
 C_{1}{\left\Vert{\omega }\right\Vert}_{L^{r}(M,\tilde w_{1}^{r})}\ ;\ $\ \par 
\quad \quad \quad $\displaystyle {\left\Vert{\nabla ^{2}v_{0}}\right\Vert}_{L^{r}(M,w^{r})}^{r}\leq
 C_{2}{\left\Vert{\omega }\right\Vert}_{L^{r}(M,\tilde w_{2}^{r})}.$\ \par 
Because\ \par 
\quad \quad \quad $\displaystyle {\left\Vert{v_{0}}\right\Vert}_{W^{2,r}(M,w^{r})}:={\left\Vert{v_{0}}\right\Vert}_{L^{r}(M,w^{r})}+{\left\Vert{\nabla
 v_{0}}\right\Vert}_{L^{r}(M,w^{r})}+{\left\Vert{\nabla ^{2}v_{0}}\right\Vert}_{L^{r}(M,w^{r})}$\
 \par 
we get\ \par 
\quad \quad \quad $\displaystyle {\left\Vert{v_{0}}\right\Vert}_{W^{2,r}(M,w^{r})}\leq
 C{\left\Vert{\omega }\right\Vert}_{L^{r}(M,\tilde w_{2}^{r})},$\ \par 
where the constant $C$ depends only on $\displaystyle n,\ \epsilon
 ,\ T$ and the constants of the weight $w$ relative to the covering
 ${\mathcal{C}}_{\epsilon }.$\ \par 
\ \par 
\quad If $\omega $ is of compact support and if $M$ is complete, by
 lemma~\ref{HCN48} we can cover $\displaystyle \Supp \ \omega
 $ by a finite set $\displaystyle \lbrace B_{j}\rbrace _{j=1,...,N}$
 and then add a layer $\displaystyle \lbrace B_{j}\rbrace _{j=N_{0}+1,...,N_{1}}$
 not intersecting $\displaystyle \Supp \ \omega ,$ to cover $\displaystyle
 \partial K'$ where $\displaystyle K'$ is a compact containing
 $K.$ This means that we can cover $\displaystyle K'$ by a finite
 set $\displaystyle \lbrace B_{j}\rbrace _{j=1,...,N_{1}}.$ 
 By linearity we get $\displaystyle \forall j=N_{0}+1,...,N_{1},\
 \omega _{j}=0\Rightarrow u_{j}=0$ and setting now $\displaystyle
 v_{0}:=\sum_{j=1}^{N_{1}}{\chi _{j}u_{j}}$ we can extend $\displaystyle
 v_{0}$ as $0$ outside $\displaystyle \bigcup_{j=1}^{N_{1}}{B_{j}}$
 hence we get that $\displaystyle v_{0}$ is compactly supported.\ \par 
\quad We set, as in lemma~\ref{HC34}, $\displaystyle B(\chi _{j},u_{j})=\Delta
 (\chi _{j}u_{j})-\chi _{j}\Delta u_{j}.$ Now consider $\Delta
 v_{0},$ we get\ \par 
\quad \quad \quad $\displaystyle \Delta v_{0}=\sum_{j\in {\mathbb{N}}}{\Delta (\chi
 _{j}u_{j})}=\sum_{j\in {\mathbb{N}}}{\chi _{j}\Delta u_{j}}+\sum_{j\in
 {\mathbb{N}}}{B(\chi _{j},u_{j})}=\omega +\omega _{1},$\ \par 
with $\displaystyle \omega _{1}:=\sum_{j\in {\mathbb{N}}}{B(\chi
 _{j},u_{j})}.$\ \par 
Clearly $\Delta v_{0}$ is linear in $\omega $ so is $\omega _{1}.$\ \par 
The $\lbrace \chi _{j}\rbrace _{j\in {\mathbb{N}}}$ being a partition
 of unity relative to the covering $\lbrace B_{j}\rbrace _{j\in
 {\mathbb{N}}},$ we have $\displaystyle \left\vert{\nabla \chi
 _{j}}\right\vert \leq \frac{1}{R_{j}}$ and $\displaystyle \left\vert{\Delta
 \chi _{j}}\right\vert \leq \frac{1}{R_{j}^{2}}.$ We also have,
 because $\displaystyle {\left\Vert{u_{j}}\right\Vert}_{W^{2,r}(B_{j})}\leq
 c{\left\Vert{\omega }\right\Vert}_{L^{r}(B_{j})},$\ \par 
\quad \quad \quad $\displaystyle {\left\Vert{\nabla u_{j}}\right\Vert}_{L^{t}(B_{j})}\leq
 cR_{j}^{-1}{\left\Vert{\omega }\right\Vert}_{L^{r}(B_{j})}$\ \par 
by lemma~\ref{3S3} \emph{(ii)}, and\ \par 
\quad \quad \quad \begin{equation} {\left\Vert{ u_{j}}\right\Vert}_{L^{s}(B_{j})}\leq
 cR_{j}^{-2}{\left\Vert{\omega }\right\Vert}_{L^{r}(B_{j})},\label{HCF38}\end{equation}\
 \par 
with $\displaystyle t=S_{1}(r),\ s=S_{2}(r)$ still by lemma~\ref{3S3}
 \emph{(i)}. Let $\displaystyle q\in \lbrack r,t\rbrack .$\ \par 
\quad By Young's inequality we get, because $\displaystyle \frac{1}{t}=\frac{1}{s}+\frac{1}{n},$\
 \par 
\quad \quad \quad ${\left\Vert{u_{j}}\right\Vert}_{L^{t}(B_{j})}={\left\Vert{{\11}_{B_{j}}u_{j}}\right\Vert}_{L^{t}(B_{j})}\leq
 {\left\Vert{u_{j}}\right\Vert}_{L^{s}(B_{j})}{\left\Vert{{\11}_{B_{j}}}\right\Vert}_{L^{n}(B_{j})}={\left\Vert{u_{j}}\right\Vert}_{L^{s}(B_{j})}\left\vert{B_{j}}\right\vert
 ^{1/n}.$\ \par 
Because $\displaystyle \left\vert{B(x,R)}\right\vert :=\mathrm{V}\mathrm{o}\mathrm{l}(B(x,R))\leq
 (1+\epsilon )^{n/2}\nu _{n}R^{n}$ by equation~(\ref{HC23}),
 we get, with $\displaystyle c_{v}=c\root{n}\of{\nu _{n}(1+\epsilon
 )^{n/2}},\ \left\vert{B_{j}}\right\vert ^{1/n}\leq R_{j}.$\ \par 
Hence\ \par 
\quad \quad \quad $\displaystyle {\left\Vert{u_{j}}\right\Vert}_{L^{t}(B_{j})}\leq
 R_{j}\ {\left\Vert{u_{j}}\right\Vert}_{L^{s}(B_{j})}\leq c_{v}R_{j}^{-1}{\left\Vert{\omega
 }\right\Vert}_{L^{r}(B_{j})},$\ \par 
the last inequality given by~(\ref{HCF38}).\ \par 
Hence a fortiori $\displaystyle {\left\Vert{u_{j}}\right\Vert}_{L^{q}(B_{j})}\leq
 c_{v}R_{j}^{-1}{\left\Vert{\omega }\right\Vert}_{L^{r}(B_{j})}.$\ \par 
By lemma~\ref{HC34} we have $\displaystyle \left\vert{B(\chi
 _{j},u_{j})}\right\vert \leq \left\vert{\Delta \chi _{j}}\right\vert
 \left\vert{u_{j}}\right\vert +2\left\vert{\nabla \chi _{j}}\right\vert
 \left\vert{\nabla u_{j}}\right\vert ,$ so we get, because $\displaystyle
 \nabla u_{j}\in L^{q}(B_{j})$ by~(\ref{5S1}),\ \par 
\quad \quad \quad ${\left\Vert{B(\chi _{j},u_{j})}\right\Vert}_{L^{q}(B_{j})}\leq
 {\left\Vert{\nabla \chi _{j}}\right\Vert}_{\infty }{\left\Vert{\nabla
 u_{j}}\right\Vert}_{L^{q}(B_{j})}+{\left\Vert{\Delta \chi _{j}}\right\Vert}_{\infty
 }{\left\Vert{u_{j}}\right\Vert}_{L^{q}(B_{j})}\leq c_{v}R_{j}^{-2}{\left\Vert{\omega
 }\right\Vert}_{L^{r}(B_{j})}.$\ \par 
Multiplying by $\displaystyle w_{j}$ we get\ \par 
\quad \quad \quad $\displaystyle w_{j}{\left\Vert{B(\chi _{j},u_{j})}\right\Vert}_{L^{q}(B_{j})}\leq
 R_{j}^{-2}w_{j}c_{v}{\left\Vert{\omega }\right\Vert}_{L^{r}(B_{j})}.$\ \par 
Set $\displaystyle \omega _{1}:=\sum_{j\in {\mathbb{N}}}{B(\chi
 _{j},u_{j})},$ then\ \par 
\quad \quad \quad ${\left\Vert{\omega _{1}}\right\Vert}_{L^{q}(M,w^{q})}^{q}\leq
 \sum_{j\in {\mathbb{N}}}{{\left\Vert{B(\chi _{j},u_{,})}\right\Vert}_{L^{q}(B_{j},w^{q})}^{q}}$\
 \par 
Notice that $\displaystyle \chi _{j}B(\chi _{j},u_{j})=B(\chi
 _{j},u_{j}),$ so again we apply lemma~\ref{5S6} with $\displaystyle
 s=q,\ I\rightarrow {\left\Vert{\omega _{1}}\right\Vert}_{L^{q}(M,w^{q})}^{q},\
 u_{j}\rightarrow B(\chi _{j},u_{j})$ and $\displaystyle \gamma
 =2\ ,\ \tilde w(x):=R(x)^{-2}w(x),$ to get\ \par 
\quad \quad \quad $\displaystyle {\left\Vert{\omega _{1}}\right\Vert}_{L^{q}(M,w^{q})}\leq
 c_{w}T^{q/r}{\left\Vert{\omega }\right\Vert}_{L^{r}(M,\tilde w^{r})}.$\ \par 
Set $\displaystyle t_{1}=t=S_{1}(r),$ we have, with $\displaystyle
 w_{1}(x)=w(x),\ w_{0}(x)=\tilde w(x):=w(x)R(x)^{-2},\ \forall
 q\in \lbrack r,t_{1}\rbrack $\ \par 
\quad \quad \quad ${\left\Vert{\omega _{1}}\right\Vert}_{L^{q}(M,w_{1}^{q})}\leq
 c_{w}T^{q/r}{\left\Vert{\omega }\right\Vert}_{L^{r}(M,w_{0}^{r})}.$\ \par 
\ \par 
\quad If $\omega $ is of compact support and if $M$ is complete, by
 lemma~\ref{HCN48} we have seen that $\displaystyle v_{0}$ is
 also of compact support hence so is $\Delta v_{0}=\omega +\omega
 _{1}.$ Which means that $\omega _{1}$ is also of compact support.\ \par 
\ \par 
\quad Now we play the same game starting with $\omega _{1}$ in place
 of $\omega $ and we get, with\ \par 
\quad \quad \quad $\displaystyle s_{2}=S_{2}(t_{1}),\ t_{2}=S_{1}(t_{1})=S_{2}(r),\
 w_{2}(x)=w(x),\ w_{1}(x)=w(x)R(x)^{-2},\ w_{0}(x)=w(x)R(x)^{-4},$\ \par 
that\ \par 
\quad \quad \quad $\displaystyle \forall q\in \lbrack r,t_{1}\rbrack ,\ \forall
 s\in \lbrack r,s_{1}\rbrack ,\ \exists v_{1}\in L^{s}(M,w_{2}^{s})\cap
 W^{2,q}(M,w_{1})::\Delta v_{1}=\omega _{1}+\omega _{2}$\ \par 
and\ \par 
\quad \quad \quad \quad $\displaystyle \forall t\in \lbrack r,t_{2}\rbrack ,\ \omega
 _{2}\in L^{t}(M,w_{2}^{t}),\ \ {\left\Vert{\omega _{2}}\right\Vert}_{L^{t}(M,w_{2}^{t})}\lesssim
 {\left\Vert{\omega _{1}}\right\Vert}_{L^{t_{1}}(M,w_{1}^{t_{1}})}\lesssim
 {\left\Vert{\omega }\right\Vert}_{L^{r}(M,w_{0}^{r})}.$\ \par 
We keep the linearity of $\displaystyle v_{1}$ w.r.t. to $\omega
 _{1}$ hence to $\omega .$ So $\omega _{2}$ is still linear w.r.t.
 $\omega .$\ \par 
\quad So by induction we have, with\ \par 
\quad \quad \quad $\displaystyle t_{k}=S_{k}(r),\ w_{k}(x):=w(x),\ w_{k-1}(x)=w(x)R(x)^{-2},...,\
 w_{0}(x)=w(x)R(x)^{-2k},$\ \par 
and, with $\displaystyle s_{j+1}=S_{2}(t_{j}),$\ \par 
\quad \quad \quad $\displaystyle \forall s\in \lbrack r,s_{j+1}\rbrack ,\ \forall
 q\in \lbrack r,\ s_{j}\rbrack \ \forall j=0,...,\ k-1,\ v_{j}\in
 L^{q}(M,w_{j+1}^{q})\cap W^{2,q}(M,w_{j}),$\ \par 
and\ \par 
\quad \quad \quad \begin{equation}  \forall q\in \lbrack r,t_{k}\rbrack ,\ \omega
 _{k}\in L^{q}(M,w_{k}^{q}),\ {\left\Vert{\omega _{k}}\right\Vert}_{L^{q}(M,w_{k}^{q})}\lesssim
 \cdot \cdot \cdot \lesssim {\left\Vert{\omega _{1}}\right\Vert}_{L^{t_{1}}(M,w_{1}^{t_{1}})}\lesssim
 {\left\Vert{\omega }\right\Vert}_{L^{r}(M,w_{0}^{r})}.\label{CL17}\end{equation}\
 \par 
Setting now $\displaystyle v:=\sum_{j=0}^{k-1}{(-1)^{j}v_{j}}$
 and $\tilde \omega :=(-1)^{k}\omega _{k},$ we have that $\Delta
 v=\omega +\tilde \omega $ and\ \par 
\quad \quad \quad $\displaystyle \forall q\in \lbrack r,s_{1}\rbrack ,\ v_{j}\in
 L^{q}(M,w_{j+1}^{q})\cap W^{2,r}(M,w),\ s_{1}=S_{2}(t_{1}),\
 w_{j+1}=w(x)R(x)^{2(j+1-k)},$\ \par 
this implies, because $\displaystyle w_{k}=w\leq w_{j+1},$\ \par 
\quad \quad \quad $\displaystyle \forall q\in \lbrack r,s_{1}\rbrack ,\ v_{j}\in
 L^{q}(M,w^{q}),\ {\left\Vert{v_{j}}\right\Vert}_{L^{q}(M,w^{q})}\leq
 c_{l}T^{1/r}{\left\Vert{\omega }\right\Vert}_{L^{r}(M,w_{0}^{r})}.$\ \par 
So we have also for $\displaystyle v:=\sum_{j=0}^{k-1}{(-1)^{j}v_{j}}\ :$\ \par 
\quad \quad \quad \begin{equation} \forall q\in \lbrack r,s_{1}\rbrack ,\ v\in
 L^{q}(M,w^{q}),\ {\left\Vert{v}\right\Vert}_{L^{q}(M,w^{q})}\leq
 kc_{l}T^{1/r}{\left\Vert{\omega }\right\Vert}_{L^{r}(M,w_{0}^{r})}\label{CF8}
 .\end{equation}\ \par 
We cannot go beyond $\displaystyle s_{1}:=S_{2}(r)$ for $v$ because
 of $\displaystyle v_{0}.$ For the same reason, we cannot go
 beyond $\displaystyle W^{2,r}(M,w).$\ \par 
For the remaining term $\displaystyle \tilde \omega ,$ we get
 a better regularity, still because we set $\displaystyle w_{k}=w,\
 \tilde \omega =(-1)^{k}\omega _{k},$\ \par 
\quad \quad \quad \begin{equation}  \forall q\in \lbrack r,t_{k}\rbrack ,\ \tilde
 \omega \in L^{q}(M,w^{q}),\ {\left\Vert{\tilde \omega }\right\Vert}_{L^{q}(M,w^{q})}\lesssim
 {\left\Vert{\omega }\right\Vert}_{L^{r}(M,w_{0}^{r})}.\label{CF7}\end{equation}\
 \par 
Clearly the linearity is kept along the induction.\ \par 
Now we choose $k$ such that the threshold $\displaystyle t_{k}:=S_{k}(r)\geq
 s.$\ \par 
\quad If $\omega $ is of compact support and if $M$ is complete, by
 lemma~\ref{HCN48} we have seen that $\displaystyle v_{0}$ and
 $\omega _{1}$ also and by induction all the $\displaystyle v_{j}$
 and $\omega _{j}$ are also of compact support. $\blacksquare $\ \par 
\quad We shall refer to this theorem as RSM for short. We notice that
 we have no completeness assumption on $M$  to get the first
 part of the result.\ \par 

\begin{lem}
~\label{HCN46} Set, for $k\in {\mathbb{N}},\ w_{k}^{q}=R(x)^{-qk},$
 we have $\displaystyle L^{q}_{p}(M,w_{k}^{q})\subset L^{q}_{p}(M)$ and\par 
\quad \quad \quad $\displaystyle \forall q>1,\ \forall f\in L^{q}_{p}(M,w_{k}^{q}),\
 {\left\Vert{f}\right\Vert}_{L^{q}_{p}(M)}\leq {\left\Vert{f}\right\Vert}_{L^{q}_{p}(M,w_{k}^{q})}.$
\end{lem}

      Proof.\ \par 
We have, because $\displaystyle \forall x\in M,\ R(x)\leq 1\Rightarrow
 w_{k}^{q}(x)\geq 1,$\ \par 
\quad \quad \quad ${\left\Vert{f}\right\Vert}_{L^{q}_{p}(M)}^{q}=\int_{M}{\left\vert{f}\right\vert
 ^{q}dv_{g}}\leq \int_{M}{\left\vert{f}\right\vert ^{q}w_{k}^{q}dv_{g}}={\left\Vert{f}\right\Vert}_{L^{q}_{p}(M,w_{k}^{q})}^{q}.\
 \blacksquare $\ \par 

\begin{rem}
~\label{CL18}We have, by inequalities~(\ref{CL17}), that $\displaystyle
 \forall q\in \lbrack r,t_{k}\rbrack ,\ \tilde \omega \in L^{q}(M,w^{q}).$
 With the choice of $\displaystyle w\equiv 1$ for the weight
 relative to the covering, with the notations of the RSM, we
 get $\displaystyle \forall q\in \lbrack r,t_{k}\rbrack ,\ \tilde
 \omega \in L^{q}(M).$\par 
\quad We also have that $\displaystyle w\equiv 1\Rightarrow \forall
 q\in \lbrack r,s_{1}\rbrack ,\ v\in L^{q}(M),$ with $\displaystyle
 s_{1}:=S_{2}(r).$
\end{rem}
\ \par 

\begin{cor}
~\label{CF9}Let $\displaystyle (M,g)$ be a complete riemannian
 manifold. For $\displaystyle r\leq 2,$ take $\displaystyle \omega
 \in L^{r}(M,w_{0}^{r})$ with $\displaystyle k\in {\mathbb{N}},\
 w_{0}(x):=R(x)^{-2k}$ and $\displaystyle s_{1}:=S_{2}(r).$\par 
Chosing $k$ big enough for the threshold $\displaystyle t_{k}:=S_{k}(r)\geq
 2$ then the orthogonal projection $H\ :\ L^{2}_{p}(M)\rightarrow
 {\mathcal{H}}^{2}_{p}(M)$ extends boundedly from $\displaystyle
 L^{r}_{p}(M,w_{0}^{r})$ to $\displaystyle {\mathcal{H}}^{2}_{p}(M).$
 This implies $\displaystyle H\omega =0\iff H\tilde \omega =0.$
\end{cor}
\quad Proof.\ \par 
For $\displaystyle \omega \in L^{r}_{p}(M,w_{0}^{r})$ and $\displaystyle
 \forall s\in \lbrack r,s_{1}\rbrack ,\ \forall q\in \lbrack
 r,t_{k}\rbrack ,$ the RSM, theorem~\ref{CF3},  gives us two
 forms $\displaystyle v\in L_{p}^{s}(M,w^{s}),\ \tilde \omega
 \in L^{q}_{p}(M,w^{q}),$ such that\ \par 
\quad \quad \quad \begin{equation}  v=T\omega ,\ \tilde \omega =A\omega ,\ \Delta
 v=\omega +\tilde \omega .\label{HCN47}\end{equation}\ \par 
where $T$ and $A$ are bounded linear operators :\ \par 
\quad \quad \quad $\displaystyle T\ :\ L^{r}_{p}(M,w_{0}^{r})\rightarrow L_{p}^{s}(M,w^{s})\
 ;\ A\ :\ L^{r}_{p}(M,w_{0}^{r})\rightarrow L^{q}_{p}(M,w^{q}).$\ \par 
Now choose $\displaystyle w\equiv 1\Rightarrow w_{0}=R(x)^{-2k}.$
 Then we have $\displaystyle \tilde \omega \in L^{q}_{p}(M),$
 and if $k$ is such that the threshold $\displaystyle t_{k}=S_{k}(r)\geq
 2,$ we have $\displaystyle \tilde \omega \in L^{2}_{p}(M).$\ \par 
\quad Hence the projection $H$ is well defined on $\displaystyle \tilde
 \omega .$ Suppose that $\displaystyle H\Delta v=0$ then we were
 done because, by~(\ref{HCN47}), we would have $\displaystyle
 0=H\Delta v=H\omega +H\tilde \omega \Rightarrow H\omega =-H\tilde
 \omega .$\ \par 
\quad We start by approximating $\displaystyle \omega $ by a sequence
 $\omega _{l}\in {\mathcal{D}}_{p}(M),\ \omega _{l}\rightarrow
 \omega $ in $\displaystyle L^{r}_{p}(M,w_{0}^{r}).$ Then apply
 the RSM to $\displaystyle \omega _{l}\ ;$ we get $\displaystyle
 v_{l}=T\omega _{l},\ \tilde \omega _{l}=A\omega _{l},\ \Delta
 v_{l}=\omega _{l}+\tilde \omega _{l}.$ We have that $\displaystyle
 v_{l},\ \tilde \omega _{l}$ have compact support and by linearity
 with~(\ref{CF8})\ \par 
\quad \quad \quad $\displaystyle \forall s\in \lbrack r,s_{1}\rbrack ,\ (v-v_{l})\in
 L^{s}(M),\ {\left\Vert{v-v_{l}}\right\Vert}_{L^{s}(M)}\leq kc_{l}T^{1/r}{\left\Vert{\omega
 -\omega _{l}}\right\Vert}_{L^{r}(M,w_{0}^{r})}$\ \par 
so $\displaystyle {\left\Vert{v-v_{l}}\right\Vert}_{L^{s}(M)}\rightarrow
 0$ and the same way with~(\ref{CF7}) we get\ \par 
\quad \quad \quad $\displaystyle \forall q\in \lbrack r,t_{k}\rbrack ,\ (\tilde
 \omega -\tilde \omega _{l})\in L^{q}(M),\ {\left\Vert{\tilde
 \omega -\tilde \omega _{l}}\right\Vert}_{L^{q}(M)}\lesssim {\left\Vert{\omega
 -\omega _{l}}\right\Vert}_{L^{r}(M,w_{0}^{r})}$\ \par 
hence $\displaystyle {\left\Vert{\tilde \omega -\tilde \omega
 _{l}}\right\Vert}_{L^{q}(M)}\rightarrow 0.$\ \par 
Then $H$ is well defined on $v_{l},\ \Delta v_{l},\ \omega _{l}$
 and $\tilde \omega _{l},$ because they are ${\mathcal{C}}^{\infty
 }$ and compactly supported hence in $\displaystyle L^{2}_{p}(M),$
 and we have\ \par 
\quad \quad \quad $\Delta v_{l}=\omega _{l}+\tilde \omega _{l}\Rightarrow H\Delta
 v_{l}=H\omega _{l}+H\tilde \omega _{l}.$\ \par 
Take $\displaystyle t_{k}\geq 2,\ h\in L^{2}_{p}(M)$ then $\displaystyle
 \ {\left\langle{H\Delta v_{l},h}\right\rangle}_{L^{2}(M)}={\left\langle{\Delta
 v_{l},Hh}\right\rangle}_{L^{2}(M)}$ because $H$ is self adjoint.
 But because $M$ is complete, $\Delta $ is essentially self adjoint
 on $\displaystyle L^{2}_{p}(M)$ by~\cite{Gaffney55} and $\displaystyle
 v_{l}$ has compact support, we have\ \par 
\quad \quad \quad $\displaystyle \ {\left\langle{\Delta v_{l},Hh}\right\rangle}_{L^{2}(M)}={\left\langle{v_{l},\Delta
 Hh}\right\rangle}_{L^{2}(M)}=0,$\ \par 
because $Hh\in {\mathcal{H}}_{p}^{2}(M).$\ \par 
So we have $\displaystyle \forall l\in {\mathbb{N}},\ H\Delta
 v_{l}=0$ and this implies\ \par 
\quad \quad \quad $\displaystyle \forall l\in {\mathbb{N}},\ H\omega _{l}+H\tilde
 \omega _{l}=0.$\ \par 
Now we have $\tilde \omega \in L^{2}_{p}(M)$ and the convergence
 $\displaystyle {\left\Vert{\tilde \omega -\tilde \omega _{l}}\right\Vert}_{L^{2}_{p}(M)}\rightarrow
 0$ by lemma~\ref{HCN46}. So, because $H$ is bounded on $\displaystyle
 L^{2}_{p}(M),$ we get $\displaystyle H\tilde \omega _{l}\rightarrow
 H\tilde \omega $ in $\displaystyle L^{2}_{p}(M),$ and this means
 $\displaystyle H\omega _{l}\rightarrow -H\tilde \omega $ also
 in $\displaystyle L^{2}_{p}(M).$ So we define, for any sequence
 $\omega _{l}\in {\mathcal{D}}_{p}(M),\ \omega _{l}\rightarrow
 \omega $ in $\displaystyle L^{r}_{p}(M,w_{0}),$ by :\ \par 
\quad \quad \quad $\displaystyle H\omega :=\lim \ _{l\rightarrow \infty }H\omega
 _{l}=-\lim \ _{l\rightarrow \infty }H\tilde \omega _{l}=-H\tilde
 \omega ,$\ \par 
so we proved that $\displaystyle H\omega _{l}$ converges in $\displaystyle
 L^{2}_{p}(M)$ to $\displaystyle -H\tilde \omega ,$ with $\tilde
 \omega $ given by the Raising Steps Method. This limit is independent
 of the sequence of approximations $\displaystyle \omega _{l},$
 and it is clearly a extension of the projection $H$ to $\displaystyle
 L^{r}(M,w_{0}^{r}).$\ \par 
This implies that $\displaystyle H\omega =0\iff H\tilde \omega
 =0.$ $\blacksquare $\ \par 

\begin{cor}
~\label{HCS44}Let $\displaystyle (M,g)$ be a complete riemannian
 manifold. We get : $\forall s\geq 2,\ {\mathcal{H}}_{p}^{2}(M)\hookrightarrow
 {\mathcal{H}}_{p}^{s}(M)$ with\par 
\quad \quad \quad $\displaystyle \forall h\in {\mathcal{H}}_{p}^{2}(M),\ h\in {\mathcal{H}}_{p}^{s}(M)$
 and $\displaystyle {\left\Vert{h}\right\Vert}_{L^{s}_{p}(M)}\leq
 C_{s}{\left\Vert{h}\right\Vert}_{L^{2}_{p}(M)}.$
\end{cor}
\quad Proof.\ \par 
Let $\omega \in {\mathcal{D}}_{p}(M)$ and $\displaystyle \varphi
 \in L^{2}_{p}(M),$ then we have $\displaystyle \ {\left\langle{H\omega
 ,\ \varphi }\right\rangle}={\left\langle{\omega ,\ H^{*}\varphi
 }\right\rangle}$ by duality ; on the other hand, because $\omega
 \in L^{2}_{p}(M),$ we get\ \par 
\quad \quad \quad $\ {\left\langle{H\omega ,\ \varphi }\right\rangle}={\left\langle{\omega
 ,\ H\varphi }\right\rangle}\ ;$\ \par 
so, against $\displaystyle {\mathcal{D}}_{p}(M),$ we have $\displaystyle
 H=H^{*}.$\ \par 
Now take $\displaystyle r\leq 2,$ and $\displaystyle \omega \in
 L^{r}(M,w_{0}^{r})$ with $\displaystyle k\in {\mathbb{N}},\
 w_{0}(x):=R(x)^{-k}$ and $\displaystyle s_{1}:=S_{2}(r).$ Chosing
 $k$ big enough for the threshold $\displaystyle t_{k}:=S_{k}(r)\geq
 2,$ then the orthogonal projection $H\ :\ L^{2}_{p}(M)\rightarrow
 L^{2}_{p}(M)$ extends boundedly from $\displaystyle L^{r}_{p}(M,w_{0}^{r})$
 to $\displaystyle L^{2}_{p}(M)$ by corollary~\ref{CF9} hence
 by duality $\displaystyle H^{*}\ :\ L^{2}_{p}(M)\rightarrow
 L^{r'}_{p}(M,w_{0}^{r}).$\ \par 
By density of ${\mathcal{D}}_{p}(M)$ in $\displaystyle L^{r}(M,w_{0}^{r})$
 we get that $\displaystyle H=H^{*}\ :\ L^{2}_{p}(M)\rightarrow
 L^{r'}_{p}(M,w_{0}^{r}).$\ \par 
\quad We also have by lemma~\ref{HCN46} $\displaystyle L^{r'}_{p}(M,w_{0}^{r})\subset
 L^{r'}_{p}(M)$ with norm less than one, hence $\displaystyle
 H\ :\ L^{2}_{p}(M)\rightarrow L^{r'}_{p}(M)$ boundedly and\ \par 
\quad \quad \quad $h\in {\mathcal{H}}_{p}^{2}(M)\Rightarrow h=Hh\in L^{r'}_{p}(M)\Rightarrow
 h\in {\mathcal{H}}_{p}^{r'}(M).$\ \par 
Now we choose $\displaystyle r=s'$ the conjugate exponent of
 $s$ to end the proof of the corollary. $\blacksquare $\ \par 
\quad We already know that harmonic forms are smooth, see for instance~[\cite{Bueler99}
 corollary 5.4], so corollary~\ref{HCS44} gives another kind
 of smoothness.\ \par 

\section{Weighted Calderon Zygmund inequalities.}
\quad In the same spirit of theorem 1.2 by Guneysu and Pigola~\cite{GuneysuPigola},
 we get the following "twisted" Calderon Zygmund inequality with
 weights and being valid directly for forms not a priori in ${\mathcal{D}}_{p}(M).$\
 \par 
\quad These CZI are twisted because there are 2 different weights in
 the inequality.\ \par 

\begin{thm}
~\label{CF11} Let $\displaystyle (M,g)$ be a complete riemannian
 manifold. Let $w$ be a weight relative to the ${\mathcal{C}}_{\epsilon
 }$ associated covering $\displaystyle \lbrace B(x_{j},5r(x_{j}))\rbrace
 _{j\in {\mathbb{N}}}$ and set $\displaystyle w_{0}:=R(x)^{-2}.$
 Let $\displaystyle u\in L^{r}_{p}(M,ww_{0}^{r})$ such that $\Delta
 u\in L^{r}_{p}(M,w)$ ; then there are constants $\displaystyle
 C_{1},C_{2}$ depending only on $n=\mathrm{d}\mathrm{i}\mathrm{m}_{{\mathbb{R}}}M,\
 r$ and $\epsilon $ such that:\par 
\quad \quad \quad $\displaystyle {\left\Vert{u}\right\Vert}_{W^{2,r}(M,w)}\leq
 C_{1}{\left\Vert{u}\right\Vert}_{L^{r}(M,ww_{0}^{r})}+C_{2}{\left\Vert{\Delta
 u}\right\Vert}_{L^{r}(M,w)}.$\par 
Moreover we have for $\displaystyle t=S_{2}(r)$ that $\displaystyle
 u\in L^{t}_{p}(M,w^{t})$ with $\displaystyle {\left\Vert{u}\right\Vert}_{L^{t}(M,w^{t})}\leq
 c{\left\Vert{u}\right\Vert}_{W^{2,r}(M,w^{r}w_{0}^{t})}.$
\end{thm}
\quad Proof.\ \par 
Let $\displaystyle u\in L^{r}(M,ww_{0}^{r}),\ \Delta u\in L^{r}(M,w).$
 Set $\displaystyle R_{j}:=5r(x_{j}),\ B_{j}:=B(x_{j},R_{j}),\
 B'_{j}=B(x_{j},2R_{j})$ and apply lemma~\ref{HF0} to get :\ \par 
there are constants $\displaystyle c_{1},c_{2}$ depending only
 on $n=\mathrm{d}\mathrm{i}\mathrm{m}_{{\mathbb{R}}}M,\ r,\ \epsilon
 $ such that\ \par 
\quad \quad \quad \begin{equation} {\left\Vert{ u}\right\Vert}_{W^{2,r}(B_{j})}\leq
 c_{1}R_{j}^{-2}{\left\Vert{u}\right\Vert}_{L^{r}(B'_{j})}+c_{2}{\left\Vert{\Delta
 u}\right\Vert}_{L^{r}(B'_{j})}.\label{7CZ0}\end{equation}\ \par 
Recall that\ \par 
\quad \quad \quad $\displaystyle {\left\Vert{u}\right\Vert}_{W^{2,r}(M,w)}:={\left\Vert{\nabla
 ^{2}u}\right\Vert}_{L^{r}(M,w)}+{\left\Vert{\nabla u}\right\Vert}_{L^{r}(M,w)}+{\left\Vert{u}\right\Vert}_{L^{r}(M,w)},$\
 \par 
so we have to compute those three terms.\ \par 
\quad \quad \quad ${\left\Vert{\nabla ^{2}u}\right\Vert}_{L^{r}(M,w)}^{r}=\int_{M}{\left\vert{\nabla
 ^{2}u}\right\vert ^{r}wdv_{g}}\leq \sum_{j\in {\mathbb{N}}}{{\left\Vert{\nabla
 ^{2}u}\right\Vert}^{r}_{L^{r}(B_{j}w)}}.$\ \par 
By~(\ref{7CZ0}) we get\ \par 
\quad \quad \quad $\displaystyle {\left\Vert{\nabla ^{2}u}\right\Vert}_{L^{r}(B_{j},w)}\leq
 (c_{1}R_{j}^{-2}c_{sw}w_{j}{\left\Vert{u}\right\Vert}_{L^{r}(B'_{j})}+c_{2}c_{sw}w_{j}{\left\Vert{\Delta
 u}\right\Vert}_{L^{r}(B'_{j})})^{r}\leq $\ \par 
\quad \quad \quad \quad \quad \quad \quad $\displaystyle \leq 2^{r/r'}c_{sw}^{r}(c_{1}R_{j}^{-2r}w_{j}{\left\Vert{u}\right\Vert}^{r}_{L^{r}(B'_{j})}+c_{2}w_{j}{\left\Vert{\Delta
 u}\right\Vert}^{r}_{L^{r}(B'_{j})}).$\ \par 
Hence\ \par 
\quad \quad \quad \begin{equation} {\left\Vert{ \nabla ^{2}u}\right\Vert}_{L^{r}(M,w)}^{r}\leq
 2^{r/r'}c_{sw}^{r}c_{iw}^{-r}c_{1}\sum_{j\in {\mathbb{N}}}{R_{j}^{-2r}{\left\Vert{u}\right\Vert}^{r}_{L^{r}(B'_{j},w)}+c_{2}w_{j}{\left\Vert{\Delta
 u}\right\Vert}^{r}_{L^{r}(B'_{j},w)}).}\label{7CZ1}\end{equation}\ \par 
Exactly as in the proof of the RSM we get\ \par 
\quad \quad \quad $\displaystyle R_{j}^{-2r}\int_{B'_{j}}{\left\vert{u(x)}\right\vert
 ^{r}w(x)dv_{g}(x)}\leq 96^{2r}\int_{B'_{j}}{\left\vert{u(x)}\right\vert
 ^{r}R(x)^{-2r}w(x)dv_{g}(x)}$\ \par 
hence, because the overlap of the Vitali covering is bounded
 by $T,$ even for the double balls $\displaystyle B'_{j},$ we get\ \par 
\quad \quad \quad $\sum_{j\in {\mathbb{N}}}{R_{j}^{-2r}\int_{B_{'j}}{\left\vert{u(x)}\right\vert
 ^{r}w(x)dv_{g}(x)}}\leq 96^{2r}T\int_{M}{\left\vert{u(x)}\right\vert
 ^{r}R(x)^{-2r}w(x)dv_{g}(x)}.$\ \par 
Easier we get\ \par 
\quad \quad \quad $\displaystyle \sum_{j\in {\mathbb{N}}}{\int_{B'_{j}}{\left\vert{\Delta
 u(x)}\right\vert ^{r}w(x)dv_{g}(x)}}\leq T\int_{M}{\left\vert{\Delta
 u(x)}\right\vert ^{r}w(x)dv_{g}(x)}.$\ \par 
So, putting in~(\ref{7CZ1}), we get\ \par 
\quad \quad \quad $\displaystyle {\left\Vert{\nabla ^{2}u}\right\Vert}_{L^{r}(M,w)}\leq
 2^{1/r'}96^{2}T^{1/r}c_{1}{\left\Vert{u}\right\Vert}_{L^{r}(M,ww_{0}^{r})}+2^{1/r'}T^{1/r}c_{2}{\left\Vert{\Delta
 u}\right\Vert}_{L^{r}(M,w)}.$\ \par 
\quad Exactly the same way we get\ \par 
\quad \quad \quad $\displaystyle {\left\Vert{\nabla u}\right\Vert}_{L^{r}(M,w)}\leq
 2^{1/r'}96^{2}T^{1/r}c_{1}{\left\Vert{u}\right\Vert}_{L^{r}(M,ww_{0}^{r})}+2^{1/r'}T^{1/r}c_{2}{\left\Vert{\Delta
 u}\right\Vert}_{L^{r}(M,w)}.$\ \par 
Hence\ \par 
\quad \quad \quad $\displaystyle {\left\Vert{u}\right\Vert}_{W^{2,r}(M,w)}\leq
 C_{1}{\left\Vert{u}\right\Vert}_{L^{r}(M,ww_{0}^{r})}+C_{2}{\left\Vert{\Delta
 u}\right\Vert}_{L^{r}(M,w)}$\ \par 
with\ \par 
\quad \quad \quad $\displaystyle C_{1}:=1+2^{1/r'}96^{2}T^{1/r}c_{1}\ ;\ C_{2}:=2^{1/r'}T^{1/r}c_{2}.$\
 \par 
\ \par 
\quad To get the "moreover" we proceed the same way. By lemma~\ref{3S2}
 \emph{(i)}, we get for the $\epsilon $ admissible ball $\displaystyle
 B_{j},$\ \par 
\quad \quad \quad $\displaystyle t=S_{2}(r),\ \forall u\in W^{2,r}(B_{j}),\ {\left\Vert{u}\right\Vert}_{L^{t}(B_{j})}\leq
 CR_{j}^{-2}\ {\left\Vert{u}\right\Vert}_{W^{2,r}(B_{j})}.$\ \par 
So, because $w$ is relative to the covering,\ \par 
\quad \quad \quad ${\left\Vert{u}\right\Vert}_{L^{t}(M,w)}^{t}\leq c_{sw}\sum_{j\in
 {\mathbb{N}}}{w_{j}\int_{B_{j}}{\left\vert{u}\right\vert ^{t}dv_{g}}}=c_{sw}\sum_{j\in
 {\mathbb{N}}}{w_{j}{\left\Vert{u}\right\Vert}_{L^{t}(B_{j})}^{t}}.$\ \par 
But, as above,\ \par 
\quad \quad \quad $\displaystyle \int_{B_{j}}{\left\vert{u}\right\vert ^{t}dv_{g}}\leq
 C^{t}R_{j}^{-2t}\ {\left\Vert{u}\right\Vert}^{t}_{W^{2,r}(B_{j})}\leq
 C^{t}R_{j}^{-2t}({\left\Vert{\nabla ^{2}u}\right\Vert}^{t}_{L^{r}(B_{j})}+{\left\Vert{\nabla
 u}\right\Vert}^{t}_{L^{r}(B_{j})}+{\left\Vert{u}\right\Vert}^{t}_{L^{r}(B_{j})}).$\
 \par 
Hence\ \par 
\quad \quad \quad $\displaystyle {\left\Vert{u}\right\Vert}_{L^{t}(M,w^{t})}^{t}\leq
 C^{t}c_{sw}\sum_{j\in {\mathbb{N}}}{w_{j}^{t}R_{j}^{-2t}({\left\Vert{\nabla
 ^{2}u}\right\Vert}^{t}_{L^{r}(B_{j})}+{\left\Vert{\nabla u}\right\Vert}^{t}_{L^{r}(B_{j})}+{\left\Vert{u}\right\Vert}^{t}_{L^{r}(B_{j})}).}$\
 \par 
But\ \par 
\quad \quad \quad $\displaystyle \sum_{j\in {\mathbb{N}}}{a_{j}^{t}}\leq (\sum_{j\in
 {\mathbb{N}}}{a_{j}^{r}})^{t/r}$ because $\displaystyle t\geq r,$\ \par 
hence we get\ \par 
\quad \quad \quad $\displaystyle A_{1}:=\sum_{j\in {\mathbb{N}}}{w_{j}^{t}R_{j}^{-2t}{\left\Vert{\nabla
 ^{2}u}\right\Vert}^{t}_{L^{r}(B_{j})}}\leq (\sum_{j\in {\mathbb{N}}}{w_{j}^{r}R_{j}^{-2r}{\left\Vert{\nabla
 ^{2}u}\right\Vert}^{r}_{L^{r}(B_{j})}})^{t/r}.$\ \par 
hence, putting the radius and the weight into the integral, which
 gives the $\displaystyle w_{0}^{t},$\ \par 
\quad \quad \quad $\displaystyle w_{j}^{r}R_{j}^{-2r}{\left\Vert{\nabla ^{2}u}\right\Vert}^{r}_{L^{r}(B_{j})}\leq
 c_{iw}^{-r}\int_{B_{j}}{\left\vert{\nabla ^{2}u}\right\vert
 ^{r}w^{r}w_{0}^{r}dv_{g}},$\ \par 
So\ \par 
\quad \quad \quad $\displaystyle A_{1}\leq c_{iw}^{-t}(\sum_{j\in {\mathbb{N}}}{\int_{B_{j}}{\left\vert{\nabla
 ^{2}u}\right\vert ^{r}w^{r}w_{0}^{r}dv_{g}}})^{t/r}.$\ \par 
The overlap of the Vitali covering is bounded by $T,$ so\ \par 
\quad \quad \quad $\displaystyle A_{1}\leq c_{iw}^{-t}(T\int_{M}{\left\vert{\nabla
 ^{2}u}\right\vert ^{r}w^{r}w_{0}^{r}dv_{g}})^{t/r}=c_{iw}^{-t}T^{t/r}{\left\Vert{\nabla
 ^{2}u}\right\Vert}_{L^{r}(M,w^{r}w_{0}^{r})}^{t}.$\ \par 
Exactly the same way, we get\ \par 
\quad \quad \quad $\displaystyle A_{2}:=\sum_{j\in {\mathbb{N}}}{w_{j}^{t}R_{j}^{-2t}{\left\Vert{\nabla
 u}\right\Vert}^{t}_{L^{r}(B_{j})}}\leq c_{iw}^{-t}T^{t/r}{\left\Vert{\nabla
 u}\right\Vert}_{L^{r}(M,w^{r}w_{0}^{r})}^{t},$\ \par 
and\ \par 
\quad \quad \quad $\displaystyle A_{3}:=\sum_{j\in {\mathbb{N}}}{w_{j}^{t}R_{j}^{-2t}{\left\Vert{u}\right\Vert}^{t}_{L^{r}(B_{j})}}\leq
 c_{iw}^{-t}T^{t/r}{\left\Vert{u}\right\Vert}_{L^{r}(M,w^{r}w_{0}^{r})}^{t},$\
 \par 
Adding we get\ \par 
\quad \quad \quad $\displaystyle {\left\Vert{u}\right\Vert}_{L^{t}(M,w^{t})}^{t}\leq
 C^{t}c_{sw}c_{iw}^{-t}(A_{1}+A_{2}+A_{3})\leq $\ \par 
\quad \quad \quad \quad \quad \quad \quad $\displaystyle \leq C^{t}c_{sw}c_{iw}^{-t}T^{t/r}({\left\Vert{\nabla
 ^{2}u}\right\Vert}_{L^{r}(M,w^{r}w_{0}^{r})}^{t}+{\left\Vert{\nabla
 u}\right\Vert}_{L^{r}(M,w^{r}w_{0}^{r})}^{t}+{\left\Vert{u}\right\Vert}_{L^{r}(M,w^{r}w_{0}^{r})}^{t})\leq
 $\ \par 
\quad \quad \quad \quad \quad \quad \quad $\displaystyle \leq C^{t}c_{sw}c_{iw}^{-t}T^{t/r}({\left\Vert{\nabla
 ^{2}u}\right\Vert}_{L^{r}(M,w^{r}w_{0}^{r})}+{\left\Vert{\nabla
 u}\right\Vert}_{L^{r}(M,w^{r}w_{0}^{r})}+{\left\Vert{u}\right\Vert}_{L^{r}(M,w^{r}w_{0}^{r})})^{t}\leq
 $\ \par 
\quad \quad \quad \quad \quad \quad \quad $\leq C^{t}c_{sw}c_{iw}^{-t}T^{t/r}({\left\Vert{u}\right\Vert}_{W^{r}(M,w^{r}w_{0}^{r})})^{t}.$\
 \par 
Taking the $t$ root we get\ \par 
\quad \quad \quad $\displaystyle {\left\Vert{u}\right\Vert}_{L^{t}(M,w^{t})}\leq
 Cc_{sw}c_{iw}^{-1}T^{1/r}{\left\Vert{u}\right\Vert}_{W^{r}(M,w^{r}w_{0}^{r})}.$\
 \par 
\quad Which ends the proof of the theorem. $\blacksquare $\ \par 

\begin{cor}
Let $\displaystyle (M,g)$ be a complete riemannian manifold.
 Set $\displaystyle w_{0}:=R(x)^{-2}.$ Let $\displaystyle u\in
 L^{r}_{p}(M,w_{0}^{r})$ such that $\Delta u\in L^{r}_{p}(M)$
 ; then there are constants $\displaystyle C_{1},C_{2}$ depending
 only on $n=\mathrm{d}\mathrm{i}\mathrm{m}_{{\mathbb{R}}}M,\
 r$ and $\epsilon $ such that :\par 
\quad \quad \quad $\displaystyle {\left\Vert{u}\right\Vert}_{W_{p}^{2,r}(M)}\leq
 C_{1}{\left\Vert{u}\right\Vert}_{L_{p}^{r}(M,w_{0}^{r})}+C_{2}{\left\Vert{\Delta
 u}\right\Vert}_{L_{p}^{r}(M)}.$\par 
Moreover we have for $\displaystyle t=S_{2}(r)$ that $\displaystyle
 u\in L^{t}_{p}(M)$ with $\displaystyle \ {\left\Vert{u}\right\Vert}_{L_{p}^{t}(M)}\leq
 c{\left\Vert{u}\right\Vert}_{W_{p}^{2,r}(M,w_{0}^{t})}.$
\end{cor}
\quad Proof.\ \par 
We choose the weight $\displaystyle w\equiv 1.$ $\blacksquare $\ \par 

\begin{cor}
~\label{7CZ2} If the complete riemannian manifold $\displaystyle
 (M,g)$ is such that the $\displaystyle \epsilon _{0}$ admissible
 radius is positive, then we get the classical Calderon Zygmund
 inequalities :\par 
\quad \quad \quad $\displaystyle \forall r,\ 1<r<\infty ,\ {\left\Vert{u}\right\Vert}_{W^{2,r}(M)}\leq
 C_{1}{\left\Vert{u}\right\Vert}_{L^{r}(M)}+C_{2}{\left\Vert{\Delta
 u}\right\Vert}_{L^{r}(M)}.$\par 
Moreover we have the classical Sobolev inequality :\par 
\quad \quad \quad for $\displaystyle t=S_{2}(r)$ we get that $\displaystyle u\in
 L^{t}_{p}(M)$ with $\displaystyle {\left\Vert{u}\right\Vert}_{L^{t}(M)}\leq
 c{\left\Vert{u}\right\Vert}_{W^{2,r}(M)}.$
\end{cor}
\quad Proof.\ \par 
If $\displaystyle \forall x\in M,\ R(x)\geq \delta >0,$ then
 $\displaystyle w_{0}(x)^{r}\simeq 1$ hence the weights disappear.
 $\blacksquare $\ \par 
\quad Recall that, by theorem 1.3 in Hebey~\cite{Hebey96}, we have
 that the harmonic radius $\displaystyle r_{H}(1+\epsilon ,\
 2,0)$ is bounded below if the Ricci curvature $\displaystyle
 Rc$ verifies $\displaystyle {\left\Vert{\nabla Rc}\right\Vert}_{\infty
 }<\infty $  and the injectivity radius is bounded below. This
 implies that the $\epsilon $ admissible radius is also bounded
 below. Hence we get the conclusion of corollary~\ref{7CZ2} in that case.\ \par 

\section{Applications.}

\begin{lem}
~\label{CF4}Let $\displaystyle t<2,$ if the weight $\alpha \in
 L^{\mu }$ with $\displaystyle \mu :=\frac{2t}{2-t},$ i.e. $\displaystyle
 \gamma (\alpha ,t)={\left\Vert{\alpha }\right\Vert}_{L^{\mu
 }(M)}^{\mu }=\int_{M}{\alpha ^{\frac{2t}{2-t}}}dv_{g}<\infty ,$ we have :\par 
\quad \quad \quad $\displaystyle \omega \in L^{2}_{p}(M)\Rightarrow \omega \in
 L^{t}_{p}(M,\alpha ).$
\end{lem}
\quad Proof.\ \par 
Young's inequality gives $\displaystyle {\left\Vert{fg}\right\Vert}_{L^{t}}\leq
 {\left\Vert{f}\right\Vert}_{L^{2}}{\left\Vert{g}\right\Vert}_{L^{q}}$
 with $\displaystyle \frac{1}{t}=\frac{1}{2}+\frac{1}{q},$ so
 let $\displaystyle \omega \in L^{2}_{p}(M),$ then, with $\displaystyle
 t<2,$ we get\ \par 
\quad \quad \quad $\displaystyle {\left({\int_{M}{\left\vert{\omega }\right\vert
 ^{t}\alpha dv_{g}}}\right)}^{1/t}\leq {\left({\int_{M}{\left\vert{\omega
 }\right\vert ^{2}dv_{g}}.}\right)}^{1/2}{\left({\int_{M}{\alpha
 ^{\frac{2t}{2-t}}}dv_{g}}\right)}^{\frac{2-t}{2t}}.$\ \par 
So if the weight $\displaystyle \alpha $ is such that $\displaystyle
 \ \gamma (\alpha ,t)=\int_{M}{\alpha ^{\frac{2t}{2-t}}}dv_{g}<\infty
 ,$ we are done. $\blacksquare $\ \par 
\quad For instance take any origin $\displaystyle 0\in M,\ M$ a complete
 riemannian manifold, and set $\rho (x):=d_{g}(0,x).$ We can
 choose a weight $\alpha ,$ function of $\rho ,\ \alpha (x):=f(\rho
 (x)),$ such that $\displaystyle \ \gamma (\alpha ,t)<\infty
 ,$ provided that $\displaystyle \alpha (x)$ goes to $0$ quickly
 enough at infinity.\ \par 
\quad Recall that $\displaystyle R(x)$ is the $\epsilon _{0}$ admissible
 radius at $\displaystyle x\in M.$\ \par 

\begin{cor}
~\label{CF5}Suppose that $\displaystyle (M,g)$ is a complete
 riemannian manifold ; let $\displaystyle r<2$ and choose a weight
 $\alpha \in L^{\infty }(M)$ verifying $\gamma (\alpha ,r)<\infty
 .$ Set $\displaystyle t:=\min (2,S_{2}(r)).$ If $\displaystyle
 t<2,$ take the weight $\alpha \in L^{\infty }(M)$ verifying
 also $\displaystyle \gamma (\alpha ,t)<\infty .$ Suppose we
 have condition (HL2,p).\par 
\quad Take $k$ big enough so that the threshold $\displaystyle S_{k}(r)\geq
 2,$ and set $\displaystyle w_{0}(x):=R(x)^{-2k},$ then for any
 $\omega \in L^{r}_{p}(M,w_{0}^{r})$ verifying $\displaystyle
 H\omega =0,$ for the orthogonal projection $H$ defined in corollary~\ref{CF9},
 there is a $\displaystyle u\in W^{2,r}_{p}(M,\alpha )\cap L^{t}_{p}(M,\alpha
 ),$ such that $\Delta u=\omega .$\par 
Moreover the solution $u$ is given linearly with respect to $\omega .$
\end{cor}
\quad Proof.\ \par 
Take $\displaystyle \omega \in L^{r}_{p}(M,w_{0}^{r}),$ with
 the choice of $\displaystyle w\equiv 1$ and $\displaystyle S_{k}(r)\geq
 2,$ the RSM theorem~\ref{CF3}, gives linear operators\ \par 
\quad \quad \quad $\displaystyle T\ :\ L^{r}_{p}(M,w_{0}^{r})\rightarrow L^{r}_{p}(M)\
 ;$ $\displaystyle A\ :\ L^{r}_{p}(M,w_{0}^{r})\rightarrow L^{2}_{p}(M),$\ \par 
such that\ \par 
\quad \quad \quad $\displaystyle v:=T\omega \in L^{r}(M)\cap L^{s}(M)\cap W^{2,r}(M)$
 verifies $\displaystyle \Delta v=\omega +\tilde \omega ,$\ \par 
with $\displaystyle s=S_{2}(r)$ and $\displaystyle \tilde \omega
 :=A\omega .$\ \par 
But\ \par 
\quad \quad \quad $\displaystyle v\in L^{t}(M)\Rightarrow v\in L^{t}(M,\alpha )$
 because $\displaystyle \alpha (x)\in L^{\infty }(M)$ is bounded :\ \par 
\quad \quad \quad ${\left\Vert{v}\right\Vert}_{L^{t}(M,\alpha )}^{t}=\int_{M}{\left\vert{v(x)}\right\vert
 ^{t}\alpha (x)dv(x)}\leq {\left\Vert{\alpha }\right\Vert}_{\infty
 }\int_{M}{\left\vert{v(x)}\right\vert ^{t}dv(x)}={\left\Vert{\alpha
 }\right\Vert}_{\infty }{\left\Vert{v}\right\Vert}_{L^{t}(M)}^{t}.$\ \par 
And the same $v\in L^{r}(M)\Rightarrow v\in L^{r}(M,\alpha ).$\ \par 
By corollary~\ref{CF9} if $\displaystyle H\omega =0$ then $\displaystyle
 H\tilde \omega =0.$\ \par 
Now we have $\displaystyle t_{k}:=S_{k}(r)\geq 2$ and we use
 the assumption (HL2,p) :\ \par 
\quad it gives the existence of a bounded linear operator $\displaystyle
 L\ :\ L^{2}_{p}(M)\rightarrow W^{2,2}_{p}(M)$ such that\ \par 
\quad \quad \quad $\Delta Lg=g,$ provided that $\displaystyle Hg=0,$\ \par 
by the spectral theorem (see, for instance, the proof of theorem
 5.10, p. 698 in Bueler~\cite{Bueler99}).\ \par 
So setting $f:=L\tilde \omega \in L^{2}_{p}(M)$ we have $\Delta
 f=\tilde \omega \in L^{2}_{p}(M).$\ \par 
\quad We set $\displaystyle u=v-f$ then $\displaystyle \Delta u=\omega
 +\tilde \omega -\tilde \omega =\omega .$ Let us see the estimates
 on $\displaystyle u.$\ \par 
Because $\displaystyle \gamma (\alpha ,r)<\infty ,$ we have by
 lemma~\ref{CF4}, $f\in L^{r}(M,\alpha ).$ If $\displaystyle
 t<2,$ we have also $\displaystyle \gamma (\alpha ,t)<\infty
 $ hence lemma~\ref{CF4} gives $f\in L^{t}(M,\alpha ).$\ \par 
\quad So in this case we have $\displaystyle u\in L^{r}(M,\alpha )\cap
 L^{t}(M,\alpha ).$\ \par 
\quad If $\displaystyle s\geq 2,$ then we have $\displaystyle t=2,\
 v\in L^{2}(M)$ by interpolation between $\displaystyle L^{r}(M)$
 and $\displaystyle L^{s}(M),$ so now we have $\displaystyle
 u\in L^{2}(M)\subset L^{t}(M,\alpha ).$\ \par 
\quad Because $\displaystyle v\in W^{2,r}(M),$ we get that $\displaystyle
 \nabla v,\ \nabla ^{2}v$ are also in $\displaystyle L^{r}(M)$
 so, the weight $\alpha $ being chosen bounded, we get that $\displaystyle
 \nabla v,\ \nabla ^{2}v$ are in $\displaystyle L^{r}(M,\ \alpha
 )$ so $\displaystyle v\in W^{2,r}(M,\ \alpha )$ We also have
 that $\displaystyle \nabla f,\ \nabla ^{2}f$ are in $\displaystyle
 L^{2}_{p}(M),$ hence because $\displaystyle \gamma (\alpha ,r)<\infty
 ,$ we get that $\displaystyle \nabla f,\ \nabla ^{2}f$ are in
 $\displaystyle L^{r}(M,\ \alpha ).$ This gives that $\displaystyle
 f\in W^{2,r}(M,\ \alpha )$ hence $\displaystyle u=v-f$$\in W^{2,r}(M,\
 \alpha ).$\ \par 
Hence in any case we get $u\in W^{2,r}_{p}(M,\alpha )\cap L^{t}_{p}(M,\alpha
 )$ and $\Delta u=\omega .$ $\blacksquare $\ \par 
\ \par 
\quad Now we shall use the linearity of our solution to get, by duality,
 results for exponents bigger than $\displaystyle 2.$ Take $\displaystyle
 r<2$ and $\displaystyle r'>2$ its conjugate.\ \par 
Let $T\ :\ L^{r}_{p}(M,w_{0}^{r})\rightarrow W_{p}^{2,r}(M)\subset
 L^{r}_{p}(M),\ A\ :\ L_{p}^{r}(M,w_{0}^{r})\rightarrow L_{p}^{2}(M)$
 be the linear operators, given by the RSM, such that\ \par 
\quad \quad \quad $\displaystyle \Delta T\omega =\omega +A\omega .$\ \par 
The hypothesis (HL2,p) gives the existence of a bounded linear
 operator $\displaystyle L\ :\ L^{2}_{p}(M)\rightarrow W^{2,2}_{p}(M)$
 such that\ \par 
\quad \quad \quad $\Delta L\tilde \omega =\tilde \omega ,$ provided that $\displaystyle
 H\tilde \omega =0\iff H\omega =0$ by corollary~\ref{CF9}.\ \par 
Hence, setting $\displaystyle C=LA\ :\ L_{p}^{r}(M,w_{0}^{r})\rightarrow
 W_{p}^{2,2}(M)$ we get\ \par 
\quad \quad \quad $\forall \omega \in L_{p}^{r}(M,w_{0}^{r}),\ \Delta (T-C)\omega
 =\omega .$\ \par 
We notice that $\forall \psi \in {\mathcal{D}}_{p}(M),$\ \par 
\quad \quad \quad $\displaystyle \Delta (T-C)\Delta \psi =\Delta \psi ,$\ \par 
just setting $\omega =\Delta \psi .$ This is possible because\ \par 
\quad \quad \quad $\forall \psi \in {\mathcal{D}}_{p}(M),\ \forall \varphi \in
 L^{2}_{p}(M),\ {\left\langle{H\Delta \psi ,\varphi }\right\rangle}={\left\langle{\Delta
 \psi ,H\varphi }\right\rangle}={\left\langle{\psi ,\Delta (H\varphi
 )}\right\rangle}=0,$\ \par 
where we used that $\Delta $ is essentially self adjoint, $M$
 being complete, and $\displaystyle \Delta (H\varphi )=0$ because
 $\displaystyle H\varphi $ is harmonic. So $\displaystyle H\Delta
 \psi =0$ and we can set $\omega =\Delta \psi $ because then
 $\displaystyle H\omega =0.$ Hence\ \par 
\quad \quad \quad \begin{equation}  (T-C)\Delta \psi =\psi +h,\label{CL13}\end{equation}\ \par 
with $h\in {\mathcal{H}}_{p}.$\ \par 
Now let $\displaystyle \varphi \in L^{2}_{p}(M)\cap L^{r'}_{p}(M)$
 and consider $\displaystyle u:=(T-C)^{*}\varphi ,$ the $*$ meaning
 the adjoint operator.\ \par 
This is meaningful because\ \par 
\quad \quad \quad $\displaystyle T^{*}\ :\ (W^{2,r}(M))'\supset L^{r'}(M)\rightarrow
 L^{r'}(M,w_{0}^{r})$\ \par 
and\ \par 
\quad \quad \quad $\displaystyle C^{*}\ :\ (W^{2,2}(M))'\supset L^{2}(M)\rightarrow
 L^{r'}(M,w_{0}^{r})$\ \par 
hence $\displaystyle u\in L^{r'}(M,w_{0}^{r}).$ We get\ \par 
\quad \quad \quad $\forall \psi \in {\mathcal{D}}(M)\cap L^{r}(M,w_{0}^{r}),\ {\left\langle{\Delta
 u,\psi }\right\rangle}_{L^{2}(M,w_{0}^{r})}={\left\langle{\Delta
 (T-C)^{*}\varphi ,\psi }\right\rangle}_{L^{2}(M,w_{0}^{r})}=$\ \par 
\quad \quad \quad \quad \quad \quad $\displaystyle =\int_{M}{\Delta {\left({(T-C)^{*}\varphi }\right)}\psi
 w_{0}^{r}dv_{g}}=\int_{M}{(T-C)^{*}\varphi \Delta (\psi w_{0}^{r})dv_{g}}={\left\langle{(T-C)^{*}\varphi
 ,\Delta (\psi w_{0}^{r})}\right\rangle}_{L^{2}(M)},$\ \par 
because $\Delta $ is essentially self adjoint and $\psi w_{0}^{r}$
 has compact support.\ \par 
\quad Hence by~(\ref{CL13})\ \par 
\quad \quad \quad $\displaystyle {\left\langle{\Delta u,\psi }\right\rangle}_{L^{2}(M,w_{0}^{r})}={\left\langle{\varphi
 ,(T-C)\Delta (\psi w_{0}^{r})}\right\rangle}_{L^{2}(M)}={\left\langle{\varphi
 ,\psi w_{0}^{r}+h}\right\rangle}_{L^{2}(M)}={\left\langle{\varphi
 ,\psi w_{0}^{r}}\right\rangle}_{L^{2}(M)},$\ \par 
provided that $\varphi \perp {\mathcal{H}},$ i.e. $\displaystyle
 H\varphi =0.$ Putting back the weight in the integral, we get\ \par 
\quad \quad \quad \begin{equation} {\left\langle{ \Delta u,\psi }\right\rangle}_{L^{2}(M,w_{0}^{r})}={\left\langle{\varphi
 ,\psi }\right\rangle}_{L^{2}(M,w_{0}^{r})}.\label{60}\end{equation}\ \par 
Now let $\psi '\in {\mathcal{D}}_{p}(M)$ and set $\displaystyle
 \psi :=\psi 'w_{0}^{-r}=\psi 'R(x)^{2kr}$ with $\displaystyle
 R(x)$ the $\epsilon $ admissible radius at the point $\displaystyle
 x\in M.$ We have seen that $\displaystyle \forall x\in M,\ R(x)>0$
 and we can smooth $\displaystyle R(x)$ to make it ${\mathcal{C}}^{\infty
 }(M)$ without changing the properties we used. For instance
 set $\tilde R(x):=\sum_{j\in {\mathbb{N}}}{\chi _{j}(x)R_{j}}$
 where $\lbrace \chi _{j}\rbrace _{j\in {\mathbb{N}}}$ is a partition
 of unity subordinated to our Vitali covering ${\mathcal{C}}_{\epsilon
 }=\lbrace B(x_{j},R_{j})\rbrace $ ; then the Lipschitz regularity
 of $\displaystyle R(x)$ contained in lemma~\ref{CF1} gives the
 existence of a constant $\displaystyle C>0$ depending only on
 $\displaystyle n,\epsilon $ such that $\displaystyle \forall
 x\in M,\ \frac{1}{C}R(x)\leq \tilde R(x)\leq CR(x).$\ \par 
So we have that $\psi \in {\mathcal{D}}_{p}(M)$ and\ \par 
\quad \quad \quad $\displaystyle {\left\langle{\Delta u,\psi }\right\rangle}_{L^{2}(M,w_{0}^{r})}={\left\langle{\varphi
 ,\psi '}\right\rangle}_{L^{2}(M)}\ ;\ $$\displaystyle \ {\left\langle{\varphi
 ,\psi }\right\rangle}_{L^{2}(M,w_{0}^{r})}={\left\langle{\varphi
 ,\psi '}\right\rangle}_{L^{2}(M)},$\ \par 
so~(\ref{60}) gives us\ \par 
\quad \quad \quad $\displaystyle {\left\langle{\Delta u,\psi '}\right\rangle}_{L^{2}(M)}={\left\langle{\varphi
 ,\psi '}\right\rangle}_{L^{2}(M)}.$\ \par 
This being true for any $\psi '\in {\mathcal{D}}_{p}(M)$ we get
 $\Delta u=\varphi $ in distributions sense, so we proved\ \par 

\begin{cor}
~\label{CF6}Suppose that $\displaystyle (M,g)$ is a complete
 riemannian manifold with (HL2,p) ; suppose we have $\displaystyle
 r<2$ with $\displaystyle k::S_{k}(r)\geq 2,$ setting $\displaystyle
 w_{0}(x):=R(x)^{-2k},$ for any $\varphi \in L^{2}_{p}(M)\cap
 L^{r'}_{p}(M),\ H\varphi =0$ we get\par 
\quad \quad \quad $u:=(T-C)^{*}\varphi ,\ u\in L^{r'}_{p}(M,w_{0}^{r})$ and $u$
 verifies $\displaystyle \Delta u=\varphi .$
\end{cor}
\quad Adding the hypothesis that the $\epsilon _{0}$ admissible radius
 is bounded below, we get more.\ \par 

\begin{cor}
~\label{62} Suppose that $\displaystyle (M,g)$ is a complete
 riemannian manifold and suppose the $\epsilon _{0}$ admissible
 radius verifies $\displaystyle \forall x\in M,\ R(x)\geq \delta
 >0,$ and suppose also hypothesis (HL2,p). Suppose we have $\displaystyle
 r<2$ with $\displaystyle k::S_{k}(r)\geq 2,$ setting $\displaystyle
 w_{0}(x):=R(x)^{-2k},$ for any $\varphi \in L^{2}_{p}(M)\cap
 L^{r'}_{p}(M),\ H\varphi =0$ we get\par 
\quad \quad \quad $\displaystyle u:=(T-C)^{*}\varphi ,\ u\in W^{2,r'}_{p}(M)$ and
 $u$ verifies $\displaystyle \Delta u=\varphi .$
\end{cor}
\quad Proof.\ \par 
Because $\displaystyle w_{0}^{r}(x)\geq 1,$ we get that $\displaystyle
 L^{r'}_{p}(M,w_{0}^{r})\subset L^{r'}_{p}(M)$ hence, applying
 this to $\displaystyle u,$ we get\ \par 
that\ \par 
\quad \quad \quad $\displaystyle u\in L^{r'}_{p}(M,w_{0}^{r})\Rightarrow u\in L^{r'}_{p}(M).$~\label{61}\
 \par 
Because the $\epsilon _{0}$ admissible radius verifies $\displaystyle
 \forall x\in M,\ R(x)\geq \delta >0,$ we have the classical
 Calderon Zygmund inequalities, corollary~\ref{7CZ2} :\ \par 
\quad \quad \quad $\displaystyle \forall r,\ 1<r<\infty ,\ {\left\Vert{u}\right\Vert}_{W^{2,r}(M)}\leq
 C_{1}{\left\Vert{u}\right\Vert}_{L^{r}(M)}+C_{2}{\left\Vert{\Delta
 u}\right\Vert}_{L^{r}(M)}.$\ \par 
The solution $u$ given by corollary~\ref{CF6}, $\displaystyle
 u:=(T-C)^{*}\varphi $ is in $\displaystyle u\in L^{r'}_{p}(M)$
 by~(\ref{61}). Because we have $\displaystyle \Delta u=\varphi
 \in L^{r'}_{p}(M),$ we get by CZI that $\displaystyle u\in $$\displaystyle
 W^{2,r'}_{p}(M),$ with control of the norms. $\blacksquare $\ \par 

\section{Non classical strong $\displaystyle L^{r}$ Hodge decomposition}
\quad We shall need :\ \par 

\begin{lem}
~\label{CL21}Let $\displaystyle r\leq 2$ and $\gamma \in W^{1,r}_{p+1}(M)\
 ;\ \beta \in W^{1,r}_{p-1}(M),\ h\in {\mathcal{H}}^{2}_{p}(M)$ then\par 
\quad \quad \quad $\displaystyle {\left\langle{d\gamma ,h}\right\rangle}={\left\langle{d^{*}\beta
 ,h}\right\rangle}=0.$
\end{lem}
\quad Proof.\ \par 
Because $\displaystyle h\in {\mathcal{H}}^{2}_{p},$ we have that
 $\displaystyle dh=d^{*}h=0$ by theorem 5.5, p. 697 in Bueler~\cite{Bueler99}.
 By the density of ${\mathcal{D}}_{k}(M)$ in $\displaystyle W^{1,r}_{k}(M)$
 which is always true in a complete riemannian manifold by theorem
 2.7, p. 13 in~\cite{Hebey96}, there is a sequence $\gamma _{k}\in
 {\mathcal{D}}_{p+1}(M)$ such that $\displaystyle {\left\Vert{\gamma
 -\gamma _{k}}\right\Vert}_{W^{1,r}(M)}\rightarrow 0$ and there
 is a sequence $\beta _{k}\in {\mathcal{D}}_{p-1}(M)$ such that
 $\displaystyle {\left\Vert{\beta -\beta _{k}}\right\Vert}_{W^{1,r}(M)}\rightarrow
 0.$\ \par 
By use of corollary~\ref{HCS44}, we have that $\displaystyle
 h\in {\mathcal{H}}^{2}_{p}\Rightarrow h\in {\mathcal{H}}^{r'}_{p}$
 because $\displaystyle r'>2,$ hence, because $\displaystyle
 d\gamma \in L^{r}_{p}(M),$\ \par 
\quad \quad \quad $\displaystyle {\left\langle{d\gamma ,h}\right\rangle}=\lim \
 _{k\rightarrow \infty }{\left\langle{d\gamma _{k},h}\right\rangle}=\lim
 \ _{k\rightarrow \infty }{\left\langle{\gamma _{k},d^{*}h}\right\rangle}=0,$\
 \par 
because $\displaystyle d^{*}$ is the formal adjoint of $d,\ \gamma
 _{k}\in {\mathcal{D}}_{p+1}(M)$ and $\displaystyle d^{*}h=0.$\ \par 
The same way we get $\displaystyle {\left\langle{d^{*}\beta ,h}\right\rangle}=0.$
 $\blacksquare $\ \par 

\begin{defin}
Let $\alpha $ be a weight on $M,$ we define the space $\displaystyle
 \tilde W^{2,r}_{p}(M,\alpha )$ to be\par 
\quad \quad \quad $\displaystyle \tilde W^{2,r}_{p}(M,\alpha ):=\lbrace u\in L^{r}_{p}(M,\alpha
 )::\Delta u\in L^{r}_{p}(M,\alpha )\rbrace $\par 
with the norm\par 
\quad \quad \quad $\displaystyle {\left\Vert{u}\right\Vert}_{\tilde W^{2,r}_{p}(M,\alpha
 )}:={\left\Vert{u}\right\Vert}_{L^{r}_{p}(M,\alpha )}+{\left\Vert{\Delta
 u}\right\Vert}_{L^{r}_{p}(M,\alpha )}.$
\end{defin}
\quad With just the hypothesis (HL2p) we get the Hodge decomposition.\ \par 

\begin{thm}
~\label{CL14}Let $\displaystyle (M,g)$ be a complete riemannian
 manifold. Let $\displaystyle r\leq 2$ and take a weight $\alpha
 \in L^{\infty }(M)$ such that $\gamma (\alpha ,r)<\infty $ ;
 with $\displaystyle k::S_{k}(r)\geq 2,$ set $\displaystyle w_{0}=R(x)^{-2k},$
 and suppose we have hypothesis (HL2,p). We have the direct decomposition
 given by linear operators :\par 
\quad \quad \quad $\ L^{r}_{p}(M,w_{0}^{r})={\mathcal{H}}_{p}^{2}\oplus \Delta
 (W^{2,r}_{p}(M,\alpha )).$\par 
With $\displaystyle r'>2,$ the conjugate exponent to $\displaystyle
 r,$ we have the weaker decomposition, still given by linear operators :\par 
\quad \quad \quad $L^{r'}_{p}(M)\cap L^{2}_{p}(M)={\mathcal{H}}_{p}^{2}\cap {\mathcal{H}}^{r'}_{p}+\Delta
 (\tilde W^{2,r'}_{p}(M)).$
\end{thm}
\quad Proof.\ \par 
Let $\omega \in L^{r}_{p}(M,w_{0}^{r})$ the remark~\ref{CL18}
 following the RSM with $\displaystyle w\equiv 1,\ w_{0}=R(x)^{-2k},$
 gives $\displaystyle u:=T\omega \in W^{2,r}_{p}(M),\ \tilde
 \omega :=A\omega \in L^{2}_{p}(M)$ such that $\Delta u=\omega
 +\tilde \omega .$ So we get\ \par 
\quad \quad \quad $\omega =\Delta u-\tilde \omega =\Delta u-(\tilde \omega -H\tilde
 \omega )-H\tilde \omega .$\ \par 
This is well defined because $\tilde \omega \in L^{2}_{p}(M)$
 and $H$ is the orthogonal projection from $\displaystyle L^{2}_{p}(M)$
 on ${\mathcal{H}}^{2}_{p}.$\ \par 
\quad Now $\displaystyle H(\tilde \omega -H\tilde \omega )=0$ hence
 by (HL2p) we get $f:=L(\tilde \omega -H\tilde \omega )$ solves
 $\Delta f=\tilde \omega -H\tilde \omega ,\ f\in W^{2,2}_{p}(M).$
 So we get\ \par 
\quad \quad \quad \begin{equation}  \omega =\Delta u-\tilde \omega =-H\tilde \omega
 +\Delta u-\Delta f,\label{CL16}\end{equation}\ \par 
with $H\tilde \omega \in {\mathcal{H}}^{2}_{p}.$ \ \par 
\quad This gives a first decomposition :\ \par 
\quad \quad \quad \begin{equation}  \omega =-H\tilde \omega +\Delta u-\Delta f,\label{CL22}\end{equation}\
 \par 
with $H\tilde \omega \in {\mathcal{H}}^{2}_{p}(M),\ u\in W^{2,r}_{p}(M)$
 and $\displaystyle f\in W^{2,2}_{p}(M).$\ \par 
With the weight $\alpha \in L^{\infty }(M)$ such that $\gamma
 (\alpha ,r)<\infty $ we have, by lemma~\ref{CF4}, $L^{2}_{p}(M)\subset
 L^{r}_{p}(M,\alpha ),$ hence the derivatives of $f$ up to second
 order are in $\displaystyle L^{2}_{p}(M)$ this implies that
 $\displaystyle f\in W^{2,r}_{p}(M,\alpha ).$ Because $\alpha
 $ is bounded, we also have $\displaystyle u\in W^{2,r}_{p}(M,\alpha ).$\ \par 
\quad It remains to set $\displaystyle v:=u-f\in W^{2,r}_{p}(M,\alpha
 )$ to get the decomposition. Because each step is linear, we
 get that this decomposition can be made linear with respect
 to $\omega .$\ \par 
\quad To get the uniqueness we consider the first decomposition~(\ref{CL22}) :\ \par 
\quad \quad \quad $\displaystyle \omega =h+\Delta (u-f)$ with $h\in {\mathcal{H}}^{2}_{p}$
 and $\displaystyle u\in W^{2,r}_{p}(M),\ f\in W^{2,2}_{p}(M).$\ \par 
If there is another one $\displaystyle \omega =h'+\Delta (u'-f')$
 then $\displaystyle 0=h-h'+\Delta (u-u'-(f-f'))\ ;$ so we have
 to show that\ \par 
\quad \quad \quad $\displaystyle 0=h+\Delta (u-f)$ with $h\in {\mathcal{H}}^{2}_{p}$
 and $\displaystyle u\in W^{2,r}_{p}(M),\ f\in W^{2,2}_{p}(M),$\ \par 
implies $\displaystyle h=0$ and $\displaystyle \Delta (u-f)=0.$\ \par 
\quad Now $\Delta u=d(d^{*}u)+d^{*}(du)=d\alpha +d^{*}\beta ,$ with
 $\alpha =d^{*}u\in W^{1,r}_{p+1}(M)$ and $\displaystyle \beta
 =du\in W^{1,r}_{p-1}(M).$ By lemma~\ref{CL21} we get $\displaystyle
 \ {\left\langle{d\alpha ,h}\right\rangle}+{\left\langle{d^{*}\beta
 ,h}\right\rangle}=0,$ so $\displaystyle \ {\left\langle{\Delta
 u,h}\right\rangle}=0.$ Exactly the same proof with $\displaystyle
 r=2$ gives $\displaystyle {\left\langle{\Delta f,h}\right\rangle}=0,$
 so, from $\displaystyle h+\Delta u-\Delta f=0,$ we get\ \par 
\quad \quad \quad $\displaystyle 0={\left\langle{h,h}\right\rangle}+{\left\langle{\Delta
 u,h}\right\rangle}+{\left\langle{\Delta f,h}\right\rangle}={\left\Vert{h}\right\Vert}_{L^{2}(M)},$\
 \par 
which implies $\displaystyle \Delta (u-f)=0$ and proves the uniqueness
 of this decomposition.\ \par 
\ \par 
\quad Now let $\displaystyle \omega \in L^{r'}_{p}(M)\cap L^{2}_{p}(M),$
 then we have\ \par 
\quad \quad \quad $\omega =H\omega +(\omega -H\omega )$ with $\displaystyle H(\omega
 -H\omega )=0.$\ \par 
We have that $H\omega \in {\mathcal{H}}^{2}_{p}(M)$ hence, by
 corollary~\ref{HCS44} because $\displaystyle \omega \in L^{2}_{p}(M),$
 we get that $\displaystyle H\omega \in {\mathcal{H}}^{r'}_{p}(M)$
 so $\displaystyle \tilde \omega :=\omega -H\omega \in L^{r'}_{p}(M)\cap
 L^{2}_{p}(M)$ and $H\tilde \omega =0.$ Now we have by corollary~\ref{CF6}
 a $\displaystyle u\in L^{r'}_{p}(M,w_{0}^{r})$ such that $\Delta
 u=\tilde \omega .$ Again this implies that $\displaystyle u\in
 L^{r'}_{p}(M)$ hence we have the decomposition\ \par 
\quad \quad \quad \begin{equation}  \forall \omega \in L^{r'}_{p}(M)\cap L^{2}_{p}(M),\
 \omega =H\omega +\Delta u=h+\Delta u,\label{CL19}\end{equation}\ \par 
with $\displaystyle h\in {\mathcal{H}}^{2}_{p}(M)\cap {\mathcal{H}}^{r'}_{p}(M)$
 and $\displaystyle u\in \tilde W^{2,r'}_{p}(M).$\ \par 
Because at each step we keep the linearity w.r.t. $\omega ,$
 we get that the decomposition is also linear w.r.t. $\omega
 .$ $\blacksquare $\ \par 
\ \par 
\quad There are two extreme cases done in the next corollaries.\ \par 

\begin{cor}
Suppose the $\epsilon _{0}$ admissible radius verifies $\displaystyle
 \forall x\in M,\ R(x)\geq \delta >0,$ and suppose also hypothesis
 (HL2,p). Take $\displaystyle r\leq 2$ and let the weight $\alpha
 \in L^{\infty }(M)$ be such that $\gamma (\alpha ,r)<\infty
 .$ Then we have the direct decomposition given by linear operators\par 
\quad \quad \quad $L^{r}_{p}(M)={\mathcal{H}}_{p}^{2}\oplus \Delta (W^{2,r}_{p}(M,\alpha )).$
\end{cor}
\quad Proof.\ \par 
In that case we have $\displaystyle \forall x\in M,\ 0<\delta
 \leq R(x)\leq 1$ hence $\displaystyle \ 1\leq w_{0}^{r}\leq
 \frac{1}{\delta ^{kr}}$ hence $\displaystyle L_{p}^{r}(M,w_{0}^{r})=L_{p}^{r}(M).$
 So we get this decomposition. $\blacksquare $\ \par 

\begin{cor}
~\label{HF1}Suppose the admissible radius verifies $\displaystyle
 \forall x\in M,\ R(x)\geq \delta >0,$ and suppose also hypothesis
 (HL2,p). Take $\displaystyle r'>2,$ then we have the direct
 decomposition given by linear operators\par 
\quad \quad \quad $\displaystyle L^{r'}_{p}(M)\cap L^{2}_{p}(M)={\mathcal{H}}_{p}^{2}\cap
 {\mathcal{H}}^{r'}_{p}\oplus \Delta (W^{2,r'}_{p}(M)).$
\end{cor}
\quad Proof.\ \par 
The classical CZI true in this case by corollary~\ref{7CZ2}, gives\ \par 
\quad \quad \quad $\displaystyle \forall r,\ 1<r<\infty ,\ {\left\Vert{u}\right\Vert}_{W^{2,r}(M)}\leq
 C_{1}{\left\Vert{u}\right\Vert}_{L^{r}(M)}+C_{2}{\left\Vert{\Delta
 u}\right\Vert}_{L^{r}(M)}.$\ \par 
So $\displaystyle u\in \tilde W^{2,r'}_{p}(M))\Rightarrow u\in
 W^{2,r'}_{p}(M))$ and we get the decomposition\ \par 
\quad \quad \quad $\displaystyle L^{r'}_{p}(M)\cap L^{2}_{p}(M)={\mathcal{H}}_{p}^{2}\cap
 {\mathcal{H}}^{r'}_{p}+\Delta (W^{2,r'}_{p}(M)).$\ \par 
\quad Now let us prove the uniqueness.\ \par 
We have the decomposition~(\ref{CL19})\ \par 
\quad \quad \quad $\displaystyle \forall \omega \in L^{r'}_{p}(M)\cap L^{2}_{p}(M),\
 \omega =h+\Delta u,$\ \par 
with $\displaystyle h\in {\mathcal{H}}^{2}_{p}(M)\cap {\mathcal{H}}^{r'}_{p}(M)$
 and $\displaystyle u\in W^{2,r'}_{p}(M).$\ \par 
\quad By (HL2,p) we have\ \par 
\quad \quad \quad $\displaystyle \exists v\in W^{2,2}_{p}(M)::\Delta v=\tilde \omega
 :=\omega -h.$\ \par 
But $\Delta v=\Delta u=\tilde \omega ,$ so if there is another
 such decomposition\ \par 
\quad \quad \quad $\omega =h'+\Delta u'=h'+\Delta v'$\ \par 
then\ \par 
\quad \quad \quad $\displaystyle 0=h-h'+\Delta (u-u')=h-h'+\Delta (v-v'),$\ \par 
Still with $\displaystyle v-v'\in W^{2,2}_{p}(M).$ So changing
 names we have\ \par 
\quad \begin{equation}  0=h+\Delta u=h+\Delta v\label{CL20}\end{equation}\ \par 
with $\displaystyle h\in {\mathcal{H}}^{2}_{p}(M)$ and $\displaystyle
 v\in W^{2,2}_{p}(M).$\ \par 
Again $\Delta v=d\alpha +d^{*}\beta $ with $\displaystyle \alpha
 =d^{*}v\in W^{1,2}_{p-1}(M)$ and $\displaystyle \beta =d^{*}v\in
 W^{1,2}_{p+1}(M)$ and by lemma~\ref{CL21} we get $\displaystyle
 {\left\langle{d\alpha ,h}\right\rangle}+{\left\langle{d^{*}\beta
 ,h}\right\rangle}=0,$ so $\displaystyle {\left\langle{\Delta
 v,h}\right\rangle}=0.$\ \par 
Hence $\displaystyle {\left\langle{\Delta u,h}\right\rangle}={\left\langle{\Delta
 v,h}\right\rangle}=0.$ But by~(\ref{CL20}) we have\ \par 
\quad \quad \quad $\displaystyle 0={\left\langle{h,h}\right\rangle}+{\left\langle{\Delta
 u,h}\right\rangle}$ so $\displaystyle {\left\Vert{h}\right\Vert}_{L^{2}(M)}=0\Rightarrow
 h=0$\ \par 
which ends the proof of uniqueness. $\blacksquare $\ \par 
\quad The admissible  radius verifies $\displaystyle \forall x\in M,\
 R(x)\geq \delta >0,$ if, for instance, the Ricci curvature of
 $M$ is bounded and the injectivity radius is strictly positive~\cite{HebeyHerzlich97}.\
 \par 
\quad We also have\ \par 

\begin{cor}
Let $\displaystyle r\leq 2,$ and, with $\displaystyle k::S_{k}(r)\geq
 2,$ set $\displaystyle w_{0}=R(x)^{-k}$ and suppose the riemannian
 volume is finite and hypothesis (HL2,p). We have the direct
 decomposition given by linear operators :\par 
\quad \quad \quad $L^{r}_{p}(M,w_{0}^{r})={\mathcal{H}}_{p}^{2}\oplus \Delta (W^{2,r}_{p}(M)).$
\end{cor}
\quad Here the weight $\alpha $ is no longer necessary because the
 volume being finite, if a form is in $\displaystyle L^{2}(M)$
 then it is already in $\displaystyle L^{r}(M).$ $\blacksquare $\ \par 

\begin{cor}
Let $\displaystyle r\leq 2$ and choose a weight $\alpha \in L^{\infty
 }(M)$ such that $\gamma (\alpha ,r)<\infty $ ; with $\displaystyle
 k::S_{k}(r)\geq 2,$ set $\displaystyle w_{0}=R(x)^{-k},$ and
 suppose we have hypothesis (HL2,p). We have the direct decompositions
 given by linear operators\par 
\quad \quad \quad $\displaystyle L^{r}_{p}(M,w_{0}^{r})={\mathcal{H}}_{p}^{2}\oplus
 d(W^{1,r}_{p}(M,\alpha ))\oplus d^{*}(W^{1,r}_{p}(M,\alpha )).$\par 
With $\displaystyle r'>2$ the conjugate exponent of $\displaystyle
 r,$ and adding the hypothesis that the $\epsilon _{0}$ admissible
 radius is bounded below, we get\par 
\quad \quad \quad $\displaystyle L^{r'}_{p}(M)\cap L^{2}_{p}(M)={\mathcal{H}}_{p}^{2}\cap
 {\mathcal{H}}^{r'}_{p}\oplus d(W^{1,r'}_{p}(M))\oplus d^{*}(W^{1,r'}_{p}(M)).$
\end{cor}
\quad Proof.\ \par 
For the first part, we have, by~(\ref{CL22}) : $\displaystyle
 \forall \omega \in L^{r}_{p}(M,w_{0}^{r}),$\ \par 
\quad \quad \quad $\displaystyle \omega =-H\tilde \omega +\Delta u-\Delta f,$\ \par 
with $H\tilde \omega \in {\mathcal{H}}^{2}_{p}(M),\ u\in W^{2,r}_{p}(M)$
 and $\displaystyle f\in W^{2,2}_{p}(M).$ Again\ \par 
\quad \quad \quad $\Delta u=d\gamma +d^{*}\beta ,$ with $\displaystyle \gamma \in
 W^{1,r}_{p-1}(M)$ and $\displaystyle \beta \in W^{1,r}_{p+1}(M),$\ \par 
and\ \par 
\quad \quad \quad $\Delta f=d\gamma '+d^{*}\beta ',$ with $\displaystyle \gamma
 '\in W^{1,2}_{p-1}(M)$ and $\displaystyle \beta '\in W^{1,2}_{p+1}(M),$\ \par 
Hence\ \par 
\quad \quad \quad $\displaystyle \omega =h+d(\gamma -\gamma ')+d^{*}(\beta -\beta ').$\ \par 
With the weight $\alpha $ we get $\displaystyle \gamma \in W^{1,r}_{p-1}(M)\Rightarrow
 \gamma \in W^{1,r}_{p-1}(M,\alpha )$ and the same for $\beta .$ And also\ \par 
$\gamma '\in W^{1,2}_{p-1}(M)\Rightarrow \gamma '\in W^{1,r}_{p-1}(M,\alpha
 )$ and the same for $\beta '.$ So, setting $\mu :=\gamma -\gamma
 ',\ \delta =\beta -\beta ',$ we have the decomposition\ \par 
\quad \quad \quad $\displaystyle \omega \in L^{r}_{p}(M,w_{0}^{r})\Rightarrow \omega
 =h+d\mu +d^{*}\delta ,$\ \par 
with $h\in {\mathcal{H}}^{2}_{p}(M)\cap {\mathcal{H}}^{r}_{p}(M,\alpha
 ),\ \mu \in W^{1,r}_{p-1}(M,\alpha ),\ \delta \in W^{1,r}_{p+1}(M,\alpha
 ).$\ \par 
\ \par 
\quad For the uniqueness, suppose that\ \par 
\quad \quad \quad $\displaystyle 0=h+d(\gamma -\gamma ')+d^{*}(\beta -\beta '),$\ \par 
by use of lemma~\ref{CL21}, we get $\displaystyle {\left\langle{d\gamma
 ,h}\right\rangle}+{\left\langle{d^{*}\beta ,h}\right\rangle}=0$
 and also $\displaystyle {\left\langle{d\gamma ',h}\right\rangle}+{\left\langle{d^{*}\beta
 ',h}\right\rangle}=0,$ so $\displaystyle h=0.$ So we have\ \par 
\quad \quad \quad $\displaystyle 0=d(\gamma -\gamma ')+d^{*}(\beta -\beta ').$\ \par 
This implies that\ \par 
\quad \quad \quad \begin{equation}  d\gamma +d^{*}\beta =d\gamma '+d^{*}\beta ',\label{CL23}\end{equation}\
 \par 
hence\ \par 
\quad \quad \quad $\displaystyle d\gamma +d^{*}\beta \in L^{r}_{p}(M)\cap L^{2}_{p}(M)\
 ;\ d\gamma '+d^{*}\beta '\in L^{r}_{p}(M)\cap L^{2}_{p}(M),$\ \par 
because \ \par 
\quad \quad \quad $\displaystyle d\gamma +d^{*}\beta \in L^{r}_{p}(M)$ and $\displaystyle
 d\gamma '+d^{*}\beta '\in L^{2}_{p}(M).$\ \par 
Now take $\varphi \in {\mathcal{D}}_{p}(M),$ because (HL2,p)
 is true we have the $\displaystyle L^{2}$ decomposition :\ \par 
\quad \quad \quad $\displaystyle \varphi =H\varphi +d\gamma \mu +d^{*}\delta $
 with $\mu ,\delta \in W^{1,2}(M).$\ \par 
We have\ \par 
\quad \quad \quad $\displaystyle {\left\langle{d(\gamma -\gamma '),\varphi }\right\rangle}={\left\langle{d(\gamma
 -\gamma '),H\varphi +d\mu +d^{*}\delta }\right\rangle}\ ;$\ \par 
by use of lemma~\ref{CL21}, we get $\displaystyle {\left\langle{d(\gamma
 -\gamma '),H\varphi }\right\rangle}=0.$ By density we have $\displaystyle
 \mu =\lim \ _{k\rightarrow \infty }\mu _{k},\ \gamma _{k}\in
 {\mathcal{D}}_{p-1}$ and $\displaystyle \delta =\lim \ _{k\rightarrow
 \infty }\delta _{k},\ \delta _{k}\in {\mathcal{D}}_{p+1},$ the
 convergence being in $\displaystyle W^{1,2}(M),$ so $\displaystyle
 d\mu =\lim \ _{k\rightarrow \infty }d\mu _{k}$ and $\displaystyle
 d^{*}\delta =\lim \ _{k\rightarrow \infty }d^{*}\delta _{k}$
 in $\displaystyle L^{2}_{p}(M).$ So we get\ \par 
\quad \quad \quad $\displaystyle {\left\langle{d(\gamma -\gamma '),d\mu +d^{*}\delta
 }\right\rangle}\ =\lim \ _{k\rightarrow \infty }{\left\langle{d(\gamma
 -\gamma '),d\mu _{k}+d^{*}\delta _{k}}\right\rangle}.$\ \par 
But\ \par 
\quad \quad \quad $\displaystyle \forall k\in {\mathbb{N}},\ {\left\langle{d(\gamma
 -\gamma '),d^{*}\delta _{k}}\right\rangle}={\left\langle{(\gamma
 -\gamma '),d^{*2}\delta _{k}}\right\rangle}=0$\ \par 
because $\displaystyle d^{*}$ is the formal adjoint of $d$ and
 $\displaystyle d^{*}\delta _{k}$ has compact support and $\displaystyle
 d^{*2}=0.$ So\ \par 
\quad \quad \quad $\displaystyle {\left\langle{d(\gamma -\gamma '),\varphi }\right\rangle}=\lim
 \ _{k\rightarrow \infty }{\left\langle{d(\gamma -\gamma '),d\mu
 _{k}}\right\rangle}.$\ \par 
With~(\ref{CL23}) we get\ \par 
\quad \quad \quad $\displaystyle \ \forall k\in {\mathbb{N}},\ {\left\langle{d(\gamma
 -\gamma '),d\mu _{k}}\right\rangle}-{\left\langle{d^{*}(\beta
 -\beta '),d\mu _{k}}\right\rangle}=0,$\ \par 
and\ \par 
\quad \quad \quad $\displaystyle \forall k\in {\mathbb{N}},\ {\left\langle{d^{*}(\beta
 -\beta '),d\mu _{k}}\right\rangle}=0,$\ \par 
because $\displaystyle d^{*}$ is the formal adjoint of $d,\ d\gamma
 _{k}$ has compact support and $\displaystyle d^{2}=0.$ So\ \par 
\quad \quad \quad $\displaystyle \forall k\in {\mathbb{N}},\ {\left\langle{d(\gamma
 -\gamma '),d\mu _{k}}\right\rangle}=0,$\ \par 
which gives\ \par 
\quad \quad \quad $\displaystyle {\left\langle{d(\gamma -\gamma '),\varphi }\right\rangle}=\lim
 \ _{k\rightarrow \infty }{\left\langle{d(\gamma -\gamma '),d\mu
 _{k}}\right\rangle}=0,$\ \par 
and this being true for any $\displaystyle \varphi \in {\mathcal{D}}_{p}(M),$
 we get $\displaystyle d(\gamma -\gamma ')=0\ ;$ this gives with~(\ref{CL23})
 $\displaystyle d^{*}(\beta -\beta ')=0.$\ \par 
\ \par 
\quad For the second case we already have	, by theorem~\ref{CL14} plus
 CZI given by corollary~\ref{7CZ2},  $\displaystyle \omega =H\omega
 +\Delta u$ with $\displaystyle u\in W^{2,r}(M).$ Now 	$\Delta
 u=d(d^{*}u)+d^{*}(du)=d\gamma +d^{*}\beta ,$ with $\gamma =d^{*}u\in
 W^{1,r'}_{p+1}(M)$ and $\displaystyle \beta =du\in W^{1,r'}_{p-1}(M).$
 This gives the decomposition.\ \par 
\quad For the uniqueness the proof is exactly the same as above, so
 we are done. $\blacksquare $\ \par 

\subsection{Non classical weak $\displaystyle L^{r}$ Hodge decomposition.~\label{HCF37}}

      Now we shall need another hypothesis :\ \par 
(HWr) if the space ${\mathcal{D}}_{p}(M)$ is dense in $\displaystyle
 W^{2,r}_{p}(M).$\ \par 
We already know that (HWr) is true if :\ \par 
\quad $\bullet $ either : the injectivity radius is strictly positive
 and the Ricci curvature is bounded~[\cite{Hebey96} theorem 2.8, p. 12].\ \par 
\quad $\bullet $ or : $M$ is geodesically complete with a bounded curvature
 tensor~[\cite{GuneysuPigola} theorem 1.1 p.3].\ \par 
\quad We have a non classical weak $\displaystyle L^{r}$ Hodge decomposition
 theorem :\ \par 

\begin{thm}
~\label{HC35}Suppose that $\displaystyle (M,g)$ is a complete
 riemannian manifold, fix $\displaystyle r\leq 2$ and choose
 a bounded weight $\alpha $ with $\displaystyle \gamma (\alpha
 ,r)<\infty .$\par 
Take $k$ with $\displaystyle S_{k}(r)\geq 2,$ and set the weight
 $\displaystyle w_{0}:=R(x)^{-2k}.$ Suppose we have (HL2,p) and
 (HW2) ; then\par 
\quad \quad \quad $L^{r}_{p}(M,\alpha )={\mathcal{H}}_{p}^{r}(M,\alpha )\oplus
 {\overline{\Delta ({\mathcal{D}}_{p}(M))}},$\par 
the closure being taken in $L^{r}(M,\alpha ).$
\end{thm}
\quad Proof.\ \par 
Take $\omega \in L^{r}_{p}(M,\alpha ).$ By density there is a
 $\omega _{\epsilon }\in {\mathcal{D}}_{p}(M)$ such that ${\left\Vert{\omega
 -\omega _{\epsilon }}\right\Vert}_{L^{r}(M,\alpha )}<\epsilon .$\ \par 
Then, because $\omega _{\epsilon }\in {\mathcal{D}}_{p}(M),$
 we have $\omega _{\epsilon }\in L^{r}_{p}(M,w_{0}^{r})$ hence by RSM :\ \par 
\quad \quad \quad $\displaystyle \forall s\geq r,\ \exists v_{\epsilon }\in L^{r}_{p}(M)\cap
 L^{s_{1}}_{p}(M)::\Delta v_{\epsilon }=\omega _{\epsilon }+\tilde
 \omega _{\epsilon },$\ \par 
with $s_{1}:=S_{2}(r),\ \tilde \omega _{\epsilon }\in L^{s}_{p}(M).$
 Moreover, because $\displaystyle \omega _{\epsilon }$ is of
 compact support, so are $\displaystyle v_{\epsilon }$ and $\tilde
 \omega _{\epsilon }.$\ \par 
\quad Taking $\displaystyle s=2,$ by (HL2,p) there is a $\displaystyle
 f_{\epsilon }\in W^{2,2}_{p}(M)::\Delta f_{\epsilon }=\tilde
 \omega _{\epsilon }-H\tilde \omega _{\epsilon }.$\ \par 
By (HW2) there is a $g_{\epsilon }\in {\mathcal{D}}_{p}(M)::{\left\Vert{f_{\epsilon
 }-g_{\epsilon }}\right\Vert}_{W^{2,2}(M)}<\epsilon $ and this implies\ \par 
\quad \quad \quad $\displaystyle {\left\Vert{\Delta f_{\epsilon }-\Delta g_{\epsilon
 }}\right\Vert}_{L^{2}(M)}<\epsilon .$\ \par 
Now we set $\displaystyle u_{\epsilon }:=v_{\epsilon }-g_{\epsilon
 },$ then $\displaystyle u_{\epsilon }$ is of compact support and we have\ \par 
\quad $\displaystyle \Delta u_{\epsilon }=\Delta v_{\epsilon }-\Delta
 g_{\epsilon }=\Delta v_{\epsilon }-\Delta f_{\epsilon }+(\Delta
 f_{\epsilon }-\Delta g_{\epsilon })=\omega _{\epsilon }+\tilde
 \omega _{\epsilon }-\tilde \omega _{\epsilon }+H\tilde \omega
 _{\epsilon }+E_{\epsilon }=\omega _{\epsilon }+H\tilde \omega
 _{\epsilon }+E_{\epsilon },$\ \par 
where we set $\displaystyle E_{\epsilon }:=\Delta f_{\epsilon
 }-\Delta g_{\epsilon }.$\ \par 
\quad So we get\ \par 
\quad \quad \quad $\displaystyle \omega =-H\tilde \omega _{\epsilon }+\Delta u_{\epsilon
 }+(\omega -\omega _{\epsilon })+E_{\epsilon }.$\ \par 
Because $\gamma (\alpha ,r)<\infty ,$ we get $\displaystyle {\left\Vert{E_{\epsilon
 }}\right\Vert}_{L^{r}(M,\alpha )}\leq C{\left\Vert{\Delta f_{\epsilon
 }-\Delta g_{\epsilon }}\right\Vert}_{L^{2}(M)}<C\epsilon .$
 For the same reason we have $H\tilde \omega _{\epsilon }\in
 {\mathcal{H}}_{p}^{2}(M)\subset {\mathcal{H}}_{p}^{r}(M,\alpha
 ),$ so we get $\displaystyle \omega \in {\mathcal{H}}_{p}^{r}(M,\alpha
 )+{\overline{\Delta ({\mathcal{D}}_{p}(M))}},$ the closure being
 taken in $\displaystyle L^{r}(M,\alpha ).$\ \par 
\quad For the uniqueness we proceed as before. We have to show that
 if $\displaystyle 0=\lim _{k\rightarrow \infty }(h_{k}+\Delta
 u_{k})$ with $\displaystyle h_{k}\in {\mathcal{H}}_{p}^{2}(M)\subset
 {\mathcal{H}}_{p}^{r}(M,\alpha )$ and $u_{k}\in {\mathcal{D}}_{p}(M),$
 the convergence in $\displaystyle L^{r}(M,\alpha ),$ then $\displaystyle
 \lim _{k\rightarrow \infty }h_{k}=0$ and $\displaystyle \lim
 _{k\rightarrow \infty }\Delta u_{k}=0.$\ \par 
We have $\displaystyle \Delta u_{k}=d\gamma _{k}+d^{*}\beta _{k},$
 with $\displaystyle \gamma _{k}=d^{*}u_{k}\in {\mathcal{D}}_{p+1}(M),$
 and $\displaystyle \beta _{k}=du_{k}\in {\mathcal{D}}_{p-1}(M).$
 So we can apply lemma~\ref{CL21} to get\ \par 
\quad \quad \quad $\displaystyle \forall k,\ {\left\langle{h_{k},d\gamma _{k}}\right\rangle}={\left\langle{h_{k},d^{*}\beta
 _{k}}\right\rangle}=0,$\ \par 
hence\ \par 
\quad \quad \quad $\displaystyle \lim _{k\rightarrow \infty }{\left\langle{h_{k},h_{k}}\right\rangle}=0\Rightarrow
 \lim _{k\rightarrow \infty }h_{k}=0$ and hence $\displaystyle
 \lim _{k\rightarrow \infty }\Delta u_{k}=0.$ $\blacksquare $\ \par 
\quad We also have a weak $\displaystyle L^{r}$ Hodge decomposition
 without hypothesis (HWr) :\ \par 

\begin{thm}
Suppose that $\displaystyle (M,g)$ is a complete riemannian manifold
 and suppose we have (HL2,p). Fix $\displaystyle r<2$ and take
 a weight $\alpha $ verifying $\gamma (\alpha ,r)<\infty .$ Then we have\par 
\quad \quad \quad $L^{r}_{p}(M,\alpha )={\mathcal{H}}_{p}^{r}(M,\alpha )\oplus
 {\overline{d({\mathcal{D}}_{p-1}(M))}}\oplus {\overline{d^{*}({\mathcal{D}}_{p+1}(M))}},$\par
 
the closures being taken in $L^{r}(M,\alpha ).$
\end{thm}
\quad Proof.\ \par 
We start exactly the same way as for theorem~\ref{HC35} to have\ \par 
\quad \quad \quad $\displaystyle v_{\epsilon }\in L^{r}_{p}(M)\cap L^{s_{1}}_{p}(M)::\Delta
 v_{\epsilon }=\omega _{\epsilon }+\tilde \omega _{\epsilon },$\ \par 
and\ \par 
\quad \quad \quad $\displaystyle f_{\epsilon }\in W^{2,2}_{p}(M)::\Delta f_{\epsilon
 }=\tilde \omega _{\epsilon }-H\tilde \omega _{\epsilon }.$\ \par 
Now we set directly $\displaystyle u_{\epsilon }:=v_{\epsilon
 }-f_{\epsilon }\Rightarrow \Delta u_{\epsilon }=\omega _{\epsilon
 }+H\tilde \omega _{\epsilon }.$ The point here is that $\displaystyle
 u_{\epsilon }$ is not of compact support because $\displaystyle
 f_{\epsilon }$ is not.\ \par 
Nevertheless we have :\ \par 
\quad \quad \quad \begin{equation}  \omega =-H\tilde \omega _{\epsilon }+(H\tilde
 \omega _{\epsilon }+\omega _{\epsilon })+(\omega -\omega _{\epsilon
 })=-H\tilde \omega _{\epsilon }+\Delta u_{\epsilon }+(\omega
 -\omega _{\epsilon }).\label{HC36}\end{equation}\ \par 
But we can approximate $\displaystyle d^{*}u_{\epsilon }$ by
 $\gamma _{\epsilon }\in {\mathcal{D}}(M)$ in $\displaystyle
 W^{1,2}(M),\ $and $\displaystyle du_{\epsilon }$ by $\beta _{\epsilon
 }\in {\mathcal{D}}(M)$ in $\displaystyle W^{1,2}(M),\ $and this
 is always possible by theorem 2.7, p. 13 in~\cite{Hebey96}. So we have\ \par 
\quad \quad \quad $\displaystyle {\left\Vert{d^{*}u_{\epsilon }-\gamma _{\epsilon
 }}\right\Vert}_{W^{1,2}(M)}<\epsilon ,\ {\left\Vert{du_{\epsilon
 }-\beta _{\epsilon }}\right\Vert}_{W^{1,2}(M)}<\epsilon .$\ \par 
And this implies\ \par 
\quad \quad \quad $\displaystyle {\left\Vert{\Delta u_{\epsilon }-d\gamma _{\epsilon
 }-d^{*}\beta _{\epsilon }}\right\Vert}_{L^{2}_{p}(M)}\leq 2\epsilon
 \Rightarrow {\left\Vert{\Delta u_{\epsilon }-d\gamma _{\epsilon
 }-d^{*}\beta _{\epsilon }}\right\Vert}_{L^{r}_{p}(M,\alpha )}\leq
 2C\epsilon ,$\ \par 
because $\displaystyle \gamma (\alpha ,r)<\infty .$ As above
 we have $H\tilde \omega _{\epsilon }\in {\mathcal{H}}_{p}^{r}(M,\alpha
 )$ so putting all this in~(\ref{HC36}) we get\ \par 
\quad \quad \quad $\omega \in {\mathcal{H}}_{p}^{r}(M,\alpha )+{\overline{d({\mathcal{D}}_{p-1}(M))}}+{\overline{d^{*}({\mathcal{D}}_{p+1}(M))}},$\
 \par 
the closure being taken in $\displaystyle L^{r}_{p}(M,\alpha ).$\ \par 
The proof of the uniqueness is exactly as in the proof of theorem~\ref{HC35},
 so we are done. $\blacksquare $\ \par 

\begin{rem}
It seems not "geometrically natural" to take the closure of 
 $d^{*}({\mathcal{D}}_{p+1}(M))$ with respect to $\displaystyle
 L^{r}_{p}(M,\alpha )$ because here the adjoint of $d,\ d^{*},$
 is taken with respect to the volume measure without any weight.
 Nevertheless this is "analytically" correct and we get nothing
 more here. This is why we call the two previous results "non classical".
\end{rem}
\quad For the case $\displaystyle r>2$ we need a stronger hypothesis,
 namely that the $\epsilon _{0}$ admissible radius is bounded
 below. Then we get a \emph{classical weak Hodge decompositions}
 for $\displaystyle r\geq 2.$\ \par 

\begin{thm}
~\label{9CW0}Suppose that $\displaystyle (M,g)$ is a complete
 riemannian manifold and suppose the $\epsilon _{0}$ admissible
 radius verifies $\displaystyle \forall x\in M,\ R(x)\geq \delta
 >0,$ suppose (HWr) and suppose also hypothesis (HL2,p). Fix
 $\displaystyle r\geq 2,$ then we have\par 
\quad \quad \quad $\displaystyle L^{r}_{p}(M)={\mathcal{H}}_{p}^{r}(M)\oplus {\overline{\Delta
 ({\mathcal{D}}_{p}(M))}}.$\par 
Without (HWr) we still get\par 
\quad \quad \quad $L^{r}_{p}(M)={\mathcal{H}}_{p}^{r}(M)\oplus {\overline{d({\mathcal{D}}_{p-1}(M))}}\oplus
 {\overline{d^{*}({\mathcal{D}}_{p+1}(M))}}.$\par 
All the closures being taken in $L^{r}(M).$
\end{thm}
\quad Proof.\ \par 
Take $\displaystyle \omega \in L^{r}_{p}(M),$ then by density
 there is a $\omega _{\epsilon }\in {\mathcal{D}}_{p}(M)$ such
 that $\displaystyle {\left\Vert{\omega -\omega _{\epsilon }}\right\Vert}_{L^{r}(M)}\leq
 \epsilon .$ This implies, because $\displaystyle r>2$ and $\omega
 _{\epsilon }$ is compactly supported, that $\displaystyle \omega
 _{\epsilon }\in L^{r}(M)\cap L^{2}(M).$ So we have $\displaystyle
 H\omega _{\epsilon }\in L^{2}(M)\Rightarrow H\omega _{\epsilon
 }\in L^{r}(M)$ by corollary~\ref{HCS44}.\ \par 
\quad So let $\varphi _{\epsilon }:=\omega _{\epsilon }-H\omega _{\epsilon
 }\in L^{r}(M)\cap L^{2}(M),$ we have $\displaystyle H\varphi
 _{\epsilon }=0$ hence by corollary~\ref{62} we have\ \par 
\quad \quad \quad $\displaystyle \exists u_{\epsilon }\in L^{r}(M,w_{0}^{r})\cap
 W^{2,r}_{p}(M)::\Delta u_{\epsilon }=\varphi _{\epsilon }.$\ \par 
So if we have (HWr) then $\exists v_{\epsilon }\in {\mathcal{D}}_{p}(M)$
 such that $\displaystyle {\left\Vert{u_{\epsilon }-v_{\epsilon
 }}\right\Vert}_{W^{2,r}(M)}<\epsilon $ and this implies\ \par 
\quad \quad \quad $\displaystyle {\left\Vert{\Delta u_{\epsilon }-\Delta v_{\epsilon
 }}\right\Vert}_{L^{r}(M)}<\epsilon .$\ \par 
Now we can write\ \par 
\quad \quad \quad $\displaystyle \omega =H\omega _{\epsilon }+(\omega -\omega _{\epsilon
 })+(\omega _{\epsilon }-H\omega _{\epsilon })=H\omega _{\epsilon
 }+(\omega -\omega _{\epsilon })+\Delta u_{\epsilon }=H\omega
 _{\epsilon }+(\omega -\omega _{\epsilon })+\Delta v_{\epsilon
 }+(\Delta u_{\epsilon }-\Delta v_{\epsilon }).$\ \par 
The term $\displaystyle E_{\epsilon }:=(\omega -\omega _{\epsilon
 })+(\Delta u_{\epsilon }-\Delta v_{\epsilon })$ is an error
 term small in $\displaystyle L^{r}(M)$ so we get\ \par 
\quad \quad \quad $\displaystyle \omega =H\omega _{\epsilon }+\Delta v_{\epsilon
 }+E_{\epsilon },$ with $H\omega _{\epsilon }\in L^{r}(M)\cap
 L^{2}(M),\ v_{\epsilon }\in {\mathcal{D}}_{p}(M),\ {\left\Vert{E_{\epsilon
 }}\right\Vert}_{L^{r}(M)}<2\epsilon .$\ \par 
\quad So we have the decomposition :\ \par 
\quad \quad \quad $\displaystyle L^{r}_{p}(M)={\mathcal{H}}_{p}^{r}(M)\oplus {\overline{\Delta
 ({\mathcal{D}}_{p}(M))}}.$\ \par 
the closures being taken in $L^{r}(M).$\ \par 
\quad Without (HWr) we approximate $\displaystyle d^{*}u_{\epsilon
 }$ by $\gamma _{\epsilon }\in {\mathcal{D}}(M)$ in $\displaystyle
 W^{1,r}(M),\ $and $\displaystyle du_{\epsilon }$ by $\beta _{\epsilon
 }\in {\mathcal{D}}(M)$ in $\displaystyle W^{1,r}(M),\ $and this
 is always possible by theorem 2.7, p. 13 in~\cite{Hebey96}. So we have\ \par 
\quad \quad \quad $\displaystyle {\left\Vert{d^{*}u_{\epsilon }-\gamma _{\epsilon
 }}\right\Vert}_{W^{1,r}(M)}<\epsilon ,\ {\left\Vert{du_{\epsilon
 }-\beta _{\epsilon }}\right\Vert}_{W^{1,r}(M)}<\epsilon .$\ \par 
And this implies\ \par 
\quad \quad \quad $\displaystyle {\left\Vert{\Delta u_{\epsilon }-d\gamma _{\epsilon
 }-d^{*}\beta _{\epsilon }}\right\Vert}_{L^{r}_{p}(M)}\leq 2\epsilon .$\ \par 
So we have\ \par 
\quad \quad \quad $\displaystyle \omega =H\omega _{\epsilon }+(\omega -\omega _{\epsilon
 })+\Delta u_{\epsilon }=H\omega _{\epsilon }+(\omega -\omega
 _{\epsilon })+d\gamma _{\epsilon }+d^{*}\beta _{\epsilon }+(\Delta
 u_{\epsilon }-d\gamma _{\epsilon }-d^{*}\beta _{\epsilon }).$\ \par 
The term $\displaystyle E_{\epsilon }:=(\omega -\omega _{\epsilon
 })+(\Delta u_{\epsilon }-d\gamma _{\epsilon }-d^{*}\beta _{\epsilon
 }).$ is an error term small in $\displaystyle L^{r}(M)$ so we get\ \par 
\quad \quad \quad $\displaystyle \omega =H\omega _{\epsilon }+d\gamma _{\epsilon
 }+d^{*}\beta _{\epsilon }+E_{\epsilon },$\ \par 
with\ \par 
\quad \quad \quad $H\omega _{\epsilon }\in L^{r}(M)\cap L^{2}(M),\ \gamma _{\epsilon
 }\in {\mathcal{D}}_{p+1}(M),\ \beta _{\epsilon }\in {\mathcal{D}}_{p-1}(M),\
 {\left\Vert{E_{\epsilon }}\right\Vert}_{L^{r}(M)}<2\epsilon .$\ \par 
\quad So we have the decomposition :\ \par 
\quad \quad \quad $\displaystyle L^{r}_{p}(M)={\mathcal{H}}_{p}^{r}(M)\oplus {\overline{d({\mathcal{D}}_{p-1}(M))}}\oplus
 {\overline{d^{*}({\mathcal{D}}_{p+1}(M))}},$\ \par 
the closures being taken in $L^{r}(M).$\ \par 
\quad The proof of the uniqueness is a slight modification of the proof
 of corollary~\ref{HF1}, so we are done. $\blacksquare $\ \par 

\begin{rem}
By theorem 1.3 in Hebey~\cite{Hebey96}, we have that the harmonic
 radius $\displaystyle r_{H}(1+\epsilon ,\ 2,0)$ is bounded below
 if the Ricci curvature $\displaystyle Rc$ verifies $\displaystyle
 {\left\Vert{\nabla Rc}\right\Vert}_{\infty }<\infty $  and the
 injectivity radius is bounded below. This implies that the $\epsilon
 $ admissible radius is also bounded below.\par 
Moreover if we add the hypothesis that the Ricci curvature $\displaystyle
 Rc$ verifies $\ \exists \delta \in {\mathbb{R}}::Rc\geq \delta
 $ then by Proposition 2.10 in Hebey~\cite{Hebey96}, we have hypothesis (HWr).
\end{rem}
\ \par 

\bibliographystyle{/usr/local/texlive/2013/texmf-dist/bibtex/bst/base/plain}

\begin{thebibliography}{10}

\bibitem{AmarSt13}
E.~Amar.
\newblock The raising steps method. {A}pplication to the $\bar\partial$
  equation in {S}tein manifolds.
\newblock {\em J. Geometric Analysis}, 26(2):898--913, 2016.

\bibitem{HodgeCompact15}
Eric Amar.
\newblock The raising steps method. {A}pplications to the ${L}^{r}$ {H}odge
  theory in a compact riemannian manifold.
\newblock {\em HAL-01158323}, 2015.

\bibitem{BerGauMaz71}
M.~Berger, P.~Gauduchon, and E.~Mazet.
\newblock {\em Le spectre d'une vari\'et\'e riemannienne}, volume 194 of {\em
  Lecture Notes in Mathematics}.
\newblock Springer-Verlag, 1971.

\bibitem{BerghLofstrom76}
J.~Bergh and J.~L\"ofStr\"om.
\newblock {\em Interpolation Spaces}, volume 223.
\newblock Grundlehren der mathematischen Wissenchaften, 1976.

\bibitem{Bueler99}
E.~L. Bueler.
\newblock The heat kernel weigted {H}odge laplacian on non compact manifolds.
\newblock {\em T.A.M.S.}, 351(2):683--713, 1999.

\bibitem{Donnelly81}
H.~Donnelly.
\newblock The diffential form spectrum of hyperbolic spaces.
\newblock {\em Manuscripta Math.}, 33:365--385, 1981.

\bibitem{EvGar92}
L.~C. Evans and R.~F. Gariepy.
\newblock {\em Measure theory and fine properties of functions.}
\newblock Studies in Advanced Mathematics. CRC Press, Boca Raton, 1992.

\bibitem{Gaffney55}
M.~P. Gaffney.
\newblock Hilbert space methods in the theory of harmonic integrals.
\newblock {\em Amer. Math. Soc.}, 78:426--444, 1955.

\bibitem{GuilbargTrudinger98}
D.~Gilbarg and N.~Trudinger.
\newblock {\em Elliptic Partial Differential equations}, volume 224 of {\em
  Grundlheren der mathematischen Wissenschaften}.
\newblock Springer, 1998.

\bibitem{Gromov91}
M.~Gromov.
\newblock K\"ahler hyperbolicity and ${L^2}$-{H}odge theory.
\newblock {\em J. Differential Geom.}, 33:253--320, 1991.

\bibitem{GuneysuPigola}
B.~Guneysu and S.~Pigola.
\newblock Calderon-{Z}ygmund inequality and {S}obolev spaces on noncompact
  riemannian manifolds.
\newblock {\em Advances in Mathematics}, 281:353--393, 2015.

\bibitem{Hebey96}
E.~Hebey.
\newblock {\em Sobolev spaces on Riemannian manifolds.}, volume 1635 of {\em
  Lecture Notes in Mathematics}.
\newblock Springer-Verlag, Berlin, 1996.

\bibitem{HebeyHerzlich97}
E.~Hebey and M.~Herzlich.
\newblock Harmonic coordinates, harmonic radius and convergence of riemannian
  manifolds.
\newblock {\em Rend. Mat. Appl. (7) 17 (1997), no. 4, 569-605 (1998)},
  17(4):569--605, 1997.

\bibitem{Kodaira49}
K.~Kodaira.
\newblock Harmonic fields in riemannian manifolds (generalized potential
  theory).
\newblock {\em Annals of Mathematics}, 50:587--665, 1949.

\bibitem{XDLi09}
X.-D. Li.
\newblock On the strong ${L}^p$-{H}odge decomposition over complete riemannian
  manifolds.
\newblock {\em Journal of Functional Analysis}, 257:3617--3646, 2009.

\bibitem{XDLi2010}
X.D. Li.
\newblock ${L}_p$-estimates and existence theorems for the $\bar\partial$
  -operator on complete k\"ahler manifolds.
\newblock {\em Advances in Mathematics}, 224:620--647, 2010.

\bibitem{XDLi010}
X.D. Li.
\newblock Riesz transforms on forms and ${L}_p$-hodge decomposition on complete
  riemannian manifolds.
\newblock {\em Rev. Mat. Iberoamericana}, 26(2):481--528, 2010.

\bibitem{XDLi10}
X.D. Li.
\newblock Sobolev inequalities on forms and ${L}_{p,q}$-cohomology on complete
  riemannian manifolds.
\newblock {\em J. Geom. Anal.}, 20:354--387, 2010.

\bibitem{Lohoue16}
N.~Lohou\'e.
\newblock L'\'equation de {P}oisson pour les formes diff\'erentielles sur un
  espace sym\'etrique et ses applications.
\newblock {\em Bulletin des sciences math\'ematiques}, 140:11--57, 2016.

\bibitem{Lott96}
J.~Lott.
\newblock The zero-in-the-spectrum question.
\newblock {\em Enseignement Math\'ematiques}, 42:341--376, 1996.

\bibitem{Morrey66}
C.~B. Morrey.
\newblock {\em Multiple Integrals in the Calculus of Variations}, volume 130 of
  {\em Die Grundlehren der mathematischen Wissenschaften}.
\newblock Springer-Verlag Berlin Heidelberg New York, 1966.

\bibitem{Scott95}
C.~Scott.
\newblock ${L}^p$ theory of differential forms on manifolds.
\newblock {\em Transactions of the Americain Mathematical Society},
  347(6):2075--2096, 1995.

\bibitem{TaylorGD}
M.~E. Taylor.
\newblock {\em Differential Geometry}.
\newblock Course of the University of North Carolina. University of North
  Carolina.
\newblock www.unc.edu/math/ Faculty/met/diffg.html.

\bibitem{Varopoulos89}
N.~Varopoulos.
\newblock Small time gaussian estimates of heat diffusion kernels.
\newblock {\em Bulletin des Sciences Math\'ematiques.}, 113:253--277, 1989.

\bibitem{Witten82}
E.~Witten.
\newblock Supersymmetry and morse theory.
\newblock {\em J. Differential Geometry}, 17:661--692, 1982.

\end{thebibliography}

\end{document}